\documentclass[aos]{imsart}

\RequirePackage{amsthm,amsmath,amsfonts,amssymb}
\RequirePackage[numbers]{natbib}
\RequirePackage[colorlinks,citecolor=blue,urlcolor=blue]{hyperref}
\RequirePackage{graphicx}

\usepackage[T1]{fontenc}
\usepackage[utf8]{inputenc}

\usepackage{fancyhdr}
\usepackage{bigints}
\usepackage{amsthm,amsmath,amssymb,bm,enumerate,mathrsfs,mathtools}
\newcommand\numberthis{\addtocounter{equation}{1}\tag{\theequation}}
\usepackage{latexsym,color,verbatim,multirow}
\usepackage{graphicx}
\usepackage{caption}
\usepackage{subcaption}
\usepackage{mathtools, enumerate}

\usepackage{graphics,amssymb,color}
\usepackage{graphicx}
\usepackage{algorithm}
\usepackage{algorithm,algorithmic, enumitem, kantlipsum}
\usepackage{xfrac}
\usepackage{subcaption}
\usepackage{float}

\usepackage{tabularx, makecell}
\usepackage{booktabs}
\usepackage[labelfont=bf,format=plain,justification=raggedright,singlelinecheck=false]{caption}

\newcommand{\tE}{\text{Exp}}
\newcommand{\real}{\mathbb{R}}

\newtheorem{proposition}{Proposition}
\newtheorem{theorem}{Theorem}
\newtheorem{lemma}{Lemma}
\newtheorem{corollary}{Corollary}
\newtheorem{remark}{Remark}
\newtheorem{assumption}{Assumption}

\newtheorem{example}{Example}[section]

\newcommand{\rl}{\mathrm{l}}

\newcommand{\rB}{\mathrm{B}}
\newcommand{\rK}{\mathrm{K}}
\newcommand{\rR}{\mathrm{R}}
\newcommand{\rg}{\mathrm{g}}
\newcommand{\rP}{\mathrm{P}}
\newcommand{\rE}{\mathrm{E}}
\newcommand{\rL}{\mathrm{L}}
\newcommand{\rH}{\mathrm{H}}
\newcommand{\rF}{\mathrm{F}}
\newcommand{\rG}{\mathrm{G}}

\newcommand{\rd}{\mathrm{d}}
\newcommand{\rD}{\mathrm{D}}
\newcommand{\rp}{\mathrm{p}}

\newcommand{\cS}{\mathcal{S}}
\newcommand{\cD}{\mathcal{D}}

\newcommand{\oE}{E_{\text{obs}}}
\newcommand{\oA}{A_{\text{obs}}}
\newcommand{\oU}{U_{n; \text{obs}}^{\;(j)}}
\newcommand{\oV}{V_{n; \text{obs}}^{\;(j)}}
\newcommand{\oT}{T_{n; \text{obs}}^{\;(j)}}
\newcommand{\oZ}{\mathcal{Z}_{n;\text{obs}}}

\allowdisplaybreaks
\makeatletter
\newcommand*{\rom}[1]{\expandafter\@slowromancap\romannumeral #1@}
\makeatother

\newcommand*{\Scale}[2][4]{\scalebox{#1}{$#2$}}

\makeatletter
\renewcommand*\env@matrix[1][\arraystretch]{%
  \edef\arraystretch{#1}%
  \hskip -\arraycolsep
  \let\@ifnextchar\new@ifnextchar
  \array{*\c@MaxMatrixCols c}}
\makeatother

\usepackage{makecell}

\startlocaldefs
\numberwithin{equation}{section}
\theoremstyle{plain}

\endlocaldefs

\def\widebreve#1{\mathop{\vbox{\m@th\ialign{##\crcr\noalign{\kern\p@}%
  \brevefill\crcr\noalign{\kern0.1\p@\nointerlineskip}%
  $\hfil\displaystyle{#1}\hfil$\crcr}}}\limits}

\makeatletter
\newcommand{\myitem}[1]{%
\item[#1]\protected@edef\@currentlabel{#1}%
}
\makeatother

\fontfamily{cmr}
\begin{document}

\begin{frontmatter}
\title{Carving model-free inference}

\begin{aug}
\author{\fnms{Snigdha} \snm{Panigrahi} \ead[label=e1]{psnigdha@umich.edu}}
\affiliation{Department of Statistics\\
University of Michigan}
%\address{
%451 West Hall, 1085 South University,\\
%Department of Statistics,\\
%University of Michigan.\\
%Ann Arbor, MI 48109-1107.\\
%\printead{e1}
%}
\end{aug}

\begin{abstract}
Complex studies involve many steps.
Selecting promising findings based on pilot data is a first step.
As more observations are collected, the investigator must decide how to combine the new data with the pilot data to construct valid selective inference.
Carving, introduced by \cite{optimal_inference}, enables the reuse of pilot data during selective inference and accounts for over-optimism from the selection process.
Currently, the justification for carving is tied to parametric models, like the commonly used Gaussian model. 
In this paper, we develop the asymptotic theory to substantiate the use of carving beyond Gaussian models.
Through both simulated and real instances, we find that carving produces valid and tight confidence intervals within a model-free setting. 
\end{abstract}

\end{frontmatter}

\section{Introduction}
\label{introduction}
Inference for a selected set of findings, also called selective inference, is a common problem in complex studies. 
Usually, the investigator begins with pilot data to select a set of promising findings.
As more observations are collected later, the investigator faces the question of how to augment the new data with the existing pilot data for valid selective inference.
On the one hand, a direct augmentation of the two datasets ignores over-optimism from the selection process.
For example, a recent article by \citep{benjamini2020selective} highlights replicability concerns with standard inference that does not account for the selection process.
On the other hand, valid selective inference, which relies only on the new data, fails to utilize observations from the pilot data. 
The latter practice is popularly known as data splitting.

\textit{Carving}, introduced by \cite{optimal_inference}, is an efficient alternative to data splitting.
Carving permits the reuse of pilot data by basing valid selective inference on a conditional distribution. 
This distribution accounts for over-optimism from the selection process by conditioning on the selection outcome seen in the pilot data.
Previous work by \cite{exact_lasso, tibshirani2016exact, randomized_response} give a recipe to construct pivots from such conditional distributions.
Applying the same recipe yields us a pivot for carving, which we call a \textit{carved pivot} in our paper.

When data is generated by a Gaussian model, the carved pivot provides exactly-valid selective inference.
But what happens when we drift away from Gaussian data? 
In model-free settings, is it still possible to use the carved pivot for asymptotically-valid selective inference?
Can we trust selective inference even when we observe rare selection outcomes in our pilot data?
This paper demonstrates that a carved pivot produces asymptotically-valid selective inference, even if our data is not from a Gaussian model.
According to our theory, selective inference with a carved pivot remains valid for rare selection outcomes.
 
\subsection{Notation}
We list some basic notations for our paper.
For $d\in \mathbb{N}$, let $[d]= \{1,2,\cdots, d\}$.
Let $|E|$ be the cardinality of set $E$ and let $E^c$ be its complement set.
Let $V^{\;(j)}$ be the $j^{\text{th}}$ component of the vector $V\in \mathbb{R}^d$.
Let $V^{\;(-j)}$ be the subvector of $V$ after we exclude the $j^{\text{th}}$ component of the original vector and let $V^{\;(E)}$ be the subvector of $V$ that collects the components in $E\subseteq [d]$.
The symbol $V'$ denotes the transpose of the vector $V$.
Unless mentioned otherwise, $\| V\|$ is understood as the $\ell_2$-norm of $V$. 
%For two vectors $U,V$ with matching dimensions, $U\geq V$ implies that the inequality holds in a componentwise sense.
$\mathcal{P}_E$ denotes a permutation matrix: $\mathcal{P}_E V$ re-orders the components of $V$ and returns the vector $\begin{pmatrix} {V^{(E)}}' & {V^{(E^c)}}'\end{pmatrix}'$.
For a positive definite matrix $M$, let $M^{1/2}$ be its principal square root.
For any matrix $M\in \mathbb{R}^{d_1\times d_2}$, $E_1\subseteq [d_1]$, $E_2\subseteq [d_2]$, let $M_{E_1, E_2}$ be the submatrix of $M$ that contains rows and columns in the sets $E_1$ and $E_2$, respectively.
Also, let $M_{E_1}$ be the submatrix of $M$ which collects its columns in the set $E_1$.
We use the notations $I_{d,d}$ and $0_{d_1,d_2}$ for the identity matrix with $d$ rows and columns and the matrix of all zeros with $d_1$ rows and $d_2$ columns, respectively.
We use the notations $0_{d}$ and $1_{d}$ to denote a vector with all $d$ components equal to zero and a vector with all $d$ components equal to one, respectively.
%Th collection $\{e_j \in \real^d: \; j\in [d]\}$ represents the standard basis vectors of $\real^d$.
%Consider $f:\mathcal{C}\to \real$ for an open set $\mathcal{C}\subset \real^d$.
%Then, $\mathcal{D}^k f(x_0)[i_1, i_2,\cdots, i_k]$ is the $k^{\text{th}}$ order partial derivative of $f$ with respect to $x^{(i_1)}$, $\cdots$, $x^{(i_k)}$ at $x_0\in \mathcal{C}$. 
%Specifically, if $d=1$, then $\mathcal{D}^{k}f(x_0)$ denotes the $k^{\text{th}}$ derivative of $f$ at $x_0$.
%We use $\sup f$ for the supremum of a bounded, real-valued function $f$.
For a positive semidefinite matrix $\Sigma \in \real^{d\times d}$ and $x\in \real^d$, let $\tE(x, \Sigma)=\exp(-\frac{1}{2}x' \Sigma x)$.
Let $\mathbf{1}_{\mathcal{E}}$ represent the indicator function, where $\mathcal{E}$ is a fixed set.
At last, let $\Phi$ be the CDF of a standard Gaussian distribution with the density function $\phi$, and let $\bar{\Phi}(x)= 1-\Phi(x)$ be the corresponding survival function at $x$. 
%In our paper, the symbol $Z$ is understood to be a standard Gaussian random vector, unless indicated otherwise.
% on $\real^d$, i.e., $Z\sim N(0_d, I_{d,d})$.

\subsection{Organization}
In Section \ref{framework}, we present a carved pivot for exactly-valid selective inference with Gaussian data.
We introduce a running example for our paper in this section, which allows us to develop the main ideas behind the asymptotic theory. 
In Section \ref{first:step:prelims}, we show that asymptotically-valid selective inference depends on the convergence of specific relative differences. 
In Section \ref{weak:convergence:Gaussian}, we build the asymptotic theory for $\real^d$-valued random variables with the identity covariance matrix.
We then extend the asymptotic theory to variables with a general covariance matrix in Section \ref{weak:convergence:Gaussian:multi}.
In Section \ref{empirical:results}, we study the empirical behavior of the carved pivot on synthetic and real data.
Section \ref{conclusion} concludes our paper with a brief discussion.
Proofs and supporting results are collected in the Appendix.

\section{Exactly-valid selective inference with carving}
\label{framework}

\subsection{Our running example}
Suppose that we observe a triangular array of independent and identically distributed $\mathbb{R}^d$-valued observations
\begin{equation}
\label{generative:P}
\zeta_{i,n}= \begin{pmatrix} \zeta^{\;(1)}_{i,n} & \zeta^{\;(2)}_{i,n} & \cdots &  \zeta^{\;(d)}_{i,n} \end{pmatrix}'\; \stackrel{\text{i.i.d.}}{\sim} \mathbb{P}_n,\ \text{ for }  i\in [n].
\end{equation}
Let 
$$\beta_n= \mathbb{E}_{\mathbb{P}_n}\left[\zeta_{1,n}\right] \in \real^d$$ 
be the unknown mean parameter. 
Let
$$\Sigma=\mathbb{E}_{\mathbb{P}_n}\Big[ (\zeta_{1,n}-\beta_n)(\zeta_{1,n}-\beta_n)' \Big]$$
be the $d\times d$ covariance matrix, which we assume is fixed and invertible.
Define
$$V_n= \sqrt{n}\bar{\zeta}_{n},$$
where $\bar{\zeta}_{n}=\frac{1}{n}\displaystyle\sum_{i=1}^n\zeta_{i,n}$.
Additionally, let 
$$\Sigma_{-j,j}= \text{Cov}_{\mathbb{P}_n}\left(V_n^{\;(-j)}, V_n^{\;(j)}\right), \ \sigma_j^2 = \text{Var}_{\mathbb{P}_n}\left(V_n^{\;(j)}\right),\; \text{ for }\; j\in [d].$$ 
Throughtout, we will assume that the distribution of $V_n$ admits a Lebesgue density.

For a fixed constant $\rho\in (0,1)$, we consider a Gaussian variable
$$W_n\sim \mathcal{N}(0_d, \rho^2\Sigma),$$
which is independent of $V_n$.
Then, using $V_n$ and $W_n$, we infer for $\beta^{\;(j)}_n$, the $j^{\text{th}}$ component of $\beta_n$, only if
\begin{equation}
\label{selection}
V_{n}^{\;(j)} + W_n^{\;(j)} >\Lambda^{\;(j)},
\end{equation}
where $\Lambda$ is a fixed vector in $\real^d$.
%That is, a perturbed version of the related Z-statistic must exceed $\Lambda$ for the selected means.
Borrowing the term \textit{randomization} from \citep{randomized_response}, we call $W_n$ a \textit{randomization variable}.
The rule used for selection, in \eqref{selection}, is called a \textit{randomized selection} rule. 
As shown afterwards, there is an asymptotic correspondence between \eqref{selection} and a similar selection on pilot data. 
We use the symbol $E_n$ to represent the indices of our selected means. 
Let
$$\oE\subseteq [d]$$
be the observed value of the random variable $E_n$. 
For brevity sake, let $|\oE|=p$.

\subsection{Exactly-valid selective inference with Gaussian data}
Suppose that our data is drawn from
 $$\mathbb{P}_n = \mathcal{N}(\beta_n, \Sigma).$$
In this case, $V_n$ is distributed as Gaussian variable with mean vector $\sqrt{n}\beta_n$ and covariance matrix $\Sigma$. 
Assume that $\oE\neq \emptyset$.
Consider $j\in \oE$. 
We obtain a conditional distribution for $V^{\;(j)}_n$, which accounts for the selection process by conditioning on a subset of the selection outcome
$$ \{E_n= \oE\}.$$
Below, we review a pivot for $\beta^{\;(j)}_n$ based on this conditional distribution.

First, we introduce some more statistics. 
Define
$$A_n= V_n^{\;(E_n^c)} + W_n^{\;(E_n^c)},\; U_n^{\;(j)} = V_n^{\;(-j)}- \frac{1}{\sigma_j^2}\Sigma_{-j,j} V_n^{\;(j)}.$$
Let $\oA$ be the observed value of $A_n$.
To draw valid selective inference, we consider the distribution of $V_n^{\;(j)}$ when conditioned on 
\begin{equation}
\label{cond:event}
\Big\{E_n= \oE, \; A_n= \oA\Big\}
\end{equation}
and the observed value of $U_n^{\;(j)}$.
Note that we condition on a subset of the selection outcome by further conditioning on some additional information $A_n$.
By adding extra conditioning, the conditional distribution becomes simpler since the selection outcome can be described as a set of straightforward sign constraints.
Additionally, we condition on $U_n^{\;(j)}$ to eliminate all parameters except $\beta^{\;(j)}_n$.
The resulting conditional distribution is then used to construct a pivot for $\beta^{\;(j)}_n$.
Proposition \ref{carved:pivot} introduces this pivot for Gaussian data.
 
Before presenting our pivot, we fix the matrices
$$R^{\;(j)}= \mathcal{P}_{\oE}\begin{bmatrix} 1 & 0 \\ \frac{1}{\sigma_j^2}\Sigma_{-j,j} & I_{d-1,d-1} \end{bmatrix},\;\; Q = \begin{bmatrix} I_{p, p} \\ 0_{d-p, p} \end{bmatrix}, \;\; r = \begin{pmatrix} \Lambda^{(\oE)} \\  \oA\end{pmatrix}.$$
Define $\rF:\real^d \to \real$ as
 \begin{align*}
& \rF(V) = \int  \tE\left(Qt -V + r, \frac{1}{\rho^2}\Sigma^{-1}\right) \cdot \mathbf{1}_{t\in \real^{p+}} dt,
 \end{align*} 
 and define
$$\rD(U; \sqrt{n}\beta_n^{\;(j)})= \int_{-\infty}^{\infty}\phi\left(\dfrac{1}{\sigma_j}(v- \sqrt{n}\beta_n^{\;(j)})\right) \cdot \rF\left(R^{\;(j)}\begin{pmatrix} v & U' \end{pmatrix}'\right) dv.$$

\begin{proposition}[Pivot]
\label{carved:pivot}
Let $\text{\normalfont Pivot}^{\;(j)}\left(\begin{pmatrix}V_n^{\;(j)} & (U_n^{\;(j)})'\end{pmatrix}'\right)$ be equal to
$$(\rD(U_n^{\;(j)}; \sqrt{n}\beta_n^{\;(j)}))^{-1}\cdot \int_{V_n^{\;(j)}}^{\infty}\phi\left(\dfrac{1}{\sigma_j}(v- \sqrt{n}\beta_n^{\;(j)})\right)\cdot \rF\left(R^{\;(j)}\begin{pmatrix} v & (U_n^{\;(j)})' \end{pmatrix}'\right) dv.$$
Conditional on the outcome in \eqref{cond:event},  $\text{\normalfont Pivot}^{\;(j)}\left(\begin{pmatrix}V_n^{\;(j)} & (U_n^{\;(j)})'\end{pmatrix}'\right)$ is distributed as a $\text{\normalfont Unif}\;(0,1)$ variable.
\end{proposition}

The pivot in Proposition \ref{carved:pivot} applies broadly to several instances of selective inference.
We provide more examples in Section \ref{empirical:results}, which includes inference after variable selection.
In each instance, we construct a pivot with a similar representation as Proposition \ref{carved:pivot}.

Before proceeding further, we turn to a special case when $\Sigma=I_{d,d}$.
We obtain a reduced form for our pivot in this special case.
To simplify further, we fix $\Lambda= 0_d$.
Notice that the components of $\beta_n$ do not share a relationship with one another.
As a result, the pivot turns out to be a univariate function in the statistic $V_n^{\;(j)}$, i.e., it is free of $U_n^{\;(j)}$.
%In particular, the pivot is derived from the Gaussian distribution of $V_n^{\;(j)}$ when conditioned on the selection outcome $\{j\in \oE\}$.

\begin{corollary}[Univariate pivot]
\label{carved:pivot:sequence}
The pivot for $\beta_n^{\;(j)}$ in Proposition \ref{carved:pivot} simplifies as
\[\text{\normalfont Pivot}^{\;(j)}\left(V_n^{\;(j)}\right) = \dfrac{\displaystyle\int_{V_n^{\;(j)}-\sqrt{n}\beta^{\;(j)}_n}^{\infty} \phi(v)\cdot \bar{\Phi}\left(-\frac{1}{\rho}{(v+\sqrt{n}\beta^{\;(j)}_n)}\right)\;dv}{\displaystyle\int_{-\infty}^{\infty}\phi(v)\cdot \bar{\Phi}\left(-\frac{1}{\rho}{(v+\sqrt{n}\beta^{\;(j)}_n)}\right)\;dv }.\]
\end{corollary}

\subsection{Contributions and related developments}
Consider an array of $\real^d$-valued observations from $\mathbb{P}_n$ as described earlier.
Suppose now, a similar selection rule is applied only to a random subsample of size $n_1 (<n)$.
These $n_1$ randomly chosen samples play the role of pilot data in our setup.
Let $n_2 = n-n_1$ and let
$$\rho^2=\frac{n_2}{n_1}$$
be the ratio of the number of observations in the new data to the number of observations in the pilot data.
We infer for $\beta^{\;(j)}_n$ only if the corresponding $Z$-test statistic exceeds a threshold $\tau^{\;(j)}$ in the pilot data, i.e., 
\[V^{\;(j)}_{n_1}> \tau^{\;(j)}.\]
%which is indicative of a positive effect for the variable in the pilot data.

The selection rule on the $n_1$ subsamples can be transformed into the randomized selection rule in \eqref{selection}, in an asymptotic sense.
To see this, we define
\begin{equation}
\label{implicit:randomization:FE}
%W_n^{\;(j)} = \frac{\sqrt{n}}{\sqrt{n_1}}V^{\;(j)}_{n_1} - V^{\;(j)}_{n}, \; \text{ for } j \in [d],
W_n^{\;(j)} = {\sqrt{1+ \rho^2}}\cdot V^{\;(j)}_{n_1} - V^{\;(j)}_{n}, \; \text{ for } j \in [d],
\end{equation}
and also let
$$\Lambda^{\;(j)} = {\sqrt{1+ \rho^2}}\cdot{\tau^{\;(j)}}.$$
Asymptotically, $W_n$ is distributed as a Gaussian variable with mean $0_d$ and covariance $\rho^2\Sigma$ and is independent of $V_n$. 

Specifically, when $\mathbb{P}_n = \mathcal{N}(\beta_n, \Sigma)$, we easily verify that 
$$W_n\sim \mathcal{N}(0_d, \rho^2\Sigma),$$
and that $W_n$ is independent of $V_n$.
That is, for Gaussian data, the selection on $n_1$ subsamples coincides exactly with the randomized selection in \eqref{selection}.
This example was provided in \cite{selective_bayesian} for $d=1$. 
In this situation, we note that Proposition \ref{carved:pivot} re-uses the pilot data to produce a carved pivot for exactly-valid selective inference.

\begin{figure}[h]
  \centering
    \hspace*{-0.5cm}
  \begin{tabular}{@{}c@{}}
    \includegraphics[height=5cm, width=11.1cm]{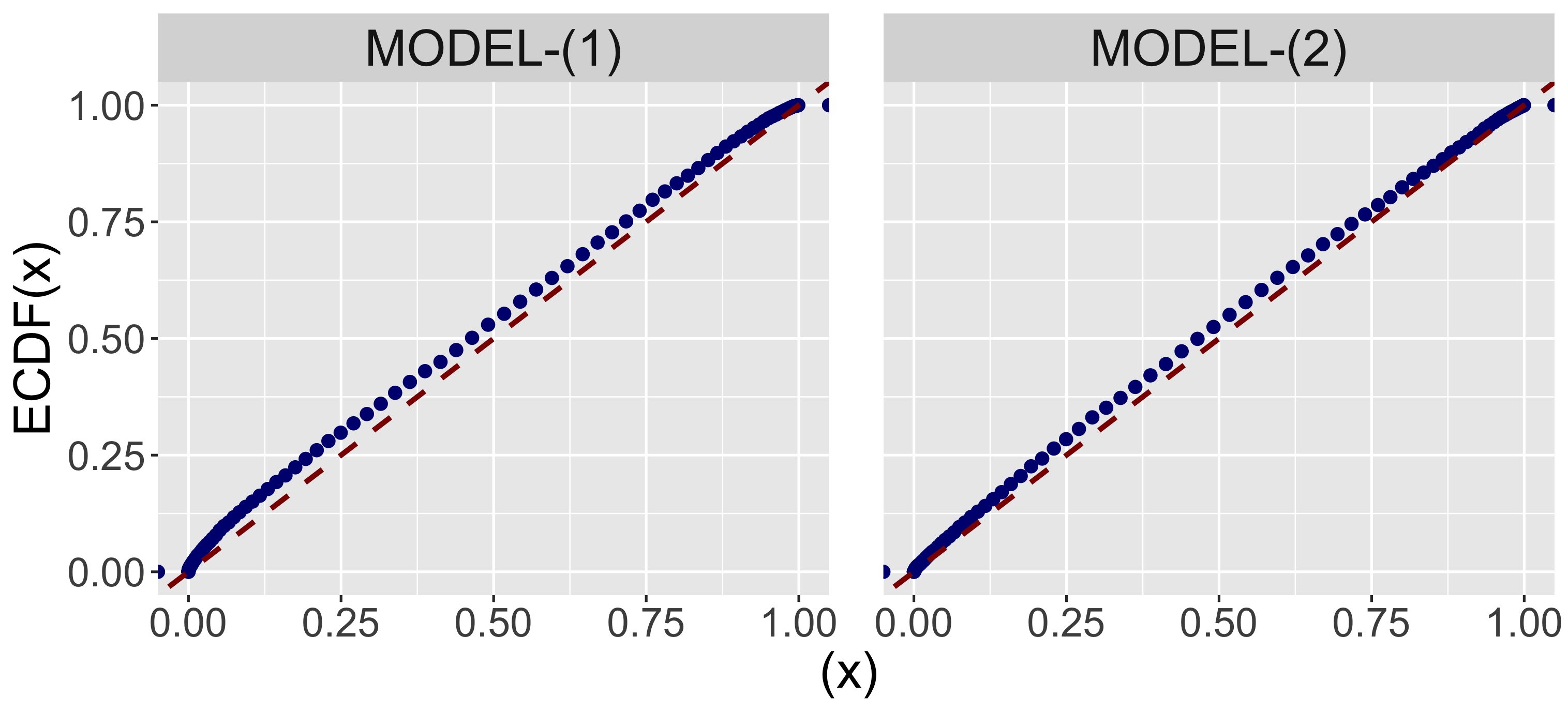} \\
     \end{tabular}
  \vspace{\floatsep}
  \hspace*{-0.2cm}
  \begin{tabular}{@{}c@{}}
    \includegraphics[height=5cm, width=11.1cm]{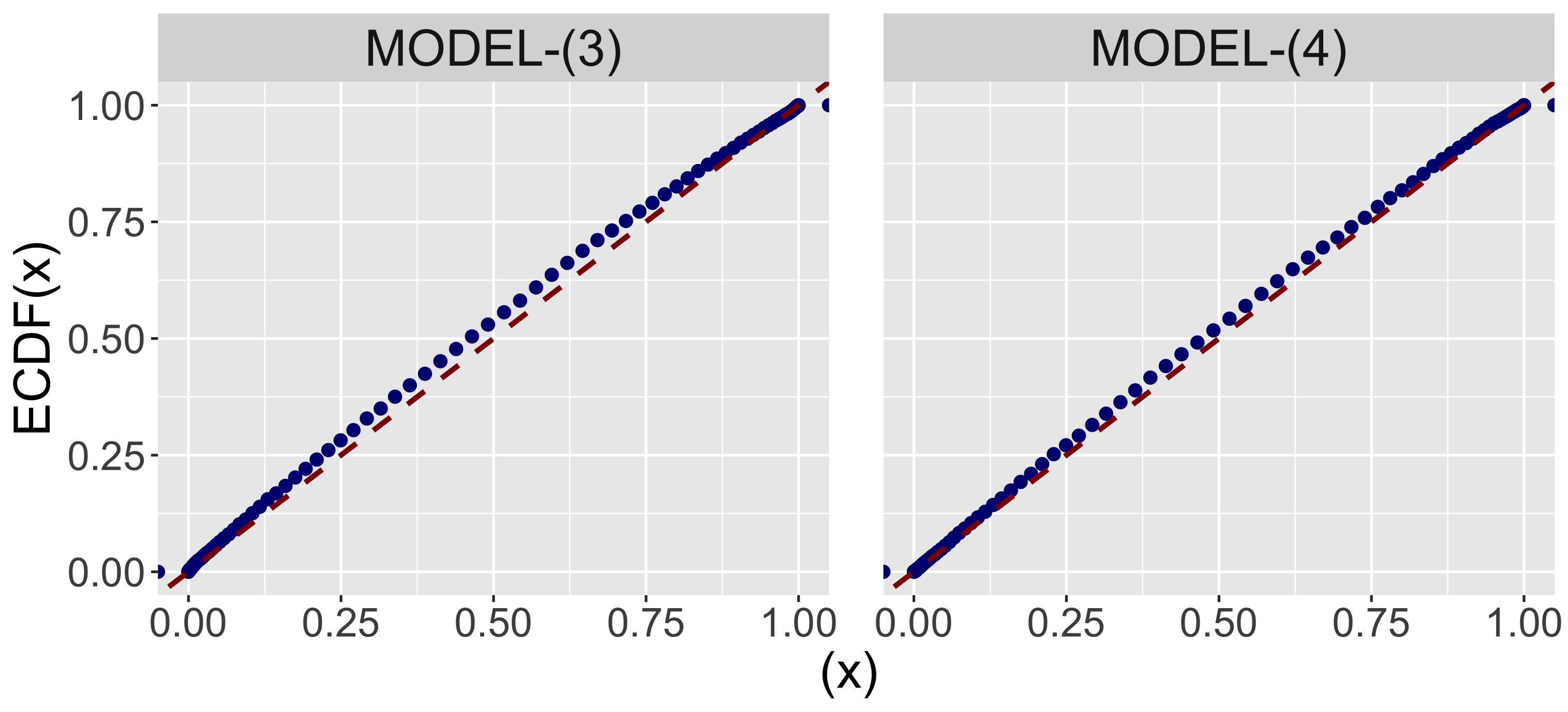} \\
      \end{tabular}
   \vspace{-8mm}
  \caption{The four panels plot the ECDF of the carved pivot when data is generated according to MODELS (1)-(4). The red dashed line represents the reference $y=x$ curve.}
  \label{fig:first}
\end{figure}

What happens when we drift away from Gaussian data?
We begin with a simple simulation.
We draw our data from four different models with non-Gaussian errors and conduct $10000$ rounds of simulation from each model.
Figure \ref{fig:first} depicts the empirical cumulative distribution function (ECDF) of the carved pivot.
A strong alignment with the $y=x$ line indicates that the carved pivot is well approximated by a $\text{Unif} \;(0,1)$ variable.

Let $n_1 = n_2=25$, i.e., $\rho^2=1$. 
Fix $\Sigma=I_{d, d}$ and fix $\tau = 0_{d}$.
Let
\begin{equation*}
\sqrt{n}\beta^{\;(j)}_n = -a_n \bar{\beta}, \text{ for} \;j \in [d],
%\label{mean:eg}
\end{equation*} 
where $a_n = n^{1/6-\delta}$, $\delta=1\mathrm{e}{-3}$ and $\bar{\beta}= 1.5$. 
We draw 
$$\zeta_{i,n} = \beta_n + \mathrm{e}_{i,n}.$$
First, each component of the error vector $\mathrm{e}_{i,n}$ is drawn independently from $\rE$.
Then we standardize each such observation to have mean $0$ and variance $1$.
MODELS (1)-(4) are based on four different choices of $\rE$.
\begin{enumerate}
\setlength\itemsep{1em}
\item MODEL-(1)\; $\rE = \text{Exponential}(1)$ with rate parameter equal to $1$ and density equal to $\mathrm{p}(x) = \exp(-x)\cdot \mathbf{1}_{x>0}$.
\item MODEL-(2)\: $\rE= \text{Exponentially Modified Gaussian distribution}(0,1,1)$ with the mean and variance of the Gaussian component equal to $0$ and $1$ respectively, and with the rate parameter of the exponential component equal to $1$, and density equal to  
$$\mathrm{p}(x) = \frac{1}{\sqrt{\pi}}\exp(0.5-x)\int_{\frac{1-x}{\sqrt{2}}}^{\infty} e^{-t^2}dt$$%$\rE = \mathcal{N}(0,1) +  \text{Exponential}(s=1)$. 
\item MODEL-(3)\: $\rE = 0.8\cdot \mathcal{N}(0,0.25) +  0.2\cdot \mathcal{N}(0,3)$, which is a mixture of two Gaussian distributions with mixing weights $0.8$ and $0.2$.
\item MODEL-(4)\: $\rE = \text{Laplace}(0,1)$ with location and scale parameters equal to $0$ and $1$ respectively, and density equal to $\mathrm{p}(x) = (2)^{-1} \exp(- |x|)$.
\end{enumerate}

We make a few interesting observations from Figure \ref{fig:first}.
First, our plot suggests that the carved pivot produces valid selective inference well beyond Gaussian data. 
Previously, \cite{randomized_response} showed that randomized selection rules with heavy tailed variables produced asymptotically-valid selective inference.
In contrast, randomization variables, based on carving, resemble Gaussian variables in the limit. 
%This difference sets apart our contributions from earlier theory.
Second, our plot shows that  selective inference remains valid for rare selection outcomes that have vanishing probabilities in the limit.
%We recognize that the selection outcome is rare, because the (log-) probability for this outcome can be bounded as
%$$\frac{1}{a_n^{2}}\log\mathbb{P}_n\left[V^{\;(j)}_{n_1} >0\right]\leq -\frac{\bar{\beta}^2}{2(1+\rho^2)}.$$
%This means that the selection outcome has a vanishing probability as $n\to \infty$. 
Prior asymptotic work such as those conducted by \cite{markovic2016bootstrap, tibshirani2018uniform, randomized_response} have only focused on selection outcomes with non-vanishing probabilities.
%Our paper investigates if  asymptotic validity of selective inference under rare outcomes has not been investigated yet.
Note that this paper extends the use of carving beyond Gaussian models through our asymptotic theory. 
Additionally, our theory confirms the validity of selective inference for a wide range of rare selection outcomes.
%We close this section with related developments in the field.
%In the past few years, many proposals have been developed to provide valid selective inference. 
%Among these proposals, the conditional approach bypasses some of the obstacles highlighted by early work \cite{leeb2003finite, leeb2005model}, and can be less conservative than a simultaneous approach to the same problem \cite{berk2013valid, bachoc2019valid}.
%As shown for our running example, carving is asymptotically equivalent to selective inference with a randomized selection rule.

Our paper is connected with the fast-growing literature on selective inference with randomization.
In recent work, \cite{kivaranovic2020tight} showed that randomized rules on Gaussian variables yield bounded confidence intervals for selective inference and \cite{panigrahi2018selection} have utilized similar rules to construct confidence intervals for the effects of selected genetic variants.
 \cite{schultheiss2021multicarving} proposed repeated carving for more stable inference in high-dimensional settings. 
 \cite{zrnic2020post} investigated the potential of randomization from an algorithmic stability perspective.
 \cite{panigrahi2022treatment} applied carving to pool summary statistics from prior studies and constructed unbiased estimators for shared parameters.
Work by \cite{panigrahi2018scalable, selective_bayesian} introduced Bayesian methods to construct inference after solving randomized variable selection algorithms.
\cite{rasines2021splitting} utilized a Gaussian randomization variable to split a dataset into two parts.
One part is utilized for selection while the other is kept aside for inference.
The theory in our paper supports the reuse of the first part when moving away from Gaussian data.

\section{Dependence on relative differences}
\label{first:step:prelims}

Our main finding in this section is that asymptotic validity of the carved pivot can be shown to depend on the convergence of specific relative differences. 
We first discuss some preliminaries.

\subsection{Some preliminaries}
We start from the randomized selection rule in \eqref{selection}, where (i) $W_n$ is distributed as a Gaussian random variable
\begin{equation*}
W_n \sim \mathcal{N}(0_d, \rho^2\Sigma), 
\end{equation*}
and (ii) $W_n$ is independent of $V_n$.    
Later, we show that asymptotic guarantees with a Gaussian randomization variable are transferable to carving under some conditions.  
We come back to this topic in Section \ref{weak:convergence:Gaussian:multi}.
 
Fixing some more notations, let 
$$\mathrm{e}_{i,n} = \Sigma^{-1/2}(\zeta_{i,n}-\beta_n),$$ 
and let $Z_{i,n} = \frac{1}{\sqrt{n}}\mathrm{e}_{i,n}$.
We assume that the components of $V_n$ in the set $\oE$ are stacked before the ones in its complement set. 
Hereafter, we find it convenient to work with a standardized version for $V_n$:
$$\mathcal{Z}_n = \Sigma^{-1/2}(V_n - \sqrt{n} \beta_n),$$
which is equal to
%We can write $\mathcal{Z}_n$ as a sum of i.i.d. variables
\begin{equation}
\mathcal{Z}_n =  \sum_{i=1}^{n}\frac{1}{\sqrt{n}}\mathrm{e}_{i,n}= \sum_{i=1}^{n}  Z_{i,n}.
\label{iid:sum}
\end{equation}  
Our pivot, in terms of the standardized variable, is now denoted by
\begin{equation}
\label{std:pivot}
\rP^{\;(j)} \left(\mathcal{Z}_n; \sqrt{n}\beta_n\right)=\text{\normalfont Pivot}^{\;(j)}\left( (R^{\;(j)})^{-1}\left(\Sigma^{1/2} \mathcal{Z}_n+ \sqrt{n}\beta_n\right)\right).
\end{equation}
This is based on noting the following equality
$$\begin{pmatrix} {V_n^{\;(j)}} & {U_n^{\;(j)}}' \end{pmatrix}' = (R^{\;(j)})^{-1}\left(\Sigma^{1/2} \mathcal{Z}_n+ \sqrt{n}\beta_n\right).$$

Our next result computes the ratio between the conditional and unconditional likelihood functions, after and before we apply the randomized selection rule.
We use the symbol
$${\text{\normalfont LR}}_{\mathbb{P}_n}(\mathcal{Z}_n; \sqrt{n}\beta_n)$$
to denote this ratio at $\beta_n$.
\begin{remark}
We stress that ${\text{\normalfont LR}}_{\mathbb{P}_n}(\mathcal{Z}_n; \sqrt{n}\beta_n)$ represents how the selection process affects the unconditional likelihood.
It is worth noting that this ratio is different from a ratio of the same likelihood function at two distinct values of $\beta_n$. 
\end{remark}
 
In the rest of the paper, we use $\mathbb{E}_{\mathbb{P}_n}[Z]$ to denote the expectation of the standardized variable $Z$ when based on the distribution $\mathbb{P}_n$. 
We use the more specific symbol $\mathbb{E}_{\mathcal{N}}[Z]$ to represent the expectation of $Z$ when $Z\sim \mathcal{N}(0_d, I_{d,d})$.
\begin{proposition}[Ratio of likelihood functions after and before selection]
\label{lik:ratio}
Let $\rF:\real^d \to \real$ be defined according to Proposition \ref{carved:pivot}. Then, the ratio between the conditional likelihood and unconditional likelihood functions is equal to
$${\text{\normalfont LR}}_{\mathbb{P}_n}(\oZ; \sqrt{n}\beta_n)=  \dfrac{\rF \left( \Sigma^{1/2}\oZ + \sqrt{n}\beta_n\right)}{{\mathbb{E}_{\mathbb{P}_n} \left[\rF\left(\Sigma^{1/2}\mathcal{Z}_n + \sqrt{n}\beta_n\right)\right]}},$$
where $\oZ$ is the observed value of $\mathcal{Z}_n$.
\end{proposition}

%The denominator of our likelihood ratio represents the probability of selection, and provides an adjustment to standard inference by accounting for the impact of selection. 

\begin{remark}
Suppose that $\mathbb{P}_n = \mathcal{N}(\beta_n, \Sigma)$. 
Equivalently, $\mathcal{Z}_n$ is distributed as $\mathcal{N}(0_d, I_{d,d})$ variable.
In this case, we utilize the subscript $\mathcal{N}$ to indicate that our likelihoods are based on Gaussian data, and the above ratio is denoted by
$${\text{\normalfont LR}}_{\mathcal{N}}(\oZ; \sqrt{n}\beta_n)= \dfrac{\rF \left( \Sigma^{1/2}\oZ + \sqrt{n}\beta_n\right)}{{\mathbb{E}_{\mathcal{N}} \left[\rF\left(\Sigma^{1/2}\mathcal{Z}_n + \sqrt{n}\beta_n\right)\right]}}.$$ 
\end{remark}

Suppose, $\mathcal{Q}$ is a real-valued measurable mapping. 
Through the ratio in Proposition \ref{lik:ratio}, we define
\begin{equation}
\label{conditional:expectation}
\widetilde{\mathbb{E}}_{\mathbb{P}_n}\left[\mathcal{Q}(\mathcal{Z}_n)\right] = \mathbb{E}_{\mathbb{P}_n}\left[\mathcal{Q}(\mathcal{Z}_n)\cdot {\text{LR}}_{\mathbb{P}_n}(\mathcal{Z}_n; \sqrt{n}\beta_n)\right].
\end{equation}
The expectation on the left-hand side of \eqref{conditional:expectation} is taken with respect to the conditional distribution of $\mathcal{Z}_n$, and the expectation on the right-hand side is taken with respect to the unconditional distribution of $\mathcal{Z}_n$. 
%The expectation on the left-hand side of \eqref{conditional:expectation} takes into account the effect of the randomized selection rule. 
Once again, for Gaussian data, we use more specific notations with the subscript $\mathcal{N}$ and define
$$\widetilde{\mathbb{E}}_{\mathcal{N}}\left[\mathcal{Q}(\mathcal{Z}_n)\right] = \mathbb{E}_{\mathcal{N}}\left[\mathcal{Q}(\mathcal{Z}_n)\cdot {\text{LR}}_{\mathcal{N}}(\mathcal{Z}_n; \sqrt{n}\beta_n)\right].$$

We are now ready to formally state our main goal in the paper.
Let $\rH\in \mathbb{C}^3(\real, \real)$ be an arbitrary function with bounded derivatives up to the third order.
Let $\mathcal{C}_n$ be a suitable collection of distributions $\mathbb{P}_n$ that we specify later.
Using our notations, we prove weak convergence of our pivot by showing that
\begin{equation}
\label{sub:qn:validity}
\lim_n \sup_{\mathbb{P}_n\in \mathcal{C}_n}\Big|\widetilde{\mathbb{E}}_{\mathbb{P}_n}\left[\rH\circ\rP^{\;(j)}(\mathcal{Z}_n;  \sqrt{n}\beta_n)\right]- \widetilde{\mathbb{E}}_{\mathcal{N}}\left[\rH\circ\rP^{\;(j)}(\mathcal{Z}_n; \sqrt{n}\beta_n)\right]\Big| = 0.
\end{equation}
The above weak convergence statement indicates that our pivot generates asymptotically-valid conditional inference even as we depart from Gaussian data.
It is important to mention that this statement assures the validity of selective inference across all distributions in the collection $\mathcal{C}_n$.

%The weak convergence guarantee in \eqref{sub:qn:validity} holds uniformly over the collection of distributions $\mathcal{C}_n$, which means the following.
%\begin{remark}
%Denote by $\text{CI}_{q}$ the $100\cdot (1-q)\%$ asymptotically-valid confidence intervals by inverting our pivot.
%Consider an arbitrarily chosen $\epsilon>0$.
%The statement in \eqref{sub:qn:validity} then guarantees the existence of a sample size $n_\epsilon$ such that the false coverage rate of $\text{CI}_{q}$, conditional on the selection outcome, is bounded above by $q+\epsilon$, whenever $n\geq n_\epsilon$.
%The above-stated guarantee holds for all the distributions $\mathbb{P}_n$ in the collection $\mathcal{C}_n$.
%\end{remark}
%A discussion on the distinction between weak convergence in a pointwise and uniform sense can be found in \cite{tibshirani2018uniform}.
%We refer readers to this paper for more on the limitations of weak convergence in a pointwise sense.

\subsection{Relative differences}

Proposition \ref{weak:convergence:relative:multi} recognizes that weak convergence of our pivot depends on the convergence of specific relative differences. 
Before we do so, define
\begin{equation}
\begin{aligned}
\rG_1(\mathcal{Z}_n; \sqrt{n}\beta_n) &= \rF \left( \Sigma^{1/2}\mathcal{Z}_n  + \sqrt{n}\beta_n\right),\\
\rG_2(\mathcal{Z}_n; \sqrt{n}\beta_n) &=\rF \left( \Sigma^{1/2}\mathcal{Z}_n + \sqrt{n}\beta_n\right) \cdot \rH\circ \rP^{\;(j)}(\mathcal{Z}_n\; ;\sqrt{n}\beta_n).
\end{aligned}
 \label{G:funcs}
\end{equation}
Let $\sup f$ denote the supremum of a bounded, real-valued function $f$.
\begin{proposition}[Relative differences]
\label{weak:convergence:relative:multi}
Suppose that
\begin{equation*}
\begin{aligned}
\rR^{(l)}_n &= \left(\mathbb{E}_{\mathcal{N}}\left[\rF \left( \Sigma^{1/2}\mathcal{Z}_n + \sqrt{n}\beta_n\right)\right]\right)^{-1}\cdot \Big|\mathbb{E}_{\mathbb{P}_n}\left[\rG_l(\mathcal{Z}_n; \sqrt{n}\beta_n)\right]- \mathbb{E}_{\mathcal{N}}\left[\rG_l(\mathcal{Z}_n; \sqrt{n}\beta_n)\right]\Big|,
\end{aligned}
\end{equation*}
for $l\in [2]$. 
Let $\sup |\rH| = \rK <\infty$.
Then, it holds that 
$$\Big|\widetilde{\mathbb{E}}_{\mathbb{P}_n}\left[\rH\circ\rP^{\;(j)}(\mathcal{Z}_n;  \sqrt{n}\beta_n)\right]- \widetilde{\mathbb{E}}_{\mathcal{N}}\left[\rH\circ\rP^{\;(j)}(\mathcal{Z}_n;  \sqrt{n}\beta_n)\right]\Big|\leq  \left(\rK\cdot \rR^{(1)}_n + \rR^{(2)}_n\right).$$
\end{proposition}

\begin{remark}
We note that the relative differences $\rR^{(l)}_n$ involve expectations that are computed with respect to the unconditional distribution of $\mathcal{Z}_n$.
\end{remark}
As a result of Proposition \ref{weak:convergence:relative:multi}, the weak convergence statement in \eqref{sub:qn:validity} follows immediately once we show that
$$\lim_n \sup_{\mathbb{P}_n\in \mathcal{C}_n} \rR^{(l)}_n = 0, \text{ for } l\in[2].$$

To close this section, we have a simplified expression for the common denominator in our relative differences.
Define
$$\bar{\Sigma} = (Q'\Sigma^{-1} Q)^{-1}, \;\bar{\mu}_n = \bar{\Sigma} Q' \Sigma^{-1} (\sqrt{n}\beta_n - r),\; \Lambda = \Sigma^{-1}- \Sigma^{-1}Q \bar{\Sigma} Q'\Sigma^{-1}.$$
\begin{proposition}
\label{rep:sel:probability}
We have
\begin{equation*}
\begin{aligned}
\mathbb{E}_{\mathcal{N}} \left[\rF\left(\Sigma^{1/2}\mathcal{Z}_n + \sqrt{n}\beta_n\right)\right]= C_0\cdot \text{\normalfont Exp}\left(\sqrt{n}\beta_n-r, \frac{1}{(1+\rho^2)}\cdot \Lambda\right) \cdot\mathbb{P}_{\mathcal{N}}\left[ T_n >0_{p}\right],
\end{aligned}
\end{equation*}
where $T_n \sim \mathcal{N}\left(\bar{\mu}_n, (1+\rho^2)\bar{\Sigma} \right)$ and $C_0$ is a constant which does not depend on $n$.
\end{proposition}
We note that $\mathbb{P}_{\mathcal{N}}\left[ T_n >0_{p}\right]$ is the probability of our selection outcome when $\mathbb{P}_n =\mathcal{N}(\beta_n, \Sigma)$.
Put another way, Proposition \ref{rep:sel:probability} states how the common denominator of our relative differences depends on this probability.
%Observe that $\mathbb{P}_{\mathcal{N}}\left[ T_n >0_{p}\right]$ represents the probability of our selection outcome when we draw Gaussian data.
%Put another way, Proposition \ref{rep:sel:probability} states how the common denominator of our relative differences depends on the probability of the selection. 

\subsection{Revisiting the univariate pivot}

We revisit our univariate pivot in Corollary \ref{carved:pivot:sequence}.
Recall that $\Sigma=I_{d,d}$ and $\Lambda =0_d$.
Consistent with our earlier notations, we represent the univariate pivot using the standardized variable through
\begin{equation}
\label{std:pivot:uni}
\rP^{\;(j)} \left(\mathcal{Z}_n^{\;(j)}; \sqrt{n}\beta_n^{\;(j)}\right) =\text{\normalfont Pivot}^{\;(j)}\left(\mathcal{Z}_n^{\;(j)}+ \sqrt{n}\beta_n^{\;(j)}\right).
\end{equation}
 
For this special case, we define
\begin{equation}
\begin{aligned}
\widetilde\rG_1(\mathcal{Z}^{\;(j)}_n; \sqrt{n}\beta^{\;(j)}_n) &=\bar{\Phi}\left( -\dfrac{1}{\rho}{(\mathcal{Z}^{\;(j)}_n+\sqrt{n}\beta^{\;(j)}_n)}\right),\\
\widetilde\rG_2(\mathcal{Z}^{\;(j)}_n; \sqrt{n}\beta^{\;(j)}_n) &=\bar{\Phi}\left( -\dfrac{1}{\rho}{(\mathcal{Z}^{\;(j)}_n+\sqrt{n}\beta^{\;(j)}_n)}\right)\cdot  \rH\circ \rP^{\;(j)}(\mathcal{Z}^{\;(j)}_n\; ;\sqrt{n}\beta_n^{\;(j)}).
\end{aligned}
 \label{G:funcs:uni}
\end{equation}
Letting $\widetilde{D}_n = \mathbb{E}_{\mathcal{N}}\left[{\bar{\Phi}\left( -\dfrac{1}{\rho}{(\mathcal{Z}^{\;(j)}_n+\sqrt{n}\beta^{\;(j)}_n)}\right)}\right]$, we now define the relevant relative differences as
\begin{equation}
\begin{aligned}
 & \widetilde\rR^{(l)}_n=  \widetilde{D}_n^{-1}\cdot \Bigg|\mathbb{E}_{\mathbb{P}_n}\left[\widetilde\rG_l(\mathcal{Z}^{\;(j)}_n; \sqrt{n}\beta^{\;(j)}_n)\right] - \mathbb{E}_{\mathcal{N}}\left[\widetilde\rG_l(\mathcal{Z}^{\;(j)}_n; \sqrt{n}\beta^{\;(j)}_n)\right]\Bigg|,
\end{aligned}
\label{carved:pivot:sequence:relative}
\end{equation}
for $l\in [2]$. 
Note that the relative differences are determined solely by the expectations of functions that involve the univariate variable $\mathcal{Z}^{\;(j)}_n$.

Suppose that $\sup |\rH| = \rK <\infty$.
Once again, we can show that
$$\Big|\widetilde{\mathbb{E}}_{\mathbb{P}_n}\left[\rH\circ\rP^{\;(j)} \left(\mathcal{Z}_n^{\;(j)}; \sqrt{n}\beta_n^{\;(j)}\right)\right]- \widetilde{\mathbb{E}}_{\mathcal{N}}\left[\rH\circ\rP^{\;(j)} \left(\mathcal{Z}_n^{\;(j)}; \sqrt{n}\beta_n^{\;(j)}\right)\right]\Big|\leq  \left(\rK\cdot \widetilde\rR^{(1)}_n + \widetilde\rR^{(2)}_n\right).$$
We use two facts here.
First, the pivot is a function of the univariate variable $\mathcal{Z}_n^{\;(j)}$. 
Second, the likelihood ratio, in Proposition \ref{lik:ratio}, is proportional to
\begin{equation*}
 \prod_{j\in \oE} {\text{\normalfont LR}}_{\mathbb{P}_n}(\oZ^{\;(j)}; \sqrt{n}\beta^{\;(j)}_n)
\end{equation*}
where 
 \begin{equation*}
\begin{aligned}
 &\Scale[0.95]{{\text{\normalfont LR}}_{\mathbb{P}_n}(\oZ^{\;(j)}; \sqrt{n}\beta^{\;(j)}_n)= \left\{{\mathbb{E}_{\mathbb{P}_n} \left[ {\bar{\Phi}\left( -\dfrac{1}{\rho}{(\mathcal{Z}_n^{\;(j)}+\sqrt{n}\beta^{\;(j)}_n)}\right)} \right]}\right\}^{-1}  {\bar{\Phi}\left( -\dfrac{1}{\rho}{(\oZ^{\;(j)}+\sqrt{n}\beta^{\;(j)}_n)}\right)}.}
\end{aligned}
\end{equation*}
It is important to note that this ratio depends on $\beta^{\;(j)}_n$ only through the univariate variable $\mathcal{Z}^{\;(j)}_n$.
For Gaussian data, the specific ratio is given by
$$
 {\text{LR}}_{\mathcal{N}}(\oZ^{\;(j)}; \sqrt{n}\beta^{\;(j)}_n)=\dfrac{\bar{\Phi}\left( -\dfrac{1}{\rho}{(\oZ^{\;(j)}+\sqrt{n}\beta^{\;(j)}_n)}\right)}{\mathbb{E}_{\mathcal{N}} \left[ {\bar{\Phi}\left( -\dfrac{1}{\rho}{(\mathcal{Z}^{\;(j)}_n+\sqrt{n}\beta^{\;(j)}_n)}\right)} \right]}.
$$
The steps in the proof of Proposition \ref{weak:convergence:relative:multi} directly lead us to the bound using our relative differences.

%\begin{corollary}[Relative differences under identity covariance]
%\label{carved:pivot:sequence:relative}
%Let $\rR^{(l)}_n$ be equal to
%\begin{equation*}
%\begin{aligned}
%$ & \Scale[0.95]{\left(\mathbb{E}_{\mathcal{N}}\left[{\bar{\Phi}\left( -\dfrac{1}{\rho}{(\mathcal{Z}^{\;(j)}_n+\sqrt{n}\beta^{\;(j)}_n)}\right)}\right]\right)^{-1}\cdot \Bigg|\mathbb{E}_{\mathbb{P}_n}\left[\rG_l(\mathcal{Z}^{\;(j)}_n; \sqrt{n}\beta^{\;(j)}_n)\right] - \mathbb{E}_{\mathcal{N}}\left[\rG_l(\mathcal{Z}^{\;(j)}_n; \sqrt{n}\beta^{\;(j)}_n)\right]\Bigg|}
%\end{aligned}
%\end{equation*}
%for $l\in [2]$. 
%\end{corollary}
At last, we note that the common denominator in our relative differences is equal to
$$\mathbb{E}_{\mathcal{N}} \left[ {\bar{\Phi}\left( -\dfrac{1}{\rho}{(\mathcal{Z}_n^{\;(j)}+\sqrt{n}\beta^{\;(j)}_n)}\right)} \right]=\mathbb{P}_{\mathcal{N}}[j\in \oE]=\bar{\Phi}\left(- \dfrac{\sqrt{n}\beta^{\;(j)}_n}{\sqrt{(1+\rho^2)}} \right),
%=\mathbb{P}_{\mathcal{N}}[j\in \oE],
$$
 which is the probability of the selection outcome on Gaussian data.

 %%%%%%%%%%%%%%%%%%%%%%%%%%%%%%%%%%%%%%%%%%%%%%%%%%%%%%%%%%%%%
 %%%%%%%%%%%%%%%%%%%%%%%%%%%%%%%%%%%%%%%%%%%%%%%%%%%%%%%%%%%%%
\section{Weak convergence of univariate pivot}
\label{weak:convergence:Gaussian}
To better understand the behavior of the multivariate pivot for a general $\Sigma$, we first analyze the simpler univariate pivot.

\subsection{Main results}
In this section, we state our main results in Theorem \ref{weak:convergence:local} and Theorem \ref{weak:convergence:rare}. 
These results demonstrate that our univariate pivot yields asymptotically-valid selective inference for two types of selection outcomes, namely bounded outcomes and rare outcomes.
We describe both types of outcomes below.

%The main results of this section, Theorem \ref{weak:convergence:local} and \ref{weak:convergence:rare}, establish weak convergence of our univariate pivot under two types of selection outcomes, namely bounded outcomes and rare outcomes.

Suppose that the components of our mean vector are bounded, i.e., 
\begin{equation}
  \label{bdd:seq}
  |\sqrt{n}\beta^{\;(j)}_n| <R, \ \text{ for each } n\in \mathbb{N} \text{ and } j\in [d].
\end{equation}
Consider the limiting case when $\mathbb{P}_n= \mathcal{N}(\beta_n, I_{d,d})$.
Recall that the probability of the selection outcome on Gaussian data is equal to:
\begin{equation*}
\mathbb{P}_{\mathcal{N}}[j\in \oE] = \bar{\Phi}\left(- \dfrac{\sqrt{n}\beta^{\;(j)}_n}{\sqrt{(1+\rho^2)}} \right).
%\label{prob:sel:out}
\end{equation*}
Clearly, this probability is bounded away from $0$ whenever the mean satisfies \eqref{bdd:seq}. 
This selection outcome is referred to as a bounded outcome.
%$$\mathbb{E}_{\mathcal{N}} \left[ {\bar{\Phi}\left( -\dfrac{1}{\rho}{(\mathcal{Z}_n^{\;(j)}+\sqrt{n}\beta^{\;(j)}_n)}\right)} \right]=\bar{\Phi}\left(- \dfrac{\sqrt{n}\beta^{\;(j)}_n}{\sqrt{(1+\rho^2)}} \right).$$

From now on, fix $j\in \oE$ and consider the relative differences $\widetilde\rR^{(l)}_n$ in \eqref{carved:pivot:sequence:relative}.

\begin{assumption}
Consider a collection of distributions $\mathcal{P}_{b,n}$ such that the mean of each distribution $\mathbb{P}_n$ in this collection satisfies \eqref{bdd:seq}.
Assume that $\mathcal{P}_{b,n}$ has uniformly bounded third moments in the following sense:
$$\displaystyle\sup_n\displaystyle\sup_{\mathbb{P}_n\in \mathcal{P}_{b,n}} \mathbb{E}_{\mathbb{P}_n}\left[ |\mathrm{e}^{\;(j)}_{1,n}|^3 \right]< \infty,$$
where $\mathrm{e}_{1,n}$ is the standardized variable which was defined in \eqref{iid:sum}.
\label{moment:bdd:uni}
\end{assumption}

\begin{theorem}[Weak convergence of under bounded outcomes]
\label{weak:convergence:local}
Under Assumption \ref{moment:bdd:uni}, we have
$$\lim_n \sup_{\mathbb{P}_n\in \mathcal{P}_{b,n}} \widetilde\rR^{(l)}_n =0\;\; \text{ for } l\in [2],$$
and as a result,
$$\lim_n \sup_{\mathbb{P}_n\in \mathcal{P}_{b,n}}\Big|\widetilde{\mathbb{E}}_{\mathbb{P}_n}\left[\rH\circ\rP^{\;(j)}(\mathcal{Z}^{\;(j)}_n;  \sqrt{n}\beta^{\;(j)}_n)\right]- \widetilde{\mathbb{E}}_{\mathcal{N}}\left[\rH\circ\rP^{\;(j)}(\mathcal{Z}^{\;(j)}_n; \sqrt{n}\beta^{\;(j)}_n)\right]\Big| = 0.$$
\end{theorem}

%Combining the assertion in Theorem \ref{weak:convergence:local} with the bound in Corollary \ref{carved:pivot:sequence:relative} leads us to conclude that
%$$\lim_n \sup_{\mathbb{P}_n\in \mathcal{P}_{b,n}}\Big|\widetilde{\mathbb{E}}_{\mathbb{P}_n}\left[\rH\circ\rP^{\;(j)}(\mathcal{Z}^{\;(j)}_n;  \sqrt{n}\beta^{\;(j)}_n)\right]- \widetilde{\mathbb{E}}_{\mathcal{N}}\left[\rH\circ\rP^{\;(j)}(\mathcal{Z}^{\;(j)}_n; \sqrt{n}\beta^{\;(j)}_n)\right]\Big| = 0.$$

Now suppose we consider the case where the mean of our distribution $\mathbb{P}_n$ grows with increasing sample size and
 $$\lim_n \sqrt{n}\beta^{\;(j)}_n = -\infty \text{ for all } j\in [d].$$
As the sample size grows bigger, the probability of the selection outcome on Gaussian data approaches $0$.
This selection outcome is referred to as a rare outcome.

From now on, we focus on a subset of these parameters that result in large deviation-type probabilities.
%Henceforth, we consider a sub-collection of these parameters that lead to probabilities of large-deviations.
Fix $\bar{\beta}>0$.
Suppose that each component of the mean vector is parameterized as
 \begin{equation}
\label{rare:seq}
 \sqrt{n}\beta_n^{\;(j)}= -a_n\bar{\beta},
\end{equation}
where $a_n \to \infty$ as $n\to \infty$ and $a_n = o(n^{1/2})$.
%, and $|\bar{\beta} -\bar{\beta_0}|\leq n^{-1}\cdot R_0$. 
 Using the Mills ratio for Gaussian tail probabilities, it is easy to see that the probability of the rare outcome vanishes to $0$ as:
 $$\mathbb{P}_{\mathcal{N}}[j\in \oE] =\bar{\Phi}\left(- \dfrac{\sqrt{n}\beta^{\;(j)}_n}{\sqrt{(1+\rho^2)}} \right) = C_0(a_n\bar{\beta})^{-1}\cdot  \phi\left(\dfrac{a_n\bar{\beta}}{\sqrt{1+\rho^2}}\right).$$
 %$$\mathbb{E}_{\mathcal{N}} \left[ {\bar{\Phi}\left( -\dfrac{1}{\rho}{(\mathcal{Z}_n^{\;(j)}+\sqrt{n}\beta^{\;(j)}_n)}\right)} \right]=\bar{\Phi}\left(- \dfrac{\sqrt{n}\beta^{\;(j)}_n}{\sqrt{(1+\rho^2)}} \right) = C_0(a_n\bar{\beta})^{-1}\cdot  \phi\left(\dfrac{a_n\bar{\beta}}{\sqrt{1+\rho^2}}\right).$$

\begin{assumption}
Consider a collection of distributions $\mathcal{P}_{r,n}$ that have means parameterized as per \eqref{rare:seq}.
Assume that the collection $\mathcal{P}_{r,n}$ has uniformly bounded exponential moments as follows
$$\displaystyle\sup_n\displaystyle\sup_{\mathbb{P}_n\in \mathcal{P}_{r,n}} \mathbb{E}_{\mathbb{P}_n}\left[ \;\exp(\chi |\mathrm{e}^{\;(j)}_{1,n}|) \;\right]< \infty $$
for some $\chi \in \real^+$.
\label{moment:rare:uni}
\end{assumption}

Let $\Psi: \mathcal{K}\to \real$ be a continuous, bounded function. 
Under the moment condition in Assumption \ref{moment:rare:uni}, the variable $\mathcal{Z}_n$ obeys Varadhan's principle of large deviations in the following sense:
\begin{equation*}
\label{VP:uni}
\frac{1}{a_n^2} \log \mathbb{E}_{\mathbb{P}_n}\left[ \exp\left(-a_n^2 \Psi\left(\frac{1}{a_n}\mathcal{Z}_n\right)\right) \cdot \mathbf{1}_{\frac{1}{a_n}\mathcal{Z}_n\in\mathcal{K}} \right] = \mathrm{r}_{\Psi,n}-\inf_{z\in \mathcal{K}} \; \left\{\frac{1}{2} z^2 + \Psi(z)\right\}, 
\end{equation*}
where $\mathrm{r}_{\Psi,n}=o(1)$.
For example, please see \cite{de1992moderate}.

\begin{assumption}
Consider $\Psi\equiv \Psi_t$ for $t\in \{0,1\}$ where $\Psi_1(z) = \frac{1}{\rho^2}(z-\bar\beta)^2$ and $\Psi_0(z) =0$.
Fix $\mathcal{K}\equiv \mathcal{K}_t$ for $t\in \{0,1\}$ where $\mathcal{K}_1=[-c_0, c_0]$ for $c_0>0$, and $\mathcal{K}_0=\mathcal{K}_1^c$.
We assume that
$$\sup_n \sup_{\mathbb{P}_n\in \mathcal{P}_{r,n}}\sup_{t\in \{0,1\}}\; a_n^2 \mathrm{r}_{\Psi_t,n} <\infty.$$
\label{moderate:assump:uni}
\end{assumption}

The conditions in Assumptions \ref{moment:rare:uni} and \ref{moderate:assump:uni} imply that 
$$
\mathbb{E}_{\mathbb{P}_n}\left[ \exp\left(-a_n^2 \Psi_t\left(\frac{1}{a_n}\mathcal{Z}_n\right)\right) \cdot \mathbf{1}_{a_n^{-1}\mathcal{Z}_n\in\mathcal{K}_t} \right] \leq \rK_0 \exp\left(-a_n^2\inf_{z\in \mathcal{K}_t} \; \left\{\frac{1}{2} z^2 + \Psi_t(z)\right\}\right),
$$
where $\rK_0$ is a constant.
As a result, we obtain the rate of decay for large-deviations type probabilities and exponentially vanishing moments.
%In other words, these conditions lead us to an  obtain the rate of decay for certain exponential moments and the probability of .

\begin{theorem}[Weak convergence under rare outcomes]
\label{weak:convergence:rare}
Suppose that the conditions in Assumptions \ref{moment:rare:uni} and \ref{moderate:assump:uni} are met.
Then, we have 
$$\lim_n \sup_{\mathbb{P}_n\in \mathcal{P}_{r,n}} \widetilde\rR^{(l)}_n =0 \;\;  \text{ for } l\in [2],$$
and
$$\lim_n \sup_{\mathbb{P}_n\in \mathcal{P}_{r,n}}\Big|\widetilde{\mathbb{E}}_{\mathbb{P}_n}\left[\rH\circ\rP^{\;(j)}(\mathcal{Z}^{\;(j)}_n;  \sqrt{n}\beta_n^{\;(j)})\right]- \widetilde{\mathbb{E}}_{\mathcal{N}}\left[\rH\circ\rP^{\;(j)}(\mathcal{Z}_n^{\;(j)}; \sqrt{n}\beta_n^{\;(j)})\right]\Big| = 0.$$
\end{theorem}

%Theorem \ref{weak:convergence:rare} allows us to claim weak convergence of our pivot for a rare outcome, i.e.,
%$$\lim_n \sup_{\mathbb{P}_n\in \mathcal{P}_{r,n}}\Big|\widetilde{\mathbb{E}}_{\mathbb{P}_n}\left[\rH\circ\rP^{\;(j)}(\mathcal{Z}^{\;(j)}_n;  \sqrt{n}\beta_n^{\;(j)})\right]- \widetilde{\mathbb{E}}_{\mathcal{N}}\left[\rH\circ\rP^{\;(j)}(\mathcal{Z}_n^{\;(j)}; \sqrt{n}\beta_n^{\;(j)})\right]\Big| = 0.$$
Relative to the conditions in Assumption \ref{moment:bdd:uni}, we impose stronger moment conditions to handle rare outcomes. 
In return, we can guarantee asymptotically-valid inference through our pivot, even when we condition on rare outcomes with large deviation-type probabilities.

\begin{remark}
We exclude the uninteresting case when 
$$\sqrt{n}\beta^{\;(j)}_n \to \infty.$$ 
This is because selection does not have an impact in large samples and standard inferences do not require an adjustment for selection. 
\end{remark}

 \subsection{Main tool for weak convergence theory}
We present the Stein bound for Gaussian approximations, which is the primary tool in our asymptotic theory. 
We then provide a brief outline of how it applies to our problem.

Fixing some more notations, we denote by
$$\mathcal{Z}_n[-i] = \mathcal{Z}_n - Z_{i,n}=  \sum_{k\in [n]\setminus i}  Z_{k,n}$$ 
the $i^{\text{th}}$ leave-one out variable.
This variable is obtained by dropping $Z_{i,n}$ from the sum defined in \eqref{iid:sum}.
Let $\mathcal{Z}^{\;(j)}_n[-i]$ be the $j^{\text{th}}$ entry of this $i^{\text{th}}$ leave-one out variable.
Consider a real-valued mapping $\rg$ that is Lebesgue-almost surely differentiable and satisfies $\mathbb{E}_{\mathcal{N}}[|\rg(Z)|]<\infty$ for $Z\sim \mathcal{N}(0,1)$.
Define
\begin{equation}
\label{Stein:function}
\cS_{\rg}(z) := \exp\left(\frac{1}{2}z^2\right)\cdot \int_{-\infty}^{z} \left\{\rg(t) -\mathbb{E}_{\mathcal{N}}(\rg(Z))\right\}\cdot \tE(t,1) dt,
\end{equation}
which is also called the Stein function for $\rg$.
For $i\in [n]$, we let
$$\mathrm{M}_i(t) = \mathbb{E}_{\mathbb{P}_n}\left[Z_{i,n}^{\;(j)} \left(\mathbf{1}_{[t, \infty)}(Z_{i,n}^{\;(j)}) \mathbf{1}_{[0, \infty)}(t)- \mathbf{1}_{(-\infty, t]}(Z_{i,n}^{\;(j)})\mathbf{1}_{(-\infty, 0)}(t) \right) \right].$$

Lemma \ref{Stein:bound} provides a bound to measure the difference between the expectations of a Gaussian variable and its non-Gaussian counterpart using these notations.
For related literature, we point out to \cite{chen2011normal}.
In this paper, we use the symbol $\mathcal{D}^{k}f(x_0)$ to denote the $k^{\text{th}}$ derivative of a differentiable function $f$ at $x_0$.
\begin{lemma}[Univariate Stein bound]
We have
$$\Big|\mathbb{E}_{\mathbb{P}_n}\left[\rg(\mathcal{Z}_n^{\;(j)})\right]- \mathbb{E}_{\mathcal{N}}\left[\rg(\mathcal{Z}_n^{\;(j)})\right]\Big|\leq \text{\normalfont SB}_{\mathbb{P}_n}(\rg)$$
where
\begin{equation*}
\begin{aligned}
\text{\normalfont SB}_{\mathbb{P}_n}(\rg) &= n\cdot \displaystyle\int_{-\infty}^{\infty}  \sup_{\alpha \in [0, 1]}  \mathbb{E}_{\mathbb{P}_n}\Bigg[\left(|t| + \frac{1}{\sqrt{n}}|{\mathrm{e}}^{\;(j)}_{1,n}|\right)\\
& \;\;\;\;\;\;\;\;\;\;\;\;\;\;\;\;\;\;\;\;\;\;\;\;\;\;\;\;\;\;\;\;\;\;\;\;\;\;\;\;\;\;\;\;\;\times\Big|\cD^2\cS_{\rg}\left(\alpha t+(1-\alpha) \frac{1}{\sqrt{n}}{\mathrm{e}}^{\;(j)}_{1,n}  + \mathcal{Z}_n^{\;(j)}[-1]\right)\Big|\Bigg] \mathrm{M}_1(t) dt.
\end{aligned}
\end{equation*}
\label{Stein:bound}
\end{lemma}

Equipped with the above bound, we review the relative differences defined in \eqref{carved:pivot:sequence:relative}.
We use the Stein bound to write the following inequality:
%We start off by noting that the common denominator of our relative differences is equal to $\mathbb{P}_{\mathcal{N}}[j\in \oE]$.
%Applying the Stein bound to the bounded functions $\rG_l$, we can write
$$
\widetilde\rR_n^{(l)} \leq (\widetilde{D}_n)^{-1}\cdot  \text{\normalfont SB}_{\mathbb{P}_n}(\widetilde\rG_l)
%\left(\mathbb{E}_{\mathcal{N}} \left[ {\bar{\Phi}\left( -\dfrac{1}{\rho}{(\mathcal{Z}_n^{\;(j)}+\sqrt{n}\beta^{\;(j)}_n)}\right)} \right]\right)^{-1} \cdot  \text{\normalfont SB}_{\mathbb{P}_n}(\widetilde\rG_l).
$$
for $l\in [2]$.

First, we consider bounded outcomes.
The probability of a bounded outcome, which is also the common denominator of our relative differences $\widetilde{D}_n$, is bounded away from $0$. 
To prove weak convergence of our pivot, we need to prove that the univariate Stein bound $\text{\normalfont SB}_{\mathbb{P}_n}(\widetilde\rG_l)$ uniformly converges to $0$ as $n$ tends to infinity.
%Thus, weak convergence of our pivot follows once we show that the univariate Stein bound $\text{\normalfont SB}_{\mathbb{P}_n}(\rG_l)$ uniformly converges to $0$ as $n\to \infty$.
When dealing with rare outcomes, the uniform convergence of the Stein bound is not enough to guarantee weak convergence of our pivot.
This is because the probability of the selection outcome also converges to $0$ at an exponentially fast rate.
To ensure weak convergence of our pivot, it is necessary for the related Stein bound to converge at a faster rate than the probability of the selection outcome.
For both types of outcomes mentioned, we investigate the large-sample behavior of the commensurate Stein bound to prove Theorem \ref{weak:convergence:local} and \ref{weak:convergence:rare}.
Detailed proofs are deferred to the Appendix. 

To conclude this section, we examine the smoothness properties of our pivot in Proposition \ref{growth:rate:pivot}. 
This result helps us study the behavior of the Stein bound in our proofs.
\begin{proposition}
Consider the univariate pivot 
$$\rP^{\;(j)} \left(\mathcal{Z}_n^{\;(j)}; \sqrt{n}\beta_n^{\;(j)}\right).$$ 
Then, the pivot's first derivative is uniformly bounded for all real-valued sequences of the mean parameter.
\label{growth:rate:pivot}
\end{proposition}

%%%%%%%%%%%%%%%%%%%%%%%%%%%%%%%%%%%%%%%%%%%%%%%%%%%%%%%%%%%%%
 %%%%%%%%%%%%%%%%%%%%%%%%%%%%%%%%%%%%%%%%%%%%%%%%%%%%%%%%%%%%%

\section{Weak convergence of multivariate pivot}
\label{weak:convergence:Gaussian:multi}
We turn to the multivariate pivot in Proposition \ref{carved:pivot}.
%For ease of reference, Table \ref{symbols} collects the symbols that are used in this section. 
Throughout the section, we will use $C_1$, $C_2$, $\cdots$ to denote constants that are free of $n$.

%The weak convergence theory for our univariate pivot is a blueprint of the analysis with the more involved multivariate pivot.
%The generalization, however, is not trivial as major differences emerge due to the dependence shared between the components of $\beta_n$.

\subsection{Main results}

In line with the preceding section, we develop our theory for bounded and rare outcomes.

We start from considering mean parameters which satisfy
  \begin{equation}
  \label{bdd:par}
  \|\sqrt{n}{\beta}_n-r\|\leq R.
  \end{equation}
Suppose that $\mathbb{P}_n =\mathcal{N}(\beta_n, \Sigma)$.
Recall that the probability of the selection outcome is equal to
$$\mathbb{P}_{\mathcal{N}}\left[ T_n >0_{p}\right],$$
where $T_n$ is a Gaussian variable as stated in Proposition \ref{rep:sel:probability}.
It is easy to see that the probability of the selection outcome is bounded away from $0$, which gives rise to bounded outcomes.

\begin{assumption}
\label{moment:bdd:multi}
We consider a collection of distributions $\mathcal{P}_{b,n}$ with bounded mean parameters as stated in \eqref{bdd:par}.
Suppose that the collection $\mathcal{P}_{b,n}$ has uniformly bounded moments as follows
$$\displaystyle\sup_n\displaystyle\sup_{\mathbb{P}_n\in \mathcal{P}_{b,n}} \mathbb{E}_{\mathbb{P}_n}\left[ \|\mathrm{e}_{1,n}\|^6 \right]< \infty.$$
\end{assumption}

Let $\rR^{(l)}_n$ be defined according to Proposition \ref{weak:convergence:relative:multi}.
Theorem \ref{weak:convergence:local:multi} assures that our pivot generates asymptotically-valid selective inference for bounded outcomes.
%proves weak convergence of the multivariate pivot under bounded outcomes.

\begin{theorem}[Weak convergence under bounded outcomes]
\label{weak:convergence:local:multi}
Under Assumption \ref{moment:bdd:multi}, we have 
$$\lim_n \sup_{\mathbb{P}_n\in \mathcal{P}_{b,n}} \rR^{(l)}_n =0 \;\;  \text{ for } l\in [2],$$
and as a result,
$$\lim_n \sup_{\mathbb{P}_n\in \mathcal{P}_{b,n}}\Big|\widetilde{\mathbb{E}}_{\mathbb{P}_n}\left[\rH\circ\rP^{\;(j)}(\mathcal{Z}_n;  \sqrt{n}\beta_n)\right]- \widetilde{\mathbb{E}}_{\mathcal{N}}\left[\rH\circ\rP^{\;(j)}(\mathcal{Z}_n; \sqrt{n}\beta_n)\right]\Big| = 0.$$
\end{theorem}

%The guarantee in \eqref{sub:qn:validity}, and subsequently weak convergence of the multivariate pivot follow immediately for bounded outcomes.

Now we turn to rare outcomes.
Fix $\bar{\beta}\in \real^d$ such that $\bar{\Sigma} Q' \Sigma^{-1}\bar{\beta} \notin (-\infty, 0]^d$.
%Fix $R_0>0$.
%, let $\bar{\mu}_{0}= \bar{\Sigma} Q' \Sigma^{-1} (\sqrt{n}\beta_0 - r)$.
%Now, fix $\bar{\beta}_0$ such that $\bar{\mu}_{0} \notin [0, \infty)^d$, 
Let the mean for our generating distribution $\mathbb{P}_n$ be parameterized as
   \begin{equation}
\begin{aligned}
& \sqrt{n}\beta_n-r = -a_n\bar{\beta},
\end{aligned}
\label{rare:par}
\end{equation}
where $a_n \to\infty$ as $n \to\infty$ and $a_n=o(n^{1/6})$.
% and $\|\bar{\beta}-\bar{\beta}_0\|\leq n^{-1/3} \cdot R_0$.
%which takes the form

For each $\beta_n$, we consider the matching parameter 
$$\bar{\mu}_{n}= \bar{\Sigma} Q' \Sigma^{-1} (\sqrt{n}\beta_n - r).$$
Based on our parameterization, note that we can write
$$\bar{\mu}_{n}=  -a_n \bar{\mu},$$
where $\bar{\mu}= \bar{\Sigma} Q' \Sigma^{-1} \bar{\beta}$.
Formalized next, we first see that the probability of the selection outcome vanishes to zero at an exponentially fast rate.  

\begin{proposition}[Probability of a rare outcome]
\label{rate:decay:Gaussian}
Consider the optimization problem
$$t_{\star} = \underset{\; t\geq \bar{\mu}}{\text{argmin}} \;\; t' \bar{\Sigma}^{-1} t.$$
Then, there exists a unique (non empty) set $\mathcal{I}\subseteq [d]$ such that the following assertions are simultaneously true:
\begin{enumerate}[label=(\roman*)]
\setlength\itemsep{1em}
  \item $t_{\star}^{(\mathcal{I})} =\bar{\mu}^{(\mathcal{I})} \neq 0_{|\mathcal{I}|}$;
  \item for $\mathcal{J}=\mathcal{I}^c$, \; $t_{\star}^{(\mathcal{J})} =\bar{\Sigma}_{\mathcal{J}, \mathcal{I}}\bar{\Sigma}_{\mathcal{I}, \mathcal{I}}^{-1} \bar{\mu}^{(\mathcal{I})} \geq \bar{\mu}^{(\mathcal{J})} \text{ whenever } \mathcal{J}\neq \emptyset;$
  \item $\left(\bar{\Sigma}_{\mathcal{I}, \mathcal{I}}^{-1} \bar{\mu}^{(\mathcal{I})}\right)^{\;(j)}>0$ for all $j\in \mathcal{I}$ \text{ and } $t'_{\star}\bar{\Sigma}^{-1}t_{\star}= (\bar{\mu}^{(\mathcal{I})})'\bar{\Sigma}^{-1}_{\mathcal{I}, \mathcal{I}}\bar{\mu}^{(\mathcal{I})}>0$.\\
\end{enumerate}
Further, we have
\begin{equation*}
\begin{aligned}
\mathbb{P}_{\mathcal{N}}\left[ T_n >0_{p}\right]= \dfrac{C_3}{(a_n)^{|\mathcal{I}|}}\cdot \text{\normalfont Exp}\left(a_n \bar{\mu}^{(\mathcal{I})}, \frac{1}{(1+\rho^2)} (\bar{\Sigma}_{\mathcal{I}, \mathcal{I}})^{-1}\right)
\end{aligned}
\end{equation*}
for sufficiently large $n$.
\end{proposition}

\begin{remark}
%Note that the probability of our selection outcome is a Gaussian tail probability. 
The proof for the above result closely follows Proposition 2.1 and Corollary 4.1 in \citep{hashorva2003multivariate}.
%because the probability of our selection outcome is a Gaussian tail probability.
Therefore, we omit further details of the proof here.
\end{remark}

%\begin{remark}
%Observe that the set $\mathcal{I}$ in the above-stated result depends on the parameter $\bar{\mu}$. 
%However, we choose to suppress this dependence here.
%This is because the probability of the selection outcome converges to zero at the same rate as long as the underlying mean vector is parameterized according to \eqref{rare:par}.
%\end{remark}

As a corollary, we observe the following. 
\begin{corollary}
\label{rep:den}
Let $\Delta = \Sigma^{-1}Q\bar{\Sigma}_{\mathcal{I}}\bar{\Sigma}_{\mathcal{I}, \mathcal{I}}^{-1} \bar{\Sigma}'_{\mathcal{I}} Q' \Sigma^{-1}$.
It holds that the common denominator of our relative differences is equal to
\begin{equation*}
\begin{aligned}
&\mathbb{E}_{\mathcal{N}}\left[\rF \left( \Sigma^{1/2}\mathcal{Z}_n + \sqrt{n}\beta_n\right)\right] =\frac{C_4}{(a_n)^{|\mathcal{I}|}}\cdot \text{\normalfont Exp}\left(\sqrt{n}\beta_n-r, \frac{1}{(1+\rho^2)}\cdot (\Lambda+ \Delta)\right).
\end{aligned}
\end{equation*}
\end{corollary}

The proof of Corollary \ref{rep:den} follows directly from the claims in Proposition \ref{rep:sel:probability} and \ref{rate:decay:Gaussian}.

\begin{assumption}
Consider a collection of distributions $\mathcal{P}_{r,n}$ such that the mean grows with $n$ as per \eqref{rare:par}. 
Assume that the collection $\mathcal{P}_{r,n}$ has a uniformly bounded exponential moment near the origin as follows:
$$\displaystyle\sup_n\displaystyle\sup_{\mathbb{P}_n\in \mathcal{P}_{r,n}}\mathbb{E}_{\mathbb{P}_n}[\; \exp(\chi\|\mathrm{e}_{1,n}\| )\;] <\infty $$
 for some $\chi\in \real^+$.
\label{moment:rare:multi}
\end{assumption}

Let $\Psi: \mathcal{K}\to \real$ be a continuous and bounded function.
Under assumption \ref{moment:rare:multi}, Varadhan's principle of large deviations for $\mathcal{Z}_n$ implies that
\begin{equation*}
\frac{1}{a_n^2}\log \mathbb{E}_{\mathbb{P}_n}\left[ \exp\left(-a_n^2 \Psi\left(\frac{1}{a_n}\mathcal{Z}_n\right)\right) \cdot \mathbf{1}_{\frac{1}{a_n}\mathcal{Z}_n\in\mathcal{K}} \right] = \mathrm{r}_{\Psi,n}-\inf_{z\in \mathcal{K}} \; \left\{\frac{1}{2} z'z + \Psi(z)\right\},
\label{VP:multi}
\end{equation*}
where $\mathrm{r}_{\Psi,n}=o(1)$.

\begin{assumption}
Consider $\Psi\equiv\Psi_t$ where $\Psi_t(z)=  \frac{1}{1-t +\rho^2}(\sqrt{t}\Sigma^{1/2}z-\bar\beta)'(\Lambda +\Delta)(\sqrt{t}\Sigma^{1/2}z-\bar\beta)$ for $t\in (0,1]$ and $\Psi_0(z)=0$.
Fix $\mathcal{K}\equiv \mathcal{K}_t$ where $\mathcal{K}_t=[-c_0 \cdot 1_d, c_0 \cdot 1_d]$ for $c_0>0$ and $t\in (0,1]$, and $\mathcal{K}_0=\mathcal{K}_1^c$.
%Suppose that $\mathcal{K}$ is a compact subset of $\mathbb{R}^d$.
We impose the condition that
$$\sup_n \sup_{\mathbb{P}_n\in \mathcal{P}_{r,n}} \sup_{t\in [0,1]}\; a_n^2 \mathrm{r}_{\Psi_t,n} <\infty.$$
\label{moderate:assump:multi}
\end{assumption}

\begin{assumption}
%Let $\mathcal{M}_{r,n}$ be the collection of mean parameters in \eqref{rare:par} for any fixed $n\in \mathbb{N}$.
Additionally, we assume that
\begin{equation*}
\sup_n\displaystyle\sup_{\mathbb{P}_n\in \mathcal{P}_{r,n}} \; \dfrac{ \mathbb{E}_{\mathbb{P}_n}[\rF(\Sigma^{1/2}\mathcal{Z}_n + \sqrt{n}\beta_n) \cdot \mathbf{1}_{\mathcal{Z}_n \in \mathcal{K}_n}]}{\mathbb{E}_{\mathcal{N}}[\rF(\Sigma^{1/2}\mathcal{Z}_n + \sqrt{n}\beta_n) \cdot \mathbf{1}_{\mathcal{Z}_n \in \mathcal{K}_n}]} <\infty
\end{equation*}
whenever
%$$\lim_n  \sup_{\beta_n \in \mathcal{M}_{r,n}}\; \dfrac{\mathbb{E}_{\mathcal{N}}\Big[\rF(\Sigma^{1/2}\mathcal{Z}_n + \sqrt{n}\beta_n) \cdot \mathbf{1}_{\mathcal{Z}_n \in \mathcal{K}_n}\Big]}{\mathbb{E}_{\mathcal{N}}\Big[\rF(\Sigma^{1/2}\mathcal{Z}_n + \sqrt{n}\beta_n)\Big]} =0.$$
%$$\lim_n  \sup_{\beta_n \in \mathcal{M}_{r,n}}\;  \dfrac{\mathbb{P}_{\mathcal{N}}\left[ E_n= \oE, \; A_n= \oA,\; \bar{\mathcal{K}}_n\right]}{\mathbb{P}_{\mathcal{N}}\left[ E_n= \oE, \; A_n= \oA\right]} =0.$$
$$\lim_n  \; \dfrac{\mathbb{E}_{\mathcal{N}}\Big[\rF(\Sigma^{1/2}\mathcal{Z}_n + \sqrt{n}\beta_n) \cdot \mathbf{1}_{\mathcal{Z}_n \in \mathcal{K}_n}\Big]}{\mathbb{E}_{\mathcal{N}}\Big[\rF(\Sigma^{1/2}\mathcal{Z}_n + \sqrt{n}\beta_n)\Big]} =0.$$
\label{tail:behavior}
\end{assumption}

Consistent with the weak convergence theory in the earlier section, we require stronger moment conditions to guarantee weak convergence of the multivariate pivot under rare outcomes.
In particular, we note the following.
\begin{remark}
Similar to our univariate analysis, the conditions in Assumptions \ref{moment:rare:multi} and \ref{moderate:assump:multi} provide a uniform bound on a set of large-deviations type probabilities and exponentially vanishing moments.
The condition in Assumption \ref{tail:behavior} controls the probability of selection outcomes that are rarer than the observed outcome on Gaussian data by imposing the restriction that these probabilities decay at an equal or faster rate than the limiting Gaussian counterpart 
More specifically, this condition allows us to establish convergence of our relative differences on a set of high probability while controlling their behavior on the complement set.
\end{remark}

Theorem \ref{weak:convergence:rare:multi} proves that our pivot offers asymptotically-valid selective inference, even when rare outcomes are observed.
%We show that the weak convergence statement in \eqref{sub:qn:validity} holds uniformly over all distributions in the collection $\mathcal{P}_{r,n}$.

\begin{theorem}[Weak convergence under rare outcomes]
\label{weak:convergence:rare:multi}
Suppose that the conditions in Assumptions \ref{moment:rare:multi}, \ref{moderate:assump:multi}, and \ref{tail:behavior} are met.
Then, we have that
$$\lim_n \sup_{\mathbb{P}_n\in \mathcal{P}_{r,n}} \rR^{(l)}_n =0,$$
and that
$$\lim_n \sup_{\mathbb{P}_n\in \mathcal{P}_{r,n}}\Big|\widetilde{\mathbb{E}}_{\mathbb{P}_n}\left[\rH\circ\rP^{\;(j)}(\mathcal{Z}_n;  \sqrt{n}\beta_n)\right]- \widetilde{\mathbb{E}}_{\mathcal{N}}\left[\rH\circ\rP^{\;(j)}(\mathcal{Z}_n; \sqrt{n}\beta_n)\right]\Big| = 0$$
for $l\in [2]$.
\end{theorem}

\subsection{Main tool for weak convergence theory}
To prove our main results in Theorem \ref{weak:convergence:local:multi} and Theorem \ref{weak:convergence:rare:multi}, we use a multivariate version of the Stein bound.
% to demonstrate weak convergence of our pivot.

Lemma 2 presents this bound for a Lebesgue-almost surely three times differentiable mapping $\rg:\real^d\to \real$, which is adopted from \cite{chatterjee2007multivariate}.
Suppose that $\mathbb{E}_{\mathcal{N}}[|\rg(Z)|]<\infty$. 
Let $Z\sim N(0_d, I_{d,d})$.
The Stein bound is defined through partial derivatives of
\begin{equation*}
\cS_{\rg}(z) = \int_0^1 \frac{1}{2t} \left(\mathbb{E}_{\mathcal{N}}\left[ \rg(\sqrt{t} z + \sqrt{1-t} Z)\right] - \mathbb{E}_{\mathcal{N}}\left[ \rg(Z)\right] \right)dt,
\end{equation*}
also called the Stein function for $\rg$.
Before stating the bound, recall that 
$$\mathcal{Z}_n[-i]= \mathcal{Z}_n - Z_{i,n}$$
denotes the $i^{\text{th}}$ leave-one out variable.
%Consider $f:\mathcal{C}\to \real$ for an open set $\mathcal{C}\subset \real^d$.
Let 
$$\mathcal{D}^k f(x_0)[i_1, i_2,\cdots, i_k]= \dfrac{\partial^k f (x_0)}{\partial x^{(i_1)} \partial x^{(i_2)} \ldots \partial x^{(i_k)}} $$
denote the $k^{\text{th}}$ order partial derivative of $f$ at $x_0$, for $i_1, i_2, \ldots, i_k \in [d]$, and 
 let $\mathrm{e}^{\star}_{i,n}$ be an independent copy of $\mathrm{e}_{i,n}$, for $i\in [d]$.
\begin{lemma}[Multivariate Stein bound]
We have that
$$|\mathbb{E}_{\mathbb{P}_n}\left[\rg(\mathcal{Z}_n)\right] - \mathbb{E}_{\mathcal{N}}\left[\rg(\mathcal{Z}_n)\right]| \leq \text{\normalfont SB}_{\mathbb{P}_n}(\rg)$$
where
\begin{equation*}
\begin{aligned}
&\text{\normalfont SB}_{\mathbb{P}_n}(\rg)=\frac{C_1}{\sqrt{n}}\sum_{\lambda, \gamma \in \{0\} \cup[3]: \lambda+\gamma\leq 3} \sum_{j, k, l} \mathbb{E}_{\mathbb{P}_n}\Big[\|\mathrm{e}_{1,n}\|^{\lambda}\|\mathrm{e}^{\star}_{1,n}\|^{\gamma} \sup_{\alpha, \kappa\in [0,1]}\Big| \cD^3\cS_{\rg}\Big(\mathcal{Z}_n[-1] \\
&\;\;\;\;\;\;\;\;\;\;\;\;\;\;\;\;\;\;\;\;\;\;\;\;\;\;\;\;\;\;\;\;\;\;\;\;\;\;\;\;\;\;\;\;\;\;\;\;\;\;\;\;\;\;\;\;\;\;\;\;\;\;\;\;\;\;\;\;\;\;\;\;\;\;\;\;\;\;\;\;\;\;\;\;\;\;\;\;\;\;\;\;\;\;\;\;\;\;\;\;\; + \frac{\alpha}{\sqrt{n}}\mathrm{e}_{1,n}+ \frac{\kappa}{\sqrt{n}}\mathrm{e}^{\star}_{1,n}\Big)[j, k, l]\Big|\Big].
\end{aligned}
\end{equation*}
\label{Stein:bound:direct}
\end{lemma}

As before, we revisit our relative differences and use the Stein bound to note that
%We can use the Stein bound on the bounded functions $\rG_l$ (see Proposition \ref{weak:convergence:relative:multi}) to bound our relative differences as
\begin{equation*}
R_n^{(l)} \leq \left(\mathbb{E}_{\mathcal{N}}\left[\rF \left( \Sigma^{1/2}\mathcal{Z}_n + \sqrt{n}\beta_n\right)\right] \right)^{-1} \cdot  \text{\normalfont SB}_{\mathbb{P}_n}(\rG_l).
\end{equation*}
To establish the weak convergence of our pivot, we analyze how the Stein bound behaves in large samples, similar to what we did for the univariate pivot.
Detailed proofs for Theorem \ref{weak:convergence:local:multi} and Theorem \ref{weak:convergence:rare:multi} are developed in the Appendix.

In line with Section \ref{weak:convergence:Gaussian}, we obtain the smoothness properties of our pivot in Proposition \ref{smoothness:property:pivot}. 
\begin{proposition}
\label{smoothness:property:pivot}
Fix $p_0\in \mathbb{N}$. We have 
\[\|\cD^{p_0} \rP^{\;(j)} \left(\mathcal{Z}_n; \sqrt{n}\beta_n\right)\|\leq \sum_{\lambda, \gamma \in \{0\} \cup[p_0]: \lambda+\gamma\leq p_0} C_2^{\lambda, \gamma}\|\mathcal{Z}_n\|^{\lambda} \|\sqrt{n}\beta_n\|^{\gamma}.\]
\end{proposition}

\begin{remark}
In contrast to the univariate theory, the multivariate version of the Stein bound involves higher order derivatives of the Stein function.
As a result, we investigate higher order smoothness properties of our multivariate pivot.
\end{remark}
%\begin{remark}
%Deviating from the weak convergence theory of our univariate pivot, the partial derivatives of the multivariate pivot are no longer uniformly bounded over their domain sets.
%\end{remark}

\subsection{Transfer of asymptotic guarantees to the carved pivot}
Having established weak convergence of our pivot for randomized rules with Gaussian variables, we come back to the selection described in \eqref{implicit:randomization:FE}. 

Following the same convention as before, we evaluate the the likelihood ratio after and before we apply the selection rule on the pilot samples.
At $\begin{pmatrix} v' & w'\end{pmatrix}'$, let the joint density for $V_n$ and $W_n$ factorize as 
$$\mathrm{p}_n(v, w)= \mathrm{p}_n(v)\cdot \bar{\mathrm{p}}_n(w\lvert v),$$
where $\mathrm{p}_n$ is the marginal density for $V_n$ and $\bar{\mathrm{p}}_n(\cdot \lvert v)$ is the conditional density of $W_n$ given $V_n=v$.
Let $\bar{\rF}_n:\real^d \to \real$ assume the value 
 \begin{align*}
& \bar\rF_n(v) = \int \bar{\mathrm{p}}_n(Qt  -v + r\lvert v) \cdot \mathbf{1}_{t\in \real^{p+}} dt.
 \end{align*} 
\begin{proposition}
\label{lik:ratio:carved}
Under the randomized selection rule in \eqref{implicit:randomization:FE}, the ratio of the conditional and unconditional likelihood functions is 
$$\overline{\text{\normalfont LR}}_{\mathbb{P}_n}(\oZ; \sqrt{n}\beta_n) = \dfrac{\bar\rF_n(\Sigma^{1/2}\oZ +\sqrt{n}\beta_n)}{\mathbb{E}_{\mathbb{P}_n}\left[\bar\rF_n(\Sigma^{1/2}\mathcal{Z}_n +\sqrt{n}\beta_n)\right]}.$$
\end{proposition}
Define
$$\overline{\mathbb{E}}_{\mathbb{P}_n}\left[\mathcal{Q}(\mathcal{Z}_n)\right] = \mathbb{E}_{\mathbb{P}_n}\left[\mathcal{Q}(\mathcal{Z}_n)\cdot {\overline{\text{LR}}}_{\mathbb{P}_n}(\mathcal{Z}_n; \sqrt{n}\beta_n)\right].$$
The expectation on the left-hand side is taken with respect to the conditional law after selection on pilot data and is expressed as an unconditional expectation on the right-hand side through the above-stated likelihood ratio.

Consider a collection of distributions $\mathcal{C}_n$.
The weak convergence of our pivot follows by proving
\begin{equation}
\label{sub:qn:validity:carved}
\lim_n \sup_{\mathbb{P}_n\in \mathcal{C}_n}\Big|\overline{\mathbb{E}}_{\mathbb{P}_n}\left[\rH\circ\rP^{\;(j)}(\mathcal{Z}_n;  \sqrt{n}\beta_n)\right]- \widetilde{\mathbb{E}}_{\mathcal{N}}\left[\rH\circ\rP^{\;(j)}(\mathcal{Z}_n; \sqrt{n}\beta_n)\right]\Big| = 0
\end{equation}
for any $\rH\in \mathbb{C}^3(\real, \real)$ with bounded derivatives up to the third order.
We substituted the first term in \eqref{sub:qn:validity} with a conditional expectation that relies on the distribution post conditioning on the selection outcome observed in the pilot data.

Our next result establishes that asymptotically-valid selective inference with Gaussian randomized rules transfers to the carved pivot.
This result holds as long as the probability of the selection outcome converges to its counterpart with Gaussian randomization.
\begin{proposition}[Transfer of asymptotic guarantees]
\label{carved:CLT}
Suppose that the conditional weak convergence statement in \eqref{sub:qn:validity} holds over a collection of distributions in $\mathcal{C}_n$.
Assume that
$$\lim_n \sup_{\mathbb{P}_n\in \mathcal{C}_n}\;\dfrac{\mathbb{E}_{\mathbb{P}_n}\left[|\bar\rF_n(\Sigma^{1/2}\mathcal{Z}_n +\sqrt{n}\beta_n) - \rF(\Sigma^{1/2}\mathcal{Z}_n +\sqrt{n}\beta_n)|\right]}{\mathbb{E}_{\mathbb{P}_n}\left[ \rF(\Sigma^{1/2}\mathcal{Z}_n +\sqrt{n}\beta_n)\right]}=0.$$
We then have the convergence in \eqref{sub:qn:validity:carved}.
\end{proposition}

\section{Empirical analysis}
\label{empirical:results}

We illustrate how our theory translates to practice in various instances of selective inference.

\begin{example}
\label{two:sample:t-test} \rm{Selectively inferring for a difference in means.}
\rm{
We selectively infer for a difference in means through the two-sample test statistic.
In alignment with the running example in our paper, we use the following scheme to draw $n$ independent and identically distributed observations with identity covariance.
For $d=2$, we draw 
$$\zeta_{i,n} = \beta_{n} + \mathrm{e}_{i,n} \text{ for } i\in [n].$$ 
Each component of $\mathrm{e}_{i,n}$ is drawn independently as 
$${\mathrm{e}}^{\;(j)}_{i,n}\stackrel{i.i.d.}{\sim} \rE$$ 
and standardized such that $$\mathbb{E}[{\mathrm{e}}^{\;(j)}_{i,n}]=0;\  \mathbb{E}[({\mathrm{e}}^{\;(j)}_{i,n})^2]= 1.$$
Note that the distribution $\rE$ is based on five different models, which include Models $(1)$-$(4)$ described in Section \ref{framework} and the baseline Gaussian Model.
We provide selective inference for $\bar{\beta}_n = \beta^{(1)}_{n}-\beta^{(2)}_{n}$ whenever the two-sample statistic
$$V_{n_1}= \frac{\sqrt{n_1}}{\sqrt{2}}(\bar{\zeta}^{(1)}_{n_1}-\bar{\zeta}^{(2)}_{n_1}),$$
 exceeds a prefixed threshold of significance.
%The exact pivot for $\bar{\beta}_n$ under the Gaussian model is no different from the running example.
We investigate the performance of our carved pivot for $\bar{\beta}_n$.

For our simulations, the difference of means is parameterized according to $\sqrt{n}\bar{\beta}_n= -a_n \bar{\beta}$ for $a_n = n^{1/6-\delta}$ and $\delta=1\mathrm{e}{-3}$.
We fix $n=50$.
We set our split proportion value at
$$\rho^2= \frac{n-n_1}{n_1}= 1/2,$$
i.e., two-thirds of our data is used to decide whether to pursue inference in the second stage. 
We vary $\bar{\beta}$ in the set $ \{2, 1, 0\}$.
For comparison, we consider asymptotic intervals based on the widely used data splitting.
The latter procedure simply uses the $n_2$ samples that were held out for inference.

We compare the $90\%$-confidence intervals from inverting the carved pivot with the $90\%$-confidence intervals from data splitting and summarize our findings in Table \ref{table:eg1:1}.
Our method is noted as ``Carve'' and  data splitting is noted as ``Split''.
%Cast into the carved setup, we recall that the univariate pivot in Corollary \ref{carved:pivot:sequence} has exact guarantees under a Gaussian generating distribution.
The cells in this table report the empirical coverage rate ``$\text{Cov}$'' of the asymptotic confidence intervals and their lengths ``Len'' when averaged over all our simulations.
The first column in the table notes the performance of the exact confidence intervals under the baseline Gaussian model.
%The remaining columns showcase the performance of the asymptotic pivot in finite samples when we depart from Gaussian data.  

As expected, both procedures approximately achieve the target coverage rate.
However, carving produces tighter intervals than data splitting across all models and all values of $\bar{\beta}$.
\begin{center}
\captionof{table}{\centering{Comparison of inference between carving and data splitting.}}
\vspace{-0.1cm}
\label{table:eg1:1}
\bgroup
\def\arraystretch{1.3}
\scalebox{0.88}{\begin{tabular}{ |c | c c |c c |c c|c c |c c | }
\hline
$\rho^2=1/2 $   & Cov & Len  & Cov  & Len  & Cov  & Len &  Cov  & Len  &  Cov  & Len \\ [0.5ex] 
\hline\hline
 $ \bar{\beta}=2$  & \multicolumn{2}{c}{Gaussian} & \multicolumn{2}{c}{Model-1} & \multicolumn{2}{c}{Model-2} & \multicolumn{2}{c}{Model-3} & \multicolumn{2}{c|}{Model-4} \\
 \hline
Carve  & $89\%$  & $0.73$ & $94\%$ & $0.74$ & $88\%$ & $0.72$ & $90\%$ & $0.67$ & $90.5\%$ & $0.68$  \\ [0.5ex] 
 \hline
Split  & $89.5\%$  & $0.95$ & $95\%$ & $0.95$ & $87\%$ & $0.95$ & $91.5\%$ & $0.95$ & $90\%$ & $0.95$\\ [0.5ex] 
 \hline\hline
$\bar{\beta}=1$  & \multicolumn{2}{c}{Gaussian} & \multicolumn{2}{c}{Model-1} & \multicolumn{2}{c}{Model-2} & \multicolumn{2}{c}{Model-3} & \multicolumn{2}{c|}{Model-4} \\
 \hline\hline
Carve  & $91.5\%$  & $0.72$ & $88.5\%$ & $0.72$ & $92\%$ & $0.74$ & $88\%$ & $0.72$ & $90.5\%$ & $0.72$  \\ [0.5ex] 
 \hline
Split  & $88\%$  & $0.95$ & $90\%$ & $0.95$ & $93\%$ & $0.96$ & $91\%$ & $0.95$ & $91\%$ & $0.95$\\ [0.5ex] 
\hline \hline
 $ \bar{\beta}=0$  & \multicolumn{2}{c}{Gaussian} & \multicolumn{2}{c}{Model-1} & \multicolumn{2}{c}{Model-2} & \multicolumn{2}{c}{Model-3} & \multicolumn{2}{c|}{Model-4} \\
  \hline\hline
Carve  & $90\%$  & $0.59$ & $90\%$ & $0.58$ & $87\%$ & $0.59$ & $88.5\%$ & $0.60$ & $91.5\%$ & $0.58$  \\ [0.5ex] 
 \hline
Split  & $91.5\%$  & $0.95$ & $90\%$ & $0.95$ & $87.5\%$ & $0.95$ & $91\%$ & $0.95$ & $91\%$ & $0.95$\\ [0.5ex] 
 \hline
\end{tabular}}
\egroup
\end{center}
}
\end{example}

\begin{example}{\rm{Selectively inferring for the $p$ largest effects.}}
\label{largest-effects}
\rm{
We consider selective inference for the effects of the $p$ largest mean statistics in our pilot data \citep{subgroupXH}. 
Let $[V_{n_1}]^{(p)}$ be the $p^{\text{th}}$ largest mean statistic using the components of $V_{n_1}$. 
We note that our selection rule in this example can be written as
\begin{equation}
\begin{aligned}
V_{n_1}^{\;(j)}&>  [V_{n_1}]^{(p+1)},\ \text{ for } j\in E_n,\\
V_{n_1}^{\;(j)}&\leq [V_{n_1}]^{(p+1)}, \ \text{ for } j \in E_n^c.
\end{aligned}
\end{equation}
%where we borrow the definition of our randomization variable from \eqref{implicit:randomization:FE}.

Suppose that $\mathbb{P}_n = \mathcal{N}(\beta_n, \Sigma)$. 
Lemma \ref{carved:pivot:eg2} gives a carved pivot after conditioning on the event
$$\{E_n = \oE, A_n = \oA\}$$
where 
$$\oA = \begin{pmatrix}  \left( \sqrt{1+\rho^2}[V_{n_1}]^{(p+1)} \cdot 1_{p}\right)' & \left(V_n^{(E_n^c)}+ W_n^{(E_n^c)}\right)' \end{pmatrix}'=  \begin{pmatrix} A_{1,n}' & A_{2,n}'\end{pmatrix}'.$$ 

To state the pivot, define the matrices
$$R^{\;(j)}= \mathcal{P}_{\oE}\begin{bmatrix} 1 & 0 \\ \frac{1}{\sigma_j^2}\Sigma_{-j,j} & I_{d-1,d-1} \end{bmatrix},\;\; Q = \begin{bmatrix} I_{p, p} \\ 0_{d-p, p} \end{bmatrix}, \;\; r = \oA.$$
\begin{proposition}
\label{carved:pivot:eg2}
Let $\text{\normalfont Pivot}^{\;(j)}\left(V_n^{\;(j)} , U_n^{\;(j)}\right)$ assume the value
$$(\rD(U_n^{\;(j)}; \sqrt{n}\beta_n^{\;(j)}))^{-1}\cdot \int_{V_n^{\;(j)}}^{\infty}\phi\left(\dfrac{1}{\sigma_j}(v- \sqrt{n}\beta_n^{\;(j)})\right)\cdot \rF\left(R^{\;(j)}\begin{pmatrix} v & (U_n^{\;(j)})' \end{pmatrix}'\right) dv.$$
Then, it holds that $\text{\normalfont Pivot}^{\;(j)}\left(V_n^{\;(j)} , U_n^{\;(j)}\right)$ is distributed as a $\text{Unif}\;(0,1)$ conditional on $\{E_n = \oE, A_n = \oA\}$.
\end{proposition}

Clearly, this pivot has the same representation as our running example.

Using the generating scheme from the preceding example, we selectively infer for the effect that corresponds to the larger sample mean.
% out of the two groups of observations, $\bar{\zeta}^{(1)}_{n_1}$ and $\bar{\zeta}^{(2)}_{n_1}$. 
A similar comparison between carving and data splitting unfolds in Table \ref{table:eg2:1} for different models.

\begin{center}
\captionof{table}{\centering{Comparison of inference between carving and data splitting.}}
\vspace{-0.1cm}
\label{table:eg2:1}
\bgroup
\def\arraystretch{1.3}
\scalebox{0.88}{\begin{tabular}{ |c | c c |c c |c c|c c |c c | }
\hline
$\rho^2=1/2$  & Cov & Len  & Cov  & Len  & Cov  & Len &  Cov  & Len  &  Cov  & Len \\ [0.5ex] 
\hline \hline
$\bar{\beta}=-2.5$  & \multicolumn{2}{c}{Gaussian} & \multicolumn{2}{c}{Model-1} & \multicolumn{2}{c}{Model-2} & \multicolumn{2}{c}{Model-3} & \multicolumn{2}{c|}{Model-4} \\
 \hline
 Carve  & $93\%$  & $0.60$ & $88\%$ & $0.64$ & $91.5\%$ & $0.62$ & $91\%$ & $0.61$ & $87\%$ & $0.59$  \\ [0.5ex] 
 \hline
Split  & $90\%$  & $0.95$ & $90.5\%$ & $0.95$ & $91.5\%$ & $0.95$ & $93.5\%$ & $0.95$ & $90.5\%$ & $0.95$\\ [0.5ex] 
 \hline \hline
 $\bar{\beta}=-1.5$  & \multicolumn{2}{c}{Gaussian} & \multicolumn{2}{c}{Model-1} & \multicolumn{2}{c}{Model-2} & \multicolumn{2}{c}{Model-3} & \multicolumn{2}{c|}{Model-4} \\
 \hline
 Carve  & $92\%$  & $0.68$ & $89\%$ & $0.69$ & $90\%$ & $0.68$ & $87.5\%$ & $0.66$ & $90\%$ & $0.66$  \\ [0.5ex] 
 \hline
Split  & $91.5\%$  & $0.95$ & $93\%$ & $0.95$ & $92.5\%$ & $0.95$ & $90.5\%$ & $0.95$ & $92.5\%$ & $0.95$\\ [0.5ex] 
 \hline \hline
 $\bar{\beta}=0$  & \multicolumn{2}{c}{Gaussian} & \multicolumn{2}{c}{Model-1} & \multicolumn{2}{c}{Model-2} & \multicolumn{2}{c}{Model-3} & \multicolumn{2}{c|}{Model-4} \\
 \hline
 Carve  & $91\%$  & $0.55$ & $87.5\%$ & $0.56$ & $87.5\%$ & $0.52$ & $87\%$ & $0.56$ & $90.5\%$ & $0.55$  \\ [0.5ex] 
 \hline
Split  & $92\%$  & $0.95$ & $88.5\%$ & $0.95$ & $88\%$ & $0.95$ & $90\%$ & $0.95$ & $88.5\%$ & $0.95$\\ [0.5ex] 
 \hline
 \end{tabular}}
\egroup
\end{center}
}
\end{example}

\begin{example}
\label{reg:example}
\rm{
We turn to inference for the selected regression coefficients after solving the LASSO. 
Let $y_n$ and $X_n$ denote our response vector and our design matrix with $d$ predictors, respectively.
%Each observation $(y_{i}, x_{i}) \in \real^{(d+1)} $ is drawn as an independent and identically distributed observation. 
% which are based on observing a triangular array of i.i.d. observations
%$$\mathrm{e}_{i,n} = (y_{i}, x_{i}) \in \real^{(d+1)} \sim \mathbb{P}_n, \  \ i\in [n].$$
We start from deriving a pivot under a randomized rule with Gaussian variables.
Consider solving 
\begin{equation}
\label{carved:elastic:net}
\underset{\beta\in\mathbb{R}^d}{\text{minimize}} \;\; \frac{1}{2\sqrt{n}}\|y_{n} - X_{n}\beta\|_2^2 + \lambda \|\beta\|_1  - W_n' \beta,
\end{equation}
where $W_n$ is a Gaussian randomization variable.
This problem has been termed as the randomized LASSO in \cite{harris2016selective}.
%Let $\begin{pmatrix}\widehat{\beta}_{n, \lambda} \\ 0_{d-p} \end{pmatrix}$ denote the LASSO coefficients. 
%This form of randomized  \citep[][]{panigrahi2018selection}.

After observing the selected set of variables $E_n=\oE$, a common model for inference is the selected model
\begin{equation*}
\label{selected:model}
y_n \sim \mathcal{N}(X_{n,\oE}\beta_{n},\sigma^2 I).
\end{equation*}
%Selection in this example is dependent on data through the statistic
%$$V_n = \frac{1}{\sqrt{n}}X_n'y_n.$$ 
Define
 $$\widehat{\beta}_n^{(\oE)}=\left((X_{n}^{(\oE)})' X_{n}^{(\oE)}\right)^{-1}(X_{n}^{(\oE)})' y_n,$$ 
the refitted least squares estimator which is obtained by regressing our response against the selected variables.
Based on the least squares estimator and the selected set of variables, let
\begin{equation}
\label{main:statistic}
\begin{pmatrix}V_n^{(\oE)} \\ V_n^{(\oE^c)} \end{pmatrix} = \begin{pmatrix} \sqrt{n} \widehat{\beta}_n^{(\oE)} \\  \frac{1}{\sqrt{n} }(X_{n}^{(\oE^c)})' (y_n - X_{n}^{(\oE)}\widehat{\beta}_n^{(\oE)})\end{pmatrix},
\end{equation}
and let
$$V^{\;(j)}_n = e_j'\sqrt{n}\widehat{\beta}_n^{(\oE)},$$
which is the $j^{\text{th}}$ regression coefficient in the selected set.

Fixing some more notations, let 
$$\begin{pmatrix}\widehat{\beta}_{n, \lambda} \\ 0_{d-p} \end{pmatrix}$$ 
denote the coefficients of the LASSO solution, where $\widehat{\beta}_{n, \lambda}$ collects its nonzero coefficients.
Let $S_n^{(E_n)}$ collect the signs of the nonzero LASSO coefficients. 
Let  $\mathcal{G}_n^{(E_n^c)}$ collect the components of the subgradient from the LASSO penalty present in the inactive set $E_n^c$ at the solution. 
Define
$$A_n =  \begin{pmatrix} A_{1,n}' & A_{2,n}'\end{pmatrix}' = \begin{pmatrix} \lambda\cdot (S_n^{(E_n)})' & (\mathcal{G}_n^{(E_n^c)})'\end{pmatrix}',$$
which we note is equal to subgradient of the LASSO penalty at the solution.
Finally, let $T_n= \text{diag}(S_n^{(E_n)})\widehat{\beta}_{n, \lambda}$ collect the magnitudes of the nonzero LASSO coefficients.

Based on these notations, fix the following matrices
$$P_n = \begin{bmatrix} \frac{1}{n}(X_{n}^{(\oE)})' X_{n}^{(\oE)} & 0_{p, d-p} \\ \frac{1}{n}(X_{n}^{(\oE^c)})' X_{n}^{(\oE)} & I_{d-p, d-p} \end{bmatrix}, Q_n = \begin{bmatrix} \frac{1}{n}(X_{n}^{(\oE)})' X_{n}^{(\oE)} \\ \frac{1}{n}(X_{n}^{(\oE^c)})' X_{n}^{(\oE)} \end{bmatrix} \text{diag}\left(S_n^{(E_n)}\right).$$ 
Let $P = \mathbb{E}_{\mathbb{P}_n}[P_n]$ and $Q = \mathbb{E}_{\mathbb{P}_n}[Q_n]$, and also let $\sigma_j^2= \sigma^2\cdot \Sigma^{(\oE)}_{j,j}$ where 
$$\Sigma^{(\oE)}= \left(\mathbb{E}_{\mathbb{P}_n}\left[\frac{1}{n}(X_{n}^{(\oE)})' X_{n}^{(\oE)}\right]\right)^{-1}.$$

%Assuming that certain properties hold exactly for now, Lemma \ref{carved:pivot:eg3} provides a pivot for exact selective inference.
%We list these properties below. 
Suppose that the randomization variable $W_n$ in \eqref{carved:elastic:net} is drawn from the Gaussian distribution $\mathcal{N}(0_d, \rho^2 \Sigma)$, independently of data, where 
$$\Sigma =\sigma^2\cdot \mathbb{E}_{\mathbb{P}_n}\left[\frac{1}{n}X'_nX_n\right].$$
For now, we assume that
\begin{enumerate}[label=(\roman*)]
 \item \label{itm:3}  the variables in \eqref{main:statistic} are distributed as Gaussian variables, where $V_n^{(\oE)}$ has mean $\sqrt{n}\beta_n$ and covariance $$ \sigma^2\cdot \left(\mathbb{E}_{\mathbb{P}_n}\left[\frac{1}{n}(X_{n}^{(\oE)})' X_{n}^{(\oE)}\right]\right)^{-1} = \sigma^2 \cdot \Sigma^{(\oE)},$$ and $V_n^{(\oE)}$ is independent of $V_n^{(\oE^c)}$.
\item \label{itm:4} the magnitudes of the nonzero LASSO coefficients satisfy:
$$
\begin{pmatrix}{W_n^{\;(\oE)}}'  &  {W_n^{\;(\oE^c)}}' \end{pmatrix}'  = Q T_n + \begin{pmatrix} A_{1,n}' & A_{2,n}'\end{pmatrix}'  - P\begin{pmatrix}{V_n^{\;(\oE)}}' & {V_n^{\;(\oE^c)}}' \end{pmatrix}'.
$$
%along with the sign constraints $T_n \in \real^{p+}$, 
\end{enumerate}
In practice, the variables $V_n$ have an asymptotic Gaussian distribution with the properties listed in \ref{itm:3}, and the equality in \ref{itm:4} holds only up to an $o_p(1)$ remainder term.

Proposition \ref{carved:pivot:eg3} gives a pivot that yields exactly-valid selective inference under the above-stated randomized rule and assumptions.

\begin{proposition}
\label{carved:pivot:eg3}
Let $\text{\normalfont Pivot}^{\;(j)}\left(V_n^{\;(j)} , U_n^{\;(j)}\right)$ assume the value
$$(\rD(U_n^{\;(j)}; \sqrt{n}\beta_n^{\;(j)}))^{-1}\cdot \int_{V_n^{\;(j)}}^{\infty}\phi\left(\dfrac{1}{\sigma_j}(v- \sqrt{n}\beta_n^{\;(j)})\right)\cdot \rF\left(PR^{\;(j)}\begin{pmatrix} v & (U_n^{\;(j)})' \end{pmatrix}'\right) dv,$$
where 
$$\rD(U; \sqrt{n}\beta_n^{\;(j)})= \int_{-\infty}^{\infty}\phi\left(\dfrac{1}{\sigma_j}(v- \sqrt{n}\beta_n^{\;(j)})\right) \cdot \rF\left(PR^{\;(j)}\begin{pmatrix} v & U' \end{pmatrix}'\right) dv.$$
Conditional on $\{E_n= \oE, \; A_n= \oA\}$, $\text{\normalfont Pivot}^{\;(j)}\left(V_n^{\;(j)} , U_n^{\;(j)}\right)$ is distributed as a $\text{Unif}\;(0,1)$ variable.
\end{proposition}

Suppose that our data contains $n$ independent and identically distributed observations.
Next, we will address the standard LASSO problem on a randomly drawn subset of data with a size of $n_1$.
For this, consider solving
\begin{equation}
\label{lasso:carved:now}
\underset{\beta\in\mathbb{R}^p}{\text{minimize}} \;\;  \frac{(1+\rho^2)}{2\sqrt{n}}\|y_{n_1} - X_{n_1}\beta\|_2^2+ \lambda \|\beta\|_1.
\end{equation}
%Yet again, we record the outcome of variable selection 
%$$\{E_n=\oE, A_n=\oA\}.$$ 

 We define
$$W_n = \dfrac{\partial }{\partial \beta}\Big\{\frac{1}{2\sqrt{n}}\|y_{n} - X_{n}\beta\|_2^2-\frac{(1+\rho^2)}{2\sqrt{n}} \|y_{n_1} - X_{n_1}\beta\|_2^2\Big\}\Big\lvert_{\widehat{\beta}^{\lambda}}.$$
Then, as shown by \citep{markovic2016bootstrap, selective_bayesian}, we can rewrite the LASSO optimization problem as \eqref{carved:elastic:net}. 
The randomization variable $W_n$ is asymptotically distributed as $\mathcal{N}(0_d, \rho^2 \Sigma)$ for $\rho^2 = \frac{n_2}{n_1}$.
Additionally, $W_n$ is asymptotically independent of $V_n$.
Also, note that the variables $V_n$ have an asymptotic Gaussian distribution with the properties listed in \ref{itm:3}.
See, for example, Proposition 4.1 in \cite{selective_bayesian} which gives the joint distribution of $W_n$ and $V_n$.

Based on our notations, we can verify that
$$
\begin{pmatrix}{W_n^{\;(\oE)}}'  &  {W_n^{\;(\oE^c)}}' \end{pmatrix}' + \mathrm{O}_n  = Q T_n + \begin{pmatrix} A_{1,n}' & A_{2,n}'\end{pmatrix}'  - P\begin{pmatrix}{V_n^{\;(\oE)}}' & {V_n^{\;(\oE^c)}}' \end{pmatrix}' ,
$$
where $\mathrm{O}_n = o_p(1)$.
In what follows, we ignore the $o_p(1)$ remainder term.
Otherwise, we can always work with the variable 
$$\widetilde{W}_n= W_n +  \mathrm{O}_n,$$
which has the same asymptotic distribution as $W_n$.

Our theory in the paper confirms that the pivot in the earlier Proposition enables us to draw asymptotically-valid inference for the selected regression coefficients.
%We proceed to selectively infer for
%$$\beta_n^{(E)}=  \left(\mathbb{E}_{\mathbb{P}_n}\left[(X_n^{(E)})' X_n^{(E)}\right]\right)^{-1}\mathbb{E}_{\mathbb{P}_n}\left[(X_n^{(E)})' y_n\right],$$ 
Below, we summarize the empirical performance of our pivot in both synthetic and real data experiments.

\noindent\textbf{Synthetic data}.\ \ Fix $n=100$ and $d=50$. 
In each round of our simulations, we draw an $n\times d$ design matrix $X$ such that the rows $x_i\sim \mathcal{N}(0_{d}, \Sigma)$ and $\Sigma_{j,k}=0.40^{|j-k|}$. 
We then draw our response according to the model
$$y_{i} = x'_{i}\beta + \sigma \cdot \mathrm{e}_{i,n},$$ 
by generating the model errors $\mathrm{e}_{i,n}$ in an i.i.d. fashion from Models $(1)$-$(4)$ and the baseline Gaussian model. 
We let $\beta \in \mathbb{R}^d$ be a sparse vector with $s=5$ signals, all of the same strength and positioned randomly in the $d$-length vector. Each signal is assigned a positive sign with probability $0.5$.
We fix $\sigma^2=1$, $\rho^2=1$, and vary $\beta$ such that the signal-to-noise ratio $\text{snr}=\frac{1}{\sigma^2}\beta' \Sigma \beta$ takes values in the set
$$ \{0.10, 0.15, 0.20\}.$$

In this example, the function $\rF$ and our pivot no longer have a closed form expression.
To alleviate this computational barrier, we use a Laplace-type probabilistic approximation proposed by \cite{panigrahi2017mcmc} to compute $\rF$. 
Inverting the approximate pivot yields asymptotic confidence intervals based on our carved pivot.
% in terms of averaged coverage across the selected coefficients and the averaged lengths of the $90\%$ confidence intervals.
The cells in Table \ref{table:eg3:0} compare the $90\%$-confidence intervals based on carving and data splitting. 
We note that our asymptotic intervals not only cover the selected regression parameters at the desired level, but also provide tighter bounds than data splitting.
Furthermore, selective inference is valid even at lower values of signal-to-noise ratio, where rare outcomes are more likely.
\begin{center}
\captionof{table}{\centering{Comparison of inference between carving and data splitting.}}
\vspace{-0.2cm}
\label{table:eg3:0}
\bgroup
\def\arraystretch{1.3}
\scalebox{0.88}{\begin{tabular}{ |c | c c |c c |c c|c c |c c | }
\hline
$\rho^2=1$  & Cov & Len  & Cov  & Len  & Cov  & Len &  Cov  & Len  &  Cov  & Len \\
\hline
$\text{snr} =0.10$  & \multicolumn{2}{c}{Gaussian} & \multicolumn{2}{c}{Model-1} & \multicolumn{2}{c}{Model-2} & \multicolumn{2}{c}{Model-3} & \multicolumn{2}{c|}{Model-4} \\ [0.5ex] 
 \hline\hline
Carve  & $88.15\%$  & $0.43$ & $88.05\%$ & $0.42$ & $89.59\%$ & $0.42$ & $90.39\%$ & $0.42$ & $90.73\%$ & $0.43$  \\ [0.5ex] 
 \hline
Split  & $88.56\%$  & $0.52$ & $88.15\%$ & $0.50$ & $90.24\%$ & $0.50$ & $88.75\%$ & $0.50$ & $88.75\%$ & $0.50$\\ [0.5ex] 
 \hline\hline
 $\text{snr}= 0.15$  & \multicolumn{2}{c}{Gaussian} & \multicolumn{2}{c}{Model-1} & \multicolumn{2}{c}{Model-2} & \multicolumn{2}{c}{Model-3} & \multicolumn{2}{c|}{Model-4} \\
 \hline\hline
Carve  & $88.61\%$  & $0.43$ & $88.28\%$ & $0.42$ & $88.56\%$ & $0.42$ & $89.06\%$ & $0.44$ & $87.69\%$ & $0.42$  \\ [0.5ex] 
 \hline
Split  & $89.27\%$  & $0.50$ & $90.26\%$ & $0.51$ & $89.78\%$ & $0.50$ & $89.00\%$ & $0.53$ & $85.51\%$ & $0.50$\\ [0.5ex] 
 \hline\hline
$\text{snr}=0.20$  & \multicolumn{2}{c}{Gaussian} & \multicolumn{2}{c}{Model-1} & \multicolumn{2}{c}{Model-2} & \multicolumn{2}{c}{Model-3} & \multicolumn{2}{c|}{Model-4} \\
 \hline\hline
Carve  & $91.22\%$  & $0.43$ & $88.53\%$ & $0.42$ & $89.18\%$ & $0.43$ & $92.88\%$ & $0.43$ & $89.94\%$ & $0.42$  \\ [0.5ex] 
 \hline
Split  & $88.65\%$  & $0.51$ & $88.15\%$ & $0.50$ & $90.54\%$ & $0.52$ & $88.72\%$ & $0.51$ & $89.97\%$ & $0.51$\\ [0.5ex] 
 \hline
\end{tabular}}
\egroup
\end{center}

\noindent\textbf{Real data}.\ \ We apply our carved pivot on real data. Our data comes from $441$ patients in the publicly available The Cancer Genome Atlas (TCGA) database  \cite{tomczak2015review}. Carving is applied to infer for the selected associations between gene expression values and log-transformed survival times for Gliomas, a common type of brain tumor.
We include $2500$ predictors with the highest variability in the observed samples and solve the LASSO on a randomly drawn subsample of the full data.
The $\ell_1$ penalty tuning parameter is fixed at a theoretical value that was sugggested by \cite{negahban2009unified}. 

We obtain confidence intervals for the selected regression coefficients by inverting the carved pivot. 
Figure \ref{len:TCGA} shows the distribution of lengths of the confidence intervals based on carving and data splitting. 
On the x-axis, we vary the ratio $1/(1+\rho^2)$. 
The plot demonstrates the advantages of conducting selective inference with the carved pivot, which re-uses data from selection steps. 
Interval estimates for both procedures grow wider when fewer holdout samples are available for inference. 
However, the benefits of carving only become more pronounced as more data is used at the selection step.

\begin{figure}[h]
\centering
\includegraphics[scale=0.65]{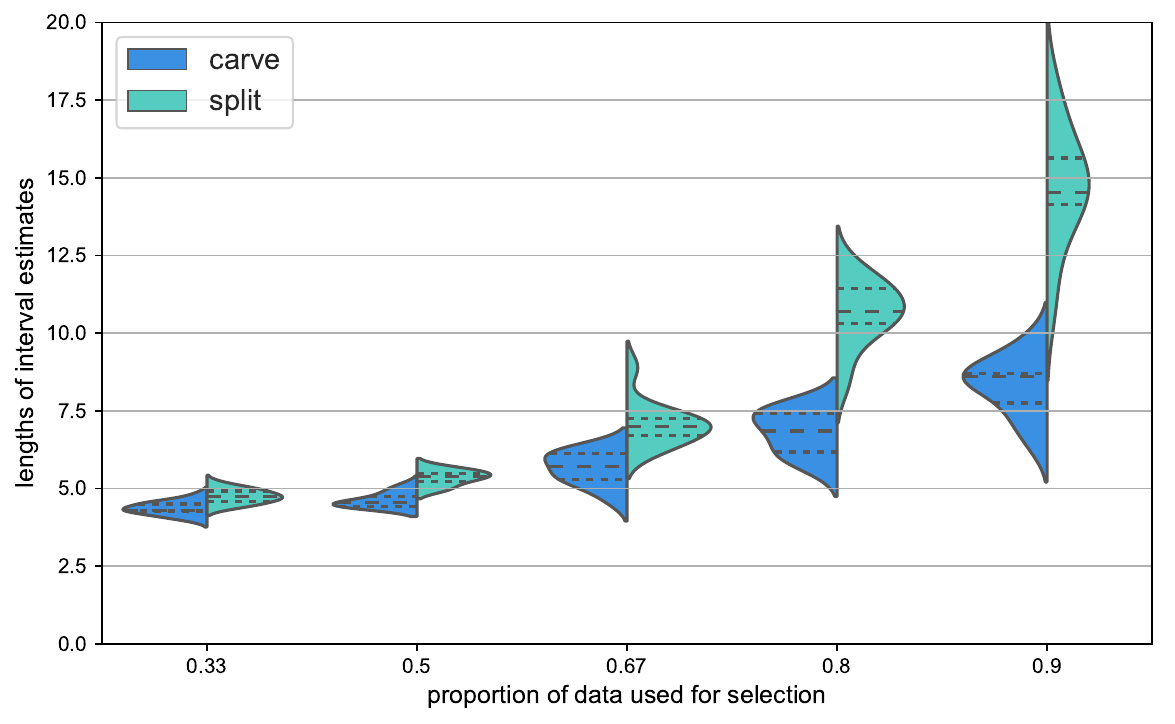}
\caption{Distribution of lengths of interval estimates for the selected regression coefficients. }
\label{len:TCGA}
\end{figure}
}
\end{example}

\section{Conclusion}
\label{conclusion}
Our paper provides an asymptotic basis for carving as we depart from Gaussian data.
Our setup considers two datasets: one of the datasets is used for selection and the other dataset is reserved for inference.
As an example, this setup is commonly encountered when the investigator selects promising findings on pilot data.
Inference for the selected findings is a natural goal when new data arrives at a later stage of the experiment.
Carving not only adjusts for overoptimism resulting from selection, but also re-uses pilot data for efficient inference. 
We show that pivots formed by conditioning on the selection outcome in the pilot data yield asymptotically-valid inference.
More generally, our theory subtantiates the use of pivots based on Gaussian randomized selection rules. 
Recent work by multiple papers, e.g., \cite{rasines2021splitting, schultheiss2021multicarving, panigrahi2022treatment, panigrahi2022integrative}, have explored the potential of randomized selection rules for improved inference, in theory and various applications. 

While in this paper we have focused on pivots based on the conditional method in \cite{exact_lasso}, for future work we will consider other types of pivot that have been developed for conditional inference.  
For example, \cite{panigrahi2022approximate} propose an approximate Gaussian pivot using the maximum likelihood estimator and \cite{liu2018more} propose pivots in the full model with strictly less conditioning than \cite{exact_lasso}. 
New theoretical results are required to study the rate of weak convergence for such pivots and investigate if asymptotically-valid selective inference continues to hold if the pivots were formed with self-normalized statistics.

\section{Acknowledgements}
S.P. acknowledges support in part by NSF grants DMS 1951980 and DMS 2113342.
S.P. would like to thank Jonathan Taylor, Liza Levina and Xuming He for their generous help in early stages of this paper. 
S.P. would also like to thank two anonymous referees for providing several insightful comments on an initial draft of the paper.

\bibliographystyle{imsart-number}
\bibliography{references.bib}

\appendix
\section*{Proofs and supporting results}
\addcontentsline{toc}{section}{Appendices}
\renewcommand{\thesubsection}{\Alph{subsection}}

\subsection{Proofs for Section 2}

\begin{proof} [Proof of Proposition 1]
First, we show that the density for the conditional distribution 
$$V_n^{\;(j)}\; \lvert \; U_n^{\;(j)}=\oU,\; E_n= \oE, \; A_n= \oA,$$
at $v$, is given by
\begin{equation}
\label{cond:law}
(\rD(\oU; \sqrt{n}\beta_n^{\;(j)}))^{-1}\cdot \phi\left(\frac{1}{\sigma_j}(v-\sqrt{n}\beta_n^{\;(j)})\right) \cdot  \rF\left(R^{\;(j)}\begin{pmatrix} v & (\oU)' \end{pmatrix}'\right).
\end{equation}
To do so, 
%observe that
%$$\begin{pmatrix}{V_n^{\;(\oE)}}' & {V_n^{\;(\oE^c)}}' \end{pmatrix}' = R^{\;(j)}\begin{pmatrix} (V_n^{\;(j)}) & (U_n^{\;(j)})' \end{pmatrix}',$$
%and that $V_n^{\;(j)}$ is independent of $U_n^{\;(j)}$.
%For $\mu \in \mathbb{R}^d$, $\Theta \in \mathbb{R}^{d\times d}$, we use $\rp(x; \mu, \Theta)$ to represent the density function of the $d$-dimensional Gaussian variable with mean $\mu$ and covariance $\Theta$ at $x$.
let us denote the mean vector and covariance matrix of $U_n^{\;(j)}$ by $\sqrt{n}\gamma_n^{\;(j)}$ and $\Sigma^{\;(j)}$, respectively. 
Observe, the joint likelihood for the variables $V_n^{\;(j)}$, $U_n^{\;(j)}$ and $W_n$ is proportional to
$$\phi\left(\dfrac{1}{\sigma_j}(V_n^{\;(j)}- \sqrt{n}\beta_n^{\;(j)})\right) \cdot \left\{  \tE(U_n^{\;(j)}-\sqrt{n}\gamma_n^{\;(j)}, (\Sigma^{\;(j)})^{-1}) \cdot \tE\left(W_n, \frac{1}{\rho^2}\Sigma^{-1}\right)\right\},$$
where we have used the independence between the variables $V_n^{\;(j)}$, $U_n^{\;(j)}$ and $W_n$.
Define a change of variables
$$\begin{pmatrix} V_n^{\;(j)}\\ U_n^{\;(j)}\\ W_n \end{pmatrix}\longrightarrow \begin{pmatrix} V_n^{\;(j)}\\ U_n^{\;(j)}\\ T_n \\ A_n  \end{pmatrix}$$
%where, for a fixed value of $V_n$, $T_n \in \real^p$ and $A_n \in \real^{d-p}$ are defined through the mapping $\pi_{V_n}$ as
where
\begin{equation}
\label{linear:map:sel}
\begin{aligned}
\begin{pmatrix}{W_n^{\;(\oE)}}'  &  {W_n^{\;(\oE^c)}}' \end{pmatrix}'  &= Q T_n + \begin{pmatrix} {\Lambda^{(\oE)}}' & {A_n}'\end{pmatrix}' - \begin{pmatrix}{V_n^{\;(\oE)}}' & {V_n^{\;(\oE^c)}}' \end{pmatrix}'\\
&= \pi_{V_n}(T_n, A_n),
\end{aligned}
\end{equation}
or equivalently, 
$$ 
 \begin{pmatrix}T_n' & A_n' \end{pmatrix}' = \pi_{V_n}^{-1}\left(\begin{pmatrix}{W_n^{\;(\oE)}}'  &  {W_n^{\;(\oE^c)}}' \end{pmatrix}'\right).
$$

We apply the above-stated change of variables to obtain a likelihood based on the density of the new variables $V_n^{\;(j)}$, $U_n^{\;(j)}$, $T_n$ and $A_n$.
It is easy to see that this likelihood is proportional to
$$\phi\left(\dfrac{1}{\sigma_j}(V_n^{\;(j)}- \sqrt{n}\beta_n^{\;(j)})\right)\cdot \left\{ \tE(U_n^{\;(j)}-\sqrt{n}\gamma_n^{\;(j)}, (\Sigma^{\;(j)})^{-1}) \cdot \tE\left(\pi_{V_n}(T_n, A_n), \frac{1}{\rho^2}\Sigma^{-1}\right) \right\},$$
after ignoring constants.
%The Jacobian from the change of variables also dissolves as a constant.

%For $V_n = R^{\;(j)}\begin{pmatrix} (V_n^{\;(j)}) & (U_n^{\;(j)})' \end{pmatrix}'$, applying this change of variables gives us the likelihood for the variables $V_n^{\;(j)}$, $U_n^{\;(j)}$, $T_n$ and $A_n$.
%Ignoring constants that include the Jacobian, this likelihood is proportional to
%$$\phi\left(\dfrac{1}{\sigma_j}(V_n^{\;(j)}- \sqrt{n}\beta_n^{\;(j)})\right)\cdot \left\{ \tE(U_n^{\;(j)}-\sqrt{n}\gamma_n^{\;(j)}, \Sigma^{-1}_j) \cdot \tE\left(\pi_{V_n}(T_n, A_n), \frac{1}{\rho^2}\Sigma^{-1}\right) \right\},$$
%after ignoring constants.

Note that the selection outcome in Equation 2.3 is equivalent to
$$\left\{T_n\in \real^{p+}, A_n= \oA\right\},$$
and that $$\begin{pmatrix}{V_n^{\;(\oE)}}' & {V_n^{\;(\oE^c)}}' \end{pmatrix}' = R^{\;(j)}\begin{pmatrix} (V_n^{\;(j)}) & (U_n^{\;(j)})' \end{pmatrix}'.$$
Thus, the distribution of $\begin{pmatrix} V_n^{\;(j)} & T'_n \end{pmatrix}'$ when conditioned on $U_n^{\;(j)}=\oU$ and the selection outcome
$$\Big\{E_n= \oE, \; A_n= \oA\Big\}$$ 
has density equal to
\begin{equation*}
\begin{aligned}
&\left(\rD(\oU; \sqrt{n}\beta_n^{\;(j)})\right)^{-1}\cdot  \phi\left(\dfrac{1}{\sigma_j}(v- \sqrt{n}\beta_n^{\;(j)})\right)\\
&\;\;\;\;\;\;\;\;\;\;\;\;\;\;\;\;\;\;\;\;\;\;\;\;\;\;\;\times \tE\left(Qt -R^{\;(j)}\begin{pmatrix} v' & {\oU}' \end{pmatrix}' + r, \frac{1}{\rho^2}\Sigma^{-1}\right) \cdot \mathbf{1}_{t\in \real^{p+}}
\end{aligned}
\end{equation*}
at $\begin{pmatrix}  v & t' \end{pmatrix}'$.
Integrating out $t$ in this joint density yields the conditional density in \eqref{cond:law}.
Our pivot is finally obtained by applying a probability integral transformation, which gives a $\text{Unif}\;(0,1)$ random variable.
\end{proof}

\begin{proof}[Proof of Corollary 1]
Note that $U_n^{\;(j)}= V_n^{\;(-j)}$ when $\Sigma=I_{d,d}$. 
Also, observe that 
 $$\tE\left(Qt -R^{\;(j)}\begin{pmatrix} v & {U_n^{\;(j)}}' \end{pmatrix}' + r, \frac{1}{\rho^2} \Sigma^{-1}\right)\propto \tE\left(t^{\;(j)} -v, \frac{1}{\rho^2}\right)\cdot \rL(V_n^{\;(-j)}, t^{\;(-j)}),$$
 where $\rL(V_n^{\;(-j)}, t^{\;(-j)})$ is a function of $V_n^{\;(-j)}$ and $t^{\;(-j)}$.
In particular, $\rL(V_n^{\;(-j)}, t^{\;(-j)})$ does not involve $v$ or $t^{\;(j)}$.
 Thus, our pivot simplifies as
 \begin{equation}
\label{pivot:pre}
 \dfrac{\displaystyle\int_{V_n^{\;(j)}}^{\infty} \phi\left(v- \sqrt{n}\beta_n^{\;(j)}\right)\cdot \bar{\Phi}\left(-\frac{1}{\rho}v\right)\;dv}{\displaystyle\int_{-\infty}^{\infty}\phi\left(v- \sqrt{n}\beta_n^{\;(j)}\right)\cdot \bar{\Phi}\left(-\frac{1}{\rho}v\right)\;dv }.
 \end{equation}
 After substituting variables
 $$\widetilde{v} = v-\sqrt{n}\beta_n^{\;(j)}$$
 in the two integrals of \eqref{pivot:pre}, we derive the simplified expression for our pivot.
\end{proof}

\subsection{Proofs for Section 3}
We provide proofs for our results under Section 3.

\begin{proof} [Proof of Proposition 2]
We follow the steps outlined in the proof of Proposition 1 to obtain the conditional likelihood for $V_n$.
It follows directly that the ratio of likelihood functions is equal to
$$\dfrac{\rF\left(V_{n; \text{obs}}\right)}{{\mathbb{E}_{\mathbb{P}_n} \left[ \rF\left(V_n\right)\right]}}.$$
Because $V_n = \Sigma^{1/2}\mathcal{Z}_n + \sqrt{n}\beta_n$, this ratio can be re-written in terms of the standardized variable as
$$\dfrac{\rF \left( \Sigma^{1/2}\oZ + \sqrt{n}\beta_n\right)}{{\mathbb{E}_{\mathbb{P}_n} \left[\rF\left(\Sigma^{1/2}\mathcal{Z}_n + \sqrt{n}\beta_n\right)\right]}}.$$
\end{proof}

\begin{proof} [Proof of Proposition 3] 
Observe that the expression on the left-hand side of the assertion is equal to
\begin{equation*}
\begin{aligned}
& \Big |\mathbb{E}_{\mathbb{P}_n}\left[\rH\circ \rP^{\;(j)}(\mathcal{Z}_n; \sqrt{n}\beta_n)\cdot{\text{LR}}_{\mathbb{P}_n}(\mathcal{Z}_n; \sqrt{n}\beta_n)\right]\\
&\;\;\;\;\;\;\;\;\;\;\;\;\;\;\;\;\;\;\;\;\;\;\;\;\;\;\;\;\;\;\;\;\;\;\;\;\;\;\;\;\;\;\;\;\;\;\;\;\;\;\;\;-\mathbb{E}_{\mathcal{N}}\left[\rH \circ \rP^{\;(j)}(\mathcal{Z}_n;\sqrt{n}\beta_n) \cdot{\text{LR}}_{\mathcal{N}}(\mathcal{Z}_n; \sqrt{n}\beta_n)\right]\Big|.
\end{aligned}
\end{equation*}
Through the triangle inequality, the difference can be bounded from above by
$$(\text{\bf T1}) + (\text{\bf T2}),$$ where
\begin{equation*}
\begin{aligned}
(\text{\bf T1}) &= \Big |\mathbb{E}_{\mathbb{P}_n}\left[\rH\circ \rP^{\;(j)}(\mathcal{Z}_n; \sqrt{n}\beta_n)\cdot {\text{LR}}_{\mathbb{P}_n}(\mathcal{Z}_n)\right]\\
&\;\;\;\;\;\;\;\;\;\;\;\;\;\;- \mathbb{E}_{\mathbb{P}_n}\left[\rH\circ \rP^{\;(j)}(\mathcal{Z}_n; \sqrt{n}\beta_n) \cdot{\text{LR}}_{\mathcal{N}}(\mathcal{Z}_n; \sqrt{n}\beta_n)\right]\Big|,
\end{aligned}
\end{equation*}
\begin{equation*}
\begin{aligned}
(\text{\bf T2})&= \Big |\mathbb{E}_{\mathbb{P}_n}\left[\rH \circ \rP^{\;(j)}(\mathcal{Z}_n;\sqrt{n}\beta_n)\cdot{\text{LR}}_{\mathcal{N}}(\mathcal{Z}_n; \sqrt{n}\beta_n)\right]\\
&\;\;\;\;\;\;\;\;\;\;\;\;\;\;- \mathbb{E}_{\mathcal{N}}\left[\rH \circ \rP^{\;(j)}(\mathcal{Z}_n;\sqrt{n}\beta_n)\cdot {\text{LR}}_{\mathcal{N}}(\mathcal{Z}_n; \sqrt{n}\beta_n)\right]\Big|.
\end{aligned}
\end{equation*}
Because $\sup |\rH| = \rK$, observe
\begin{align*}
(\text{\bf T1}) &\leq \displaystyle\int |\rH\circ \rP^{\;(j)}(z;\sqrt{n}\beta_n)|\cdot \Bigg |\dfrac{\rF \left( \Sigma^{1/2}z+ \sqrt{n}\beta_n\right)}{{\mathbb{E}_{\mathbb{P}_n} \left[\rF\left(\Sigma^{1/2}\mathcal{Z}_n + \sqrt{n}\beta_n\right)\right]}}  \\
&\;\;\;\;\;\;\;\;\;\;\;\;\;\;\;\;\;\;\;\;\;\;\;\;\;\;\;\;\;\;\;\;\;\;\;\;\;\;\;\;\;\;\;\;\;\;\;\;\;\;\;\;\;\;\;\;\;\;\;\;\;\;\;\;\;\;\;\;\;\;\;-\dfrac{\rF \left( \Sigma^{1/2}z + \sqrt{n}\beta_n\right)}{{\mathbb{E}_{\mathcal{N}} \left[\rF\left(\Sigma^{1/2}\mathcal{Z}_n + \sqrt{n}\beta_n\right)\right]}}\Bigg | \rd\mathbb{P}_n(z)\\
&\leq \rK \cdot \dfrac{ \Big|\mathbb{E}_{\mathbb{P}_n}\left[\rF \left( \Sigma^{1/2}\mathcal{Z}_n + \sqrt{n}\beta_n\right)\right] - \mathbb{E}_{\mathcal{N}}\left[\rF \left( \Sigma^{1/2}\mathcal{Z}_n + \sqrt{n}\beta_n\right)\right]\Big |}{\mathbb{E}_{\mathcal{N}}\left[\rF \left( \Sigma^{1/2}\mathcal{Z}_n + \sqrt{n}\beta_n\right)\right]}= \rK\cdot \rR^{(1)}_n .
\end{align*}
Next, we plug the value of ${\text{LR}}_{\mathcal{N}}(\mathcal{Z}_n; \sqrt{n}\beta_n)$ into the expression for (\text{\bf T2}).
We observe that the term in (\text{\bf T2}) assumes the following expression
\begin{align*}
\label{T:2}
&\Big(\mathbb{E}_{\mathcal{N}}\left[\rF \left( \Sigma^{1/2}\mathcal{Z}_n + \sqrt{n}\beta_n\right)\right]\Big)^{-1}\Bigg|\mathbb{E}_{\mathbb{P}_n}\left[\rH\circ \rP^{\;(j)}(\mathcal{Z}_n;\sqrt{n}\beta_n)\cdot\rF \left( \Sigma^{1/2}\mathcal{Z}_n + \sqrt{n}\beta_n\right)\right]\\
&\;\;\;\;\;\;\;\;\;\;\;\;\;\;\;\;\;\;\;\;\;\;\;\;\;\;\;\;\;\;\;\;\;\;\;\;\;\;\;\;\;\;\;\;\;\;\;\;\;\;\;\;\;\;\;\;\;\;\;\;\;\;\;\;\;\;\; - \mathbb{E}_{\mathcal{N}}\left[\rH\circ \rP^{\;(j)}(\mathcal{Z}_n; \sqrt{n}\beta_n)\cdot\rF \left( \Sigma^{1/2}\mathcal{Z}_n + \sqrt{n}\beta_n\right)\right] \Bigg|,
\end{align*}
which is equal to $\rR^{(2)}_n$.
This proves our claim.
%Note that $\rR^{(2)}_n$ is equal to
\end{proof}

\begin{proof} [Proof of Proposition 4]
This assertion follows by noticing that
$$\mathbb{E}_{\mathcal{N}} \left[\rF\left(\Sigma^{1/2}\mathcal{Z}_n + \sqrt{n}\beta_n\right)\right]$$ 
is equal to
\begin{align*}
& \int \frac{1}{(2\pi)^{d/2}} \tE\left(Qt -(\Sigma^{1/2}z + \sqrt{n}\beta_n) + r, \frac{1}{\rho^2}\Sigma^{-1}\right) \cdot \mathbf{1}_{t\in \real^{p+}} \cdot \tE(z, I_d) dzdt\\
&=C_0\cdot  \tE\left(\sqrt{n}\beta_n-r, \frac{1}{(1+\rho^2)} \cdot \Lambda \right)\\
&\;\;\;\;\;\;\;\;\;\;\;\;\;\;\;\;\;\;\;\;\;\times \int \frac{\text{det}((1+\rho^2)\bar\Sigma)^{-1/2}}{(2\pi)^{(p)/2}} \tE\left(t-\bar{\mu}_n, \frac{\bar{\Sigma}^{-1}}{(1+\rho^2)}\right)\cdot \mathbf{1}_{t\in \real^{p+}} dt\\
&=C_0\cdot  \tE\left(\sqrt{n}\beta_n-r, \frac{1}{(1+\rho^2)} \cdot \Lambda \right) \cdot \mathbb{P}_{\mathcal{N}}\left[ T_n >0_{p}\right],
 \end{align*} 
 where $C_0$ is a constant free of $n$.
\end{proof}

\subsection{Proofs for Section 4}
Throughout our proofs, uppercase letters with subscripts represent constants that do not depend on the sample size $n$.
We use the symbol $\mathcal{D}^{k}f(x_0)$ to denote the $k^{\text{th}}$ derivative of a differentiable function $f$ at $x_0$, and simply use $\mathcal{D}f(x_0)$ to denote the first derivative of $f$.
The symbol $Z$ is understood to be a standard Gaussian variable.

\subsubsection{Proofs for Lemma 1 and Proposition 5}
We first derive the Stein bound and properties of our univariate pivot before proceeding to proofs for our main results.

\begin{proof} [Proof of Lemma 1]
The Stein function $\cS_{\rg}$ satisfies the following equality 
\begin{equation*}
 \mathbb{E}_{\mathbb{P}_n}\left[\rg(\mathcal{Z}_n^{\;(j)})\right]- \mathbb{E}_{\mathcal{N}}\left[\rg(\mathcal{Z}_n^{\;(j)})\right]=\mathbb{E}_{\mathbb{P}_n}\left[\cD\cS_{\rg}(\mathcal{Z}_n^{\;(j)})\right]- \mathbb{E}_{\mathbb{P}_n}\left[\mathcal{Z}_n^{\;(j)} \cS_{\rg}(\mathcal{Z}_n^{\;(j)})\right];
 \label{St:eqn}
 \end{equation*}
see Equation 3 in \cite{chen2021stein}.
 Simplifying the two terms on the right-hand side expression, we note that
 $$\mathbb{E}_{\mathbb{P}_n}\left[\cD\cS_{\rg}(\mathcal{Z}_n^{\;(j)})\right] = \sum_{i=1}^n \int_{-\infty}^{\infty} \mathbb{E}_{\mathbb{P}_n}\left[\cD\cS_{\rg}(\mathcal{Z}_n^{\;(j)})\right]\cdot \mathrm{M}_i(t) dt,$$
 and
 \begin{equation*}
 \begin{aligned}
  \mathbb{E}_{\mathbb{P}_n}\left[\mathcal{Z}_n^{\;(j)} \cS_{\rg}(\mathcal{Z}_n^{\;(j)})\right] &= \sum_{i=1}^n  \mathbb{E}_{\mathbb{P}_n}\left[Z_{i,n}^{\;(j)}\int_{0}^{Z_{i,n}^{\;(j)}}\cD\cS_{\rg}(\mathcal{Z}_n^{\;(j)}[-i] + t)dt \right]\\
  &=\sum_{i=1}^n  \mathbb{E}_{\mathbb{P}_n}\Big[\int_{-\infty}^{\infty}\cD\cS_{\rg}(\mathcal{Z}_n^{\;(j)}[-i] + t) Z_{i,n}^{\;(j)} \Big(\mathbf{1}_{[t, \infty)}(Z_{i,n}^{\;(j)}) \mathbf{1}_{[0, \infty)}(t)\\
  &\;\;\;\;\;\;\;\;\;\;\;\;\;\;\;\;\;\;\;\;\;\;\;\;\;\;\;\;\;\;\;\;\;\;\;\;\;\;\;\;\;\;\;\;\;\;\;\;\;\;\;\;\;\;\;\;\;\;\;\;\;\;\;\;\;\;\;\;\;\;\;\;- \mathbf{1}_{(-\infty, t]}(Z_{i,n}^{\;(j)})\mathbf{1}_{(-\infty, 0)}(t) \Big)dt \Big]\\
  &= \sum_{i=1}^n  \int_{-\infty}^{\infty}\mathbb{E}_{\mathbb{P}_n}\left[\cD\cS_{\rg}(\mathcal{Z}_n^{\;(j)}[-i] + t) \right] \mathrm{M}_i(t)  dt.
 \end{aligned}
 \end{equation*}
Combining the two terms, we write
  \begin{equation*}
 \begin{aligned}
& \mathbb{E}_{\mathbb{P}_n}\left[\cD\cS_{\rg}(\mathcal{Z}_n^{\;(j)})- \mathcal{Z}_n^{\;(j)} \cS_{\rg}(\mathcal{Z}_n^{\;(j)})\right] \\
&= \sum_{i=1}^n \int_{-\infty}^{\infty} \mathbb{E}_{\mathbb{P}_n}\left[\cD\cS_{\rg}(\mathcal{Z}_n^{\;(j)})- \cD\cS_{\rg}(\mathcal{Z}_n^{\;(j)}[-i] + t)\right]\cdot \mathrm{M}_i(t) dt.
 \end{aligned}
 \end{equation*}
After applying a Taylor series expansion to the integrand in the previous display, we observe that
\begin{equation*}
\begin{aligned}
&\Big|\mathbb{E}_{\mathbb{P}_n}\left[\rg(\mathcal{Z}_n^{\;(j)})\right]- \mathbb{E}_{\mathcal{N}}\left[\rg(\mathcal{Z}_n^{\;(j)})\right]\Big|\\
&\leq \sum_{i=1}^n \displaystyle\int_{-\infty}^{\infty}  \sup_{\alpha \in [0, 1]}  \mathbb{E}_{\mathbb{P}_n}\Bigg[\left(|t| +\frac{1}{\sqrt{n}} |\mathrm{e}_{i,n}^{\;(j)}|\right) \\
&\;\;\;\;\;\;\;\;\;\;\;\;\;\;\;\;\;\;\;\;\;\;\;\;\;\;\;\;\;\;\;\;\;\;\;\;\;\;\;\;\;\;\;\;\;\;\;\;\;\;\;\;\;\;\;\;\;\;\;\;\;\;\;\;\;\times\Big|\cD^2\cS_{\rg}\left(\alpha t+(1-\alpha) \frac{1}{\sqrt{n}} \mathrm{e}_{i,n}^{\;(j)} + \mathcal{Z}_n^{\;(j)}[-i]\right)\Big|\Bigg] \mathrm{M}_i(t) dt\\
&=n\cdot \displaystyle\int_{-\infty}^{\infty}  \sup_{\alpha \in [0, 1]}  \mathbb{E}_{\mathbb{P}_n}\Bigg[\left(|t| +\frac{1}{\sqrt{n}} |\mathrm{e}_{1,n}^{\;(j)}|\right)\\
&\;\;\;\;\;\;\;\;\;\;\;\;\;\;\;\;\;\;\;\;\;\;\;\;\;\;\;\;\;\;\;\;\;\;\;\;\;\;\;\;\;\;\;\;\;\;\;\;\;\;\;\;\;\;\;\;\;\;\;\;\;\;\;\;\;\times \Big|\cD^2\cS_{\rg}\left(\alpha t+(1-\alpha) \frac{1}{\sqrt{n}} \mathrm{e}_{1,n}^{\;(j)}+ \mathcal{Z}_n^{\;(j)}[-1]\right)\Big|\Bigg] \mathrm{M}_1(t) dt.
\end{aligned}
\end{equation*}
The last equality is based on the fact that $\left\{\mathrm{e}_{i,n}^{\;(j)}, i\in [n]\right\}$ is a collection of identically distributed observations.
\end{proof}

\begin{proof}[Proof of Proposition 5]
To prove this result, we obtain an alternate representation of our pivot.
Note that our pivot is equal to
\begin{equation*}
\Big(\mathbb{P}_{\mathcal{N}}[T_n >0]\Big)^{-1}\; {\mathbb{E}_{\mathcal{N}}\Big[\mathbb{E}_{\mathcal{N}}\left[\mathbf{1}_{\{\mathcal{Z}_n^{\;(j)} >\oZ^{\;(j)}\}} \cdot \mathbf{1}_{\{T_n\in (0,\infty)\}} \;\lvert \; T_n = t\right]\Big]},
\end{equation*}
where 
%Consider variables $\mathcal{Z}_n^{\;(j)}$ and  $T_n$. 
%Suppose that the density of these variables at $z$, $t$ is proportional to
%\[\tE(z, 1) \cdot \tE\left(t-(z+\sqrt{n}\beta_n^{\;(j)}), \frac{1}{\rho^2}\right).\]
%Based on the distribution described above, we have
$$T_n\sim \mathcal{N}(\sqrt{n}\beta^{\;(j)}_n, 1+\rho^2),\; \mathcal{Z}_n^{\;(j)} \lvert T_n=t \sim \mathcal{N}\left(\; \frac{1}{(1+\rho^2)}(t-\sqrt{n}\beta_n^{\;(j)}), \frac{\rho^2}{(1+\rho^2)}\;\right).$$
Using the above-stated representation, our pivot can be re-written as
\begin{equation*}
\begin{aligned}
%&\Big(\mathbb{P}_{\mathcal{N}}[T_n >0]\Big)^{-1}\; {\mathbb{E}_{\mathcal{N}}\Big[\mathbb{P}_{\mathcal{N}}\left[\mathcal{Z}_n^{\;(j)} >\oZ^{\;(j)}  \;\lvert \; T_n = t\right]\; \cdot \mathbf{1}_{t\in (0,\infty)}\Big]}\\
&\dfrac{\bigints \bar{\Phi}\left(\frac{1}{\rho}\sqrt{1+\rho^2} \oZ^{\;(j)}  -\frac{1}{\rho \sqrt{1+\rho^2}}(t-\sqrt{n}\beta_n^{\;(j)})\right)\cdot \tE\left(t-\sqrt{n}\beta^{\;(j)}_n,\frac{1}{1+\rho^2}\right) \cdot \mathbf{1}_{t\in (0, \infty)}\; dt}{\bigints\tE\left(t-\sqrt{n}\beta^{\;(j)}_n,\frac{1}{1+\rho^2}\right) \cdot \mathbf{1}_{t \in (0, \infty)}\;dt}.
\end{aligned}
\end{equation*}
It is straightforward to see that the first derivative of the pivot is bounded above by a constant.
\end{proof}

\subsubsection{Supporting Results}
\label{supp:univariate}

We state some supporting results to prove the weak convergence of the univariate pivot.
%Our supporting results s for the Stein function in Section 4.2. 
%Consider $f:\mathcal{C}\to \real$ for an open set $\mathcal{C}\subset \real^d$.
%Then, $\mathcal{D}^k f(x_0)[i_1, i_2,\cdots, i_k]$ is the $k^{\text{th}}$ order partial derivative of $f$ with respect to $x^{(i_1)}$, $\cdots$, $x^{(i_k)}$ at $x_0\in \mathcal{C}$. 

\begin{lemma}
\label{fact:2}
Let $\widetilde\rG_l$ be defined according to Equation (3.7). 
Denote by $\mathcal{M}$ the collection of all real-valued sequences $\{\beta_n^{\;(j)}: n\in \mathbb{N}\}$.
For $l\in [2]$, it holds that
$$\sup_{\mathcal{M}} \sup |\cD^2 \cS_{\widetilde\rG_l}|<\infty.$$
\end{lemma}

\begin{proof}[Proof of Lemma \ref{fact:2}]
We note that 
\[\cD^2\cS_{\widetilde\rG_l}(z) = (1+ z^2) \cS_{\widetilde\rG_l}(z) + z\cdot (\widetilde\rG_l(z; \sqrt{n}\beta_n^{\;(j)}) -\mathbb{E}_{\mathcal{N}}[\widetilde\rG_l(Z; \sqrt{n}\beta_n^{\;(j)})]) + \cD\widetilde\rG_l(z; \sqrt{n}\beta_n^{\;(j)}),\]
where 
\begin{equation*}
\begin{aligned}
&\widetilde\rG_l(z; \sqrt{n}\beta_n^{\;(j)}) -\mathbb{E}_{\mathcal{N}}[\widetilde\rG_l(Z; \sqrt{n}\beta_n^{\;(j)})] \\
&\;\;\;\;\;\;\;\;\;\;\;\;\;\;\;\;\;\;\;\;\;\;\;\;\;\;\;\;\;\;\;\;\;\;\;= \int_{-\infty}^{z} \cD\widetilde\rG_l(t; \sqrt{n}\beta_n^{\;(j)}) \Phi(t)dt -  \int_{z}^{\infty} \cD\widetilde\rG_l(t; \sqrt{n}\beta_n^{\;(j)}) \bar{\Phi}(t)dt.
\end{aligned}
\end{equation*}
Also, observe that 
\begin{equation*}
\begin{aligned}
\cS_{\widetilde\rG_l}(z) = -\sqrt{2\pi} \exp\left(\frac{1}{2}z^2\right) & \Big(\Phi(z) \cdot \int_{z}^{\infty}\cD\widetilde\rG_l(t; \sqrt{n}\beta_n^{\;(j)})\bar{\Phi}(t)dt  \\
&+ \bar{\Phi}(z) \cdot \int_{-\infty}^{z} \cD\widetilde\rG_l(t; \sqrt{n}\beta_n^{\;(j)})\Phi(t)dt\Big).
\end{aligned}
\end{equation*}
Therefore, we have 
\begin{equation*}
\begin{aligned}
 |\cD^2\cS_{\widetilde\rG_l}(z)| &\leq |(1+z^2)\cS_{\widetilde\rG_l}(z) + z(\widetilde\rG_l(z; \sqrt{n}\beta_n^{\;(j)}) -\mathbb{E}[\widetilde\rG_l(Z; \sqrt{n}\beta_n^{\;(j)})]| + |\cD\widetilde\rG_l(z; \sqrt{n}\beta_n^{\;(j)})| \\
&\leq \Big|(z-\sqrt{2\pi}(1+z^2)\exp\left(\frac{1}{2}z^2\right)\bar{\Phi}(z))\cdot \int_{-\infty}^{z} \cD\widetilde\rG_l(t; \sqrt{n}\beta_n^{\;(j)}) \Phi(t)dt \Big|\\
&\;\;+ \Big|(-z-\sqrt{2\pi}(1+z^2)\exp\left(\frac{1}{2}z^2\right)\Phi(z)) \cdot \int_{z}^{\infty} \cD\widetilde\rG_l(t; \sqrt{n}\beta_n^{\;(j)}) \bar{\Phi}(t)dt \Big|\\
&\;\;+ |\cD\widetilde\rG_l(z; \sqrt{n}\beta_n^{\;(j)})|\\
& \leq 2 \cdot \displaystyle\sup_z |\cD\widetilde\rG_l(z; \sqrt{n}\beta_n^{\;(j)})|.
\end{aligned}
\end{equation*}
At last, observe that: (i) $\rH$ is uniformly bounded and has a uniformly bounded derivative, and (ii) $\rP^{\;(j)}(\cdot; \sqrt{n}\beta^{\;(j)}_n)$ has a uniformly bounded derivative and the bound is uniform over $\mathcal{M}$. 
This proves our claim.
\end{proof}

\begin{lemma}
\label{bar:G}
Define the set $\mathcal{S}_n = \left[-c a_n \bar{\beta}, ca_n \bar{\beta}\right]$, where $c\in (0,1)$, and define the functions
$\rg_1(z; \sqrt{n}\beta^{\;(j)}_n) = 1$ and $\rg_2(z; \sqrt{n}\beta^{\;(j)}_n) = \rH\circ \rP^{\;(j)}(z\; ;\sqrt{n}\beta^{\;(j)}_n)$.
Let $\mathcal{M}_r= \{\beta_n^{(j)}: \ n\in \mathbb{N}\}$ be the collection of parameters that are parameterized as per Equation (4.2).
We note that there exists a differentiable, real-valued function $\mathcal{A}_n(\cdot; \sqrt{n}\beta_n^{\;(j)}):\real \to \real$ such that
$$\bar{\rG}_l(z; \sqrt{n}\beta_n^{\;(j)}) = \rg_l(z; \sqrt{n}\beta^{\;(j)}_n) \cdot  {\normalfont\text{Exp}}\left(z+ \sqrt{n}\beta_n^{\;(j)}, \frac{1}{\rho^2}\right) \cdot \mathcal{A}_n(z; \sqrt{n}\beta_n^{\;(j)})\cdot \mathbf{1}_{z\in \mathcal{S}_n}$$
agrees with $\widetilde\rG_l(z; \sqrt{n}\beta_n^{\;(j)})$ on its support set.
Additionally, this function satisfies
\begin{equation*}
\begin{aligned}
& \displaystyle\sup_{\mathcal{M}_r} \displaystyle\sup_{z\in \mathcal{S}_n} |\sqrt{n}\beta_n^{\;(j)}| \cdot |\mathcal{A}_n(z; \sqrt{n}\beta_n^{\;(j)})|<\infty,\\
& \displaystyle\sup_{\mathcal{M}_r} \displaystyle\sup_{z\in \mathcal{S}_n} |\sqrt{n}\beta_n^{\;(j)}| \cdot  |\mathcal{D}\mathcal{A}_n(z; \sqrt{n}\beta_n^{\;(j)})|<\infty.
\end{aligned}
\end{equation*}
\end{lemma}

\begin{proof}
This proof is immediate once we define
$$\mathcal{A}_n(z; \sqrt{n}\beta_n^{\;(j)}) =\dfrac{1}{\sqrt{2\pi}\rho} \int_{s>0} \exp\left(-\frac{1}{2\rho^2}s^2 + \frac{1}{\rho^2}s(z +\sqrt{n}\beta_n^{\;(j)})\right) ds.$$
\end{proof}

\begin{lemma}
\label{growth:Stein:function}
Let $\bar{\rG}^{+}_l(\cdot; \sqrt{n}\beta_n^{\;(j)})= \max(\bar{\rG}_l(\cdot; \sqrt{n}\beta_n^{\;(j)}), 0).$
Then, it holds that
\[\sup |\cS_{\bar{\rG}^{+}_l}| \leq \frac{C_0}{\sqrt{n}|\beta_n^{\;(j)}|}\cdot\text{\normalfont Exp}\left(\sqrt{n}\beta_n^{\;(j)}, \frac{1}{1+\rho^2}\right),\]
where $C_0$ is a constant that does not depend on $n$.
\end{lemma}

\begin{proof}[Proof of Lemma \ref{growth:Stein:function}]
The Stein function in Equation (4.3) satisfies 
\[\cS_{\bar{\rG}^{+}_l}(z)\leq \begin{dcases} 
       \exp\left(\frac{1}{2}z^2\right)\bigintsss_{-\infty}^{z} \bar{\rG}^{+}_l(t; \sqrt{n}\beta_n^{\;(j)})\tE(t, 1)dt & \text{ if } z< 0, \\ 
      \sqrt{2\pi} \cdot \exp\left(\frac{1}{2}z^2\right)\mathbb{E}_{\mathcal{N}}\left[\bar{\rG}^{+}_l(Z; \sqrt{n}\beta_n^{\;(j)})\right]\bar{\Phi}(z) & \text{ if } z \geq 0. \;
   \end{dcases}\]
Observe, $$\int_{-\infty}^{z} \bar{\rG}^{+}_l(t; \sqrt{n}\beta_n^{\;(j)})\tE(t, 1)dt$$ is bounded from above by 
\[ \frac{C_1}{\sqrt{n}|\beta_n^{\;(j)}|} \cdot \tE\left(\sqrt{n}\beta_n^{\;(j)}, \frac{1}{1+\rho^2}\right)\cdot \bar{\Phi}\left(-\frac{1}{\rho}{\sqrt{1+\rho^2} \cdot z} - \frac{1}{{\rho\sqrt{1+\rho^2}}}{\sqrt{n}\beta_n^{\;(j)}}\right).\]
Also, note that
\[\mathbb{E}_{\mathcal{N}}\left[\bar{\rG}^{+}_l(Z; \sqrt{n}\beta_n^{\;(j)})\right] \leq  \frac{C_2}{\sqrt{n}|\beta_n^{\;(j)}|}\cdot \tE\left(\sqrt{n}\beta_n^{\;(j)}, \frac{1}{1+\rho^2}\right).\]
Then, it follows that
\begin{equation*}
\begin{aligned}
 \cS_{\bar{\rG}^{+}_l}(z) &\leq C_3\cdot\Bigg( \exp\left(\frac{1}{2}z^2\right)\cdot \bar{\Phi}\left(-\frac{1}{\rho}{\sqrt{1+\rho^2} \cdot z} - \frac{1}{{\rho\sqrt{1+\rho^2}}}{\sqrt{n}\beta_n^{\;(j)}}\right)\cdot \mathbf{1}_{(-\infty,0)}(z)\\
&\;\;\;\;\;\;\;\;\;\;\;\;\;\;\;\;\;\;+   \exp\left(\frac{1}{2}z^2\right)\cdot \bar\Phi(z) \cdot \mathbf{1}_{(0, \infty)}(z)\Bigg) \cdot \frac{1}{\sqrt{n}|\beta_n^{\;(j)}|}\cdot \tE\left(\sqrt{n}\beta_n^{\;(j)}, \frac{1}{1+\rho^2}\right)\\
&\leq \frac{C_0}{\sqrt{n}|\beta_n^{\;(j)}|}\cdot\text{\normalfont Exp}\left(\sqrt{n}\beta_n^{\;(j)}, \frac{1}{1+\rho^2}\right).
\end{aligned}
\end{equation*}
\end{proof}

\begin{lemma}
Under Assumptions 2 and 3, we have
$$\sup_n\sup_{\mathbb{P}_n\in \mathcal{P}_{r,n}} \left({\normalfont \tE}\left(\sqrt{n}\beta_n^{\;(j)}, \frac{1}{1+\rho^2}\right)\right)^{-1}\cdot\mathbb{E}_{\mathbb{P}_n}\left[\exp\left(-\frac{1}{2\rho^2}{(\mathcal{Z}_n^{\;(j)}[-1]+\sqrt{n}\beta_n^{\;(j)})^2}\right)\right]<\infty.$$
\label{fact:00}
\end{lemma}

\begin{proof}
Let $\mathcal{K}_1= [-c_0, c_0]$ where $c_0> \frac{(\bar\beta + 1)}{\sqrt{(1+\rho^2)}}$, and let $\Psi_1(z) =  \frac{1}{2\rho^2}(z-\bar{\beta})^2$.
Under the conditions in Assumptions 2 and 3, we first note that
\begin{equation*}
\begin{aligned}
& \sup_n\sup_{\mathbb{P}_n\in \mathcal{P}_{r,n}} \left(\tE\left(\sqrt{n}\beta_n^{\;(j)}, \frac{1}{1+\rho^2}\right)\right)^{-1}\\
 &\;\;\;\;\;\;\;\;\;\;\;\;\;\;\times\mathbb{E}_{\mathbb{P}_n}\left[\exp\left(-\frac{1}{2\rho^2}{(\mathcal{Z}_n^{\;(j)}[-1]+\sqrt{n}\beta_n^{\;(j)})^2}\right)\cdot \mathbf{1}_{\frac{1}{a_n}\mathcal{Z}_n^{\;(j)}[-1] \in \mathcal{K}_1}\right]\\
 &= \sup_n \sup_{\mathbb{P}_n\in \mathcal{P}_{r,n}} \left(\tE\left(a_n\bar{\beta}, \frac{1}{1+\rho^2}\right)\right)^{-1}\\
 &\;\;\;\;\;\;\;\;\;\;\;\;\;\;\times\mathbb{E}_{\mathbb{P}_n}\left[\exp\left(-a_n^2\Psi_{1}\left(\frac{1}{a_n}\mathcal{Z}_n^{\;(j)}[-1]\right)\right)\cdot \mathbf{1}_{\frac{1}{a_n}\mathcal{Z}_n^{\;(j)}[-1] \in \mathcal{K}_1}\right] <\infty.
  \end{aligned}
\end{equation*}
Suppose that $\mathcal{K}_0=\mathcal{K}_1^c$. 
Under the same assumptions, we observe that
\begin{equation*}
\begin{aligned}
& \sup_n\sup_{\mathbb{P}_n\in \mathcal{P}_{r,n}} \left(\tE\left(\sqrt{n}\beta_n^{\;(j)}, \frac{1}{1+\rho^2}\right)\right)^{-1}\\
 &\;\;\;\;\;\;\;\;\;\;\;\;\;\;\times\mathbb{E}_{\mathbb{P}_n}\left[\exp\left(-\frac{1}{2\rho^2}{(\mathcal{Z}_n^{\;(j)}[-1]+\sqrt{n}\beta_n^{\;(j)})^2}\right)\cdot \mathbf{1}_{\frac{1}{a_n}\mathcal{Z}_n^{\;(j)}[-1] \in \mathcal{K}_0}\right]\\
 %&=  \sup_n\sup_{\mathbb{P}_n\in \mathcal{P}_{r,n}} \left(\tE\left(\sqrt{n}\beta_n^{\;(j)}, \frac{1}{1+\rho^2}\right)\right)^{-1}\\
% &\;\;\;\;\;\;\;\;\;\;\;\;\;\;\;\;\;\;\;\;\;\;\;\;\times \mathbb{E}_{\mathbb{P}_n}\left[\exp\left(-a_n^2\Psi_1\left(\frac{1}{a_n}\mathcal{Z}_n^{\;(j)}[-1]\right)\right) \mathbf{1}_{\frac{1}{a_n}\mathcal{Z}_n^{\;(j)}[-1] \in \mathcal{R}^c}\right]  \\
 &\leq \sup_n \sup_{\mathbb{P}_n\in \mathcal{P}_{r,n}} \left(\tE\left(\sqrt{n}\beta_n^{\;(j)}, \frac{1}{1+\rho^2}\right)\right)^{-1}\cdot\mathbb{P}_{\mathbb{P}_n}\left[\mathbf{1}_{\frac{1}{a_n}\mathcal{Z}_n^{\;(j)}[-1] \in \mathcal{K}_0}\right]<\infty.
 \end{aligned}
\end{equation*}
Combining both observations proves our claim.
\end{proof}

\begin{theorem}
\label{fact:4}
Suppose that the conditions under Assumptions 2 and 3 are met.
Then, we have
$$\text{\normalfont SB}_{\mathbb{P}_n}(\bar{\rG}_l^{+})\leq \frac{E_{l,2}}{\sqrt{n}}\cdot \text{\normalfont Exp}\left(\sqrt{n}\beta_n^{\;(j)}, \frac{1}{1+\rho^2}\right).
%,\text{\normalfont SB}_{\mathbb{P}_n}(\bar{\rG}_l^{-})\leq \frac{E_{l,2}}{\sqrt{n}}\cdot \text{\normalfont Exp}\left(\sqrt{n}\beta_n^{\;(j)}, \frac{1}{1+\rho^2}\right).
$$
\end{theorem}

\begin{proof}[Proof of Theorem \ref{fact:4}]
%We will prove the claim for $\text{\normalfont SB}_{\mathbb{P}_n}(\bar{\rG}_l^{+})$, omitting the exactly similar proof for $\text{\normalfont SB}_{\mathbb{P}_n}(\bar{\rG}_l^{-})$.
We divide our proof into two main steps.
In Step 1, we will establish that $\text{\normalfont SB}_{\mathbb{P}_n}(\bar{\rG}^{+}_l)$ is bounded by
\begin{equation}
A_1\cdot \Bigg(\text{Bd}_n(\beta^{\;(j)}_n) + \frac{1}{\sqrt{n}}{\tE\left(\sqrt{n}\beta_n^{\;(j)}, \frac{1}{1+\rho^2}\right) } \Bigg)
\label{part:I}
\end{equation}
where 
\begin{align*}
& \text{Bd}_n(\beta^{\;(j)}_n)  =n \cdot \displaystyle\int_{-\infty}^{\infty} \sup_{\alpha\in [0,1]}\mathbb{E}_{\mathbb{P}_n}\Big[\left(|t|+\frac{1}{\sqrt{n}}|{\mathrm{e}}_{1, n}^{\;(j)}|\right)\cdot \tE\left(W_{t,\alpha}+ \sqrt{n}\beta_n^{\;(j)}, \frac{1}{\rho^2}\right) \\ 
&\;\;\;\;\;\;\;\;\;\;\;\;\;\;\;\;\;\;\;\; \;\;\;\;\;\;\;\;\;\;\;\;\;\;\;\;\;\;\;\; \;\;\;\;\;\;\;\;\;\;\;\;\;\;\;\;\;\;\;\; \;\;\;\;\;\;\;\;\;\;\;\;\;\;\;\;\;\;\;\; \;\;\;\;\;\;\;\;\;\;\;\;\;\;\;\;\;\;\;\; \times  \mathbf{1}_{\left\{W_{t,\alpha} \in \mathcal{S}_n\right\}}\Big]\mathrm{M}_1(t) dt,
\end{align*}
and
$$
W_{t,\alpha}= \alpha t+(1-\alpha) Z_{1,n}^{\;(j)} + \mathcal{Z}_n^{\;(j)}[-1].
$$
In Step 2, we will show
\begin{equation}
\text{Bd}_n(\beta^{\;(j)}_n) \leq  \frac{A_2}{\sqrt{n}}\cdot{\tE\left(\sqrt{n}\beta_n^{\;(j)}, \frac{1}{1+\rho^2}\right) }
\label{part:II}
\end{equation}
to finish our proof.
\smallskip

\textbf{Step 1}. 
Letting 
$$ 
\rB(z;\sqrt{n}\beta^{\;(j)}_n) = \tE\left(z+ \sqrt{n}\beta_n^{\;(j)}, \frac{1}{\rho^2}\right)  \cdot \mathcal{A}_n(z; \sqrt{n}\beta_n^{\;(j)}), 
$$
we note that
\begin{equation*}
\begin{aligned}
\cD^2\cS_{\bar{\rG}^{+}_l}(z) &\leq (1+ z^2) |\cS_{\bar{\rG}^{+}_l}(z)| +D_1| \sqrt{n}\beta_n^{\;(j)}|\cdot \mathbb{E}_{\mathcal{N}}\left[\rB(Z;\sqrt{n}\beta_n^{\;(j)})\cdot\mathbf{1}_{Z\in\mathcal{S}_n}\right]   \\
&\;\;\;\;\;\;\;\;\;\;\;\;\;\;\;\;\;\;\;\;\;\;\;\;\;\;\;\;\;\;\;\;\;\;\;\;\;\;\;\;+D_2\cdot \left(|\sqrt{n}\beta_n^{\;(j)}| \rB(z;\sqrt{n}\beta_n^{\;(j)})+\cD\rB(z;\sqrt{n}\beta_n^{\;(j)})\right),
\end{aligned}
\end{equation*}
for $z\in\mathcal{S}_n$.
Plugging this bound into the expression of the Stein bound (see Lemma 1), we observe that
\[\text{\normalfont SB}_{\mathbb{P}_n}(\bar{\rG}^{+}_l) \leq \text{Bd}_{1,n} + D_1\cdot \text{Bd}_{2, n} + D_2\cdot \text{Bd}_{3, n},\]
where
\begin{equation*}
\begin{aligned}
&\text{Bd}_{1, n} = n \cdot\displaystyle\int_{-\infty}^{\infty}\sup_{\alpha \in [0,1]} \mathbb{E}_{\mathbb{P}_n}\Big[(|t|+|Z_{1,n}^{\;(j)}|)\left(1+ W^2_{t,\alpha}\right)  |\cS_{\bar{\rG}_l}(W_{t,\alpha})|\Big]\mathrm{M}_1(t) dt; 
\end{aligned}
\end{equation*}
\begin{equation*}
\begin{aligned}
& \text{Bd}_{2, n} = n^{3/2}|\beta_n^{\;(j)}|\cdot\mathbb{E}_{\mathcal{N}}\left[\rB(Z;\sqrt{n}\beta_n^{\;(j)})\cdot \mathbf{1}_{Z\in\mathcal{S}_n}\right]\cdot\displaystyle\int_{-\infty}^{\infty} \mathbb{E}_{\mathbb{P}_n}\left[|t|+|Z^{\;(j)}_{1,n}|\right] \mathrm{M}_1(t) dt; 
\end{aligned}
\end{equation*}
\begin{equation*}
\begin{aligned}
& \text{Bd}_{3, n} = n\cdot \displaystyle\int_{-\infty}^{\infty} \sup_{\alpha\in [0,1]}\mathbb{E}_{\mathbb{P}_n}\Bigg[(|t|+|Z_{1,n}^{\;(j)}|)\cdot \tE\left(W_{t,\alpha}+ \sqrt{n}\beta_n^{\;(j)}, \frac{1}{\rho^2}\right)  \\
&\;\;\;\;\;\;\;\;\;\;\;\;\;\;\;\;\;\;\;\;\;\;\;\;\;\;\;\;\;\;\;\;\;\;\;\;\;\;\;\;\;\;\;\;\;\;\;\;\times\Big(\sqrt{n}|\beta_n^{\;(j)}|\mathcal{A}_n(W_{t,\alpha}) + |\cD\mathcal{A}_n(W_{t,\alpha})|\Big)\cdot \mathbf{1}_{\left\{W_{t,\alpha}\in \mathcal{S}_n\right\}}\Bigg]\mathrm{M}_1(t) dt.
\end{aligned}
\end{equation*}
\smallskip

\noindent The conclusion in Lemma \ref{growth:Stein:function} when combined with the condition in Assumption 2 allows us to write
\begin{equation*}
\begin{aligned}
&\text{Bd}_{1, n} \leq  \frac{\sqrt{n}}{|\beta_n^{\;(j)}|}{\tE\left(\sqrt{n}\beta_n^{\;(j)}, \frac{1}{1+\rho^2}\right)} \bigintsss_{-\infty}^{\infty} \mathbb{E}_{\mathbb{P}_n}\Big[(|t|+|Z^{\;(j)}_{1,n}|)\\
&\;\;\;\;\;\;\;\;\;\;\;\;\;\;\;\;\;\;\;\;\;\;\;\;\;\;\;\;\;\;\;\;\;\;\;\;\;\;\;\;\;\;\;\;\;\;\;\;\;\;\;\;\;\;\;\;\;\;\;\;\;\;\;\;\;\;\;\;\;\;\;\;\;\;\;\times \left(1+ (|t|+|Z^{\;(j)}_{1,n}| + \mathcal{Z}_n^{\;(j)}[-1])^2\right)\Big] \mathrm{M}_1(t) dt\\
&\;\;\;\;\;\;\;\;\;\;= \frac{D_3}{n|\beta_n^{\;(j)}|}\cdot{\tE\left(\sqrt{n}\beta_n^{\;(j)}, \frac{1}{1+\rho^2}\right) }.
\end{aligned}
\end{equation*}
Additionally, we observe that
$$\text{Bd}_{2, n} \leq  \frac{D_4}{\sqrt{n}} \tE\left(\sqrt{n}\beta_n^{\;(j)}, \frac{1}{1+\rho^2}\right).$$
At last, the properties of the function $\mathcal{A}_n(\cdot; \sqrt{n}\beta_n^{\;(j)})$ (see Lemma \ref{bar:G}) lead us to note that
\begin{equation*}
\begin{aligned}
\text{Bd}_{3,n} 
%&\leq D_5\cdot \max\left(\displaystyle\sup_{\mathcal{M}_r} \displaystyle\sup_{z\in \mathcal{S}_n} |\sqrt{n}\beta_n^{\;(j)}||\mathcal{A}_n(z; \sqrt{n}\beta_n^{\;(j)})|, \displaystyle\sup_{\mathcal{M}_r} \displaystyle\sup_{z\in \mathcal{S}_n} |\sqrt{n}\beta_n^{\;(j)}| |\cD \mathcal{A}_n(z; \sqrt{n}\beta_n^{\;(j)})|\right) \\
%&\;\;\;\times n\left(1+\frac{1}{|\sqrt{n}\beta_n^{\;(j)}|}\right)\cdot \displaystyle\int_{-\infty}^{\infty} \sup_{\alpha\in [0,1]}\mathbb{E}_{\mathbb{P}_n}\Big[(|t|+|Z_{1,n}^{\;(j)}|)\\
%&\;\;\;\;\;\;\;\;\;\;\;\;\;\;\;\;\;\;\;\;\;\;\;\;\;\;\;\;\;\;\;\;\;\;\;\;\;\;\;\;\;\;\;\;\;\;\;\;\;\;\;\;\;\;\times  \tE\left(W_{t,\alpha}+ \sqrt{n}\beta_n^{\;(j)}, \frac{1}{\rho^2}\right)\cdot \mathbf{1}_{\left\{W_{t,\alpha}\in \mathcal{S}_n\right\}}\Big]\mathrm{M}_1(t) dt\\
&\leq D_5 \cdot n\left(1+\frac{1}{|\sqrt{n}\beta_n^{\;(j)}|}\right)\displaystyle\int_{-\infty}^{\infty} \sup_{\alpha\in [0,1]}\mathbb{E}_{\mathbb{P}_n}\Big[(|t|+|Z_{1,n}^{\;(j)}|)\cdot \tE\left(W_{t,\alpha}+ \sqrt{n}\beta_n^{\;(j)}, \frac{1}{\rho^2}\right)\\
& \;\;\;\;\;\;\;\;\;\;\;\;\;  \;\;\;\;\;\;\;\;\;\;\;\;\;  \;\;\;\;\;\;\;\;\;\;\;\;\;  \;\;\;\;\;\;\;\;\;\;\;\;\;  \;\;\;\;\;\;\;\;\;\;\;\;\;  \;\;\;\;\;\;\;\;\;\;\;\;\;   \;\;\;\;\;\;\;\;\;\;\;\;\;  \;\;\;\;\;\;\;\;\;\;\;\;\;   \;\;\;\;\;\;\;\;\;\;\;\;\;    \times\mathbf{1}_{\left\{W_{t,\alpha}\in \mathcal{S}_n\right\}}\Big]\mathrm{M}_1(t) dt.
\end{aligned}
\end{equation*}
This gives the bound in \eqref{part:I}.

\medskip

\textbf{Step 2}. 
We start from the following bound
\begin{equation*}
\begin{aligned}
\text{Bd}_n(\beta^{\;(j)}_n) &\leq n\displaystyle\bigintsss  \mathbb{E}_{\mathbb{P}_n}\Big[(|t| + |Z_{1,n}^{\;(j)}|)\Big\{\exp\Big(\frac{1}{\rho^2}{(|t|+|Z_{1,n}^{\;(j)}|)(|t|+|Z_{1,n}^{\;(j)}| +(c+1)|\sqrt{n}\beta_n^{\;(j)}|)}\\
&\;\;\;\;\;\;\;\;\;\;\;\;\;\;\;\;\;\;\;\;\;\;\;\;\;\;\;\;\;\;\;\;\;\;\;\;\;\;\;\;\;\;\;\;\;\;\;\;\;\;\;\;\;\;\;\;\;\;\;\;\;\;\;\;\;-\frac{1}{2\rho^2}{(\mathcal{Z}_n^{\;(j)}[-1]+\sqrt{n}\beta_n^{\;(j)})^2}\Big)\Big\}\Big]\mathrm{M}_1(t)dt,
%&\leq  D_7\cdot n \tE\left(\sqrt{n}\beta_n^{\;(j)}, \frac{1}{1+\rho^2}\right)\cdot\displaystyle\int  \mathbb{E}_{\mathbb{P}_n}\Big[(|t| + |Z_{1,n}^{\;(j)}|)\\
%&\;\;\;\;\;\;\;\;\;\;\;\;\;\;\;\;\;\;\;\;\;\times\exp\Big(\frac{1}{\rho^2}(c+1)\cdot|{\sqrt{n}\beta_n^{\;(j)}|(|t| + |Z_{1,n}^{\;(j)}|)}+\frac{1}{\rho^2}{(|t| + |Z_{1,n}^{\;(j)}|)^2}\Big)\Big] \mathrm{M}_1(t)dt.
\end{aligned}
\label{step}
\end{equation*}
where we have used the fact that
$$\alpha t+(1-\alpha) Z^{\;(j)}_{i,n}\leq |t| + |Z^{\;(j)}_{i,n}|.$$
Note that $Z_{1,n}^{\;(j)}$ and $\mathcal{Z}_n^{\;(j)}[-1]$ are independent variables and that
$$\sup_n\sup_{\mathbb{P}_n\in \mathcal{P}_{r,n}} \left(\tE\left(\sqrt{n}\beta_n^{\;(j)}, \frac{1}{1+\rho^2}\right)\right)^{-1}\cdot\mathbb{E}_{\mathbb{P}_n}\left[\exp\left(-\frac{1}{2\rho^2}{(\mathcal{Z}_n^{\;(j)}[-1]+\sqrt{n}\beta_n^{\;(j)})^2}\right)\right]<\infty$$
(see Lemma \ref{fact:00}).
As a result, we can state that
\begin{equation*}
\begin{aligned}
\text{Bd}_n(\beta^{\;(j)}_n) &\leq  D_6\cdot n \cdot\tE\left(\sqrt{n}\beta_n^{\;(j)}, \frac{1}{1+\rho^2}\right)\cdot\displaystyle\int  \mathbb{E}_{\mathbb{P}_n}\Big[(|t| + |Z_{1,n}^{\;(j)}|)\\
&\;\;\;\;\;\;\;\;\;\;\;\;\;\;\;\times\exp\Big(\frac{1}{\rho^2}(c+1)\cdot|{\sqrt{n}\beta_n^{\;(j)}|(|t| + |Z_{1,n}^{\;(j)}|)}+\frac{1}{\rho^2}{(|t| + |Z_{1,n}^{\;(j)}|)^2}\Big)\Big] \mathrm{M}_1(t)dt.
\end{aligned}
\label{step}
\end{equation*}

Define
 $$
\textbf{(TA)}=\displaystyle\int \exp\left( \frac{1}{\rho^2}(c+1)\cdot|\sqrt{n}\beta_n^{\;(j)}| |t| +  \frac{2}{\rho^2}t^2\right)\mathrm{M}_1(t)dt,
$$
$$
\textbf{(TB)}= \displaystyle\int |t|\exp\left(\frac{1}{\rho^2}(c+1)\cdot|\sqrt{n}\beta_n^{\;(j)}||t| +  \frac{2}{\rho^2}t^2\right)\mathrm{M}_1(t)dt.
$$
%Noting that
%$$(|t| + |Z_{1,n}^{\;(j)}|)^2 \leq 2t^2 + 2(Z_{1,n}^{\;(j)})^2,$$
Note, the previous bound on $\text{Bd}_n(\beta^{\;(j)}_n)$ simplifies as 
\begin{equation*}
\begin{aligned}
\text{Bd}_n(\beta^{\;(j)}_n)  &\leq
 D_6\cdot n\cdot \tE\left(\sqrt{n}\beta_n^{\;(j)}, \frac{1}{1+\rho^2}\right)\cdot \mathbb{E}_{\mathbb{P}_n}\Big[|Z_{1,n}^{\;(j)}|\exp\Big(\frac{1}{\rho^2}(c+1)\cdot|{\sqrt{n}\beta_n^{\;(j)}||Z_{1,n}^{\;(j)}|}\\
&+  \frac{2}{\rho^2} (Z_{1,n}^{\;(j)})^2\Big)\Big]\times  \displaystyle\int \exp\left( \frac{1}{\rho^2}(c+1)\cdot|\sqrt{n}\beta_n^{\;(j)}| |t| +  \frac{2}{\rho^2}t^2\right)\mathrm{M}_1(t)dt\\
&+  D_6\cdot n \cdot\tE\left(\sqrt{n}\beta_n^{\;(j)}, \frac{1}{1+\rho^2}\right)\cdot \mathbb{E}_{\mathbb{P}_n}\Big[\exp\Big(\frac{1}{\rho^2}(c+1)\cdot |\sqrt{n}\beta_n^{\;(j)}| |Z_{1,n}^{\;(j)}| +  \frac{2}{\rho^2}(Z_{1,n}^{\;(j)})^2\Big)\Big]\\
&\times  \displaystyle\int |t|\exp\left(\frac{1}{\rho^2}(c+1)\cdot|\sqrt{n}\beta_n^{\;(j)}||t| +  \frac{2}{\rho^2}t^2\right)\mathrm{M}_1(t)dt\\
&\leq  D_7\cdot  n \cdot\tE\left(\sqrt{n}\beta_n^{\;(j)}, \frac{1}{1+\rho^2}\right)\cdot \left(\frac{1}{ \sqrt{n}}\textbf{(TA)} +  \textbf{(TB)}\right).
\end{aligned}
\label{bound:inter}
\end{equation*}

Left to analyze the terms $\textbf{(TA)}$, $\textbf{(TB)}$, we observe that the integrands in both terms are symmetric functions about $0$ and are increasing on the positive axis. 
Therefore, we simplify $\textbf{(TA)}$ as follows:
\begin{equation*}
\begin{aligned}
\textbf{(TA)} &=\int_{0}^{\infty} z  \int_{0}^{z} \exp\left(\frac{1}{\rho^2}(c+1)\cdot|\sqrt{n}\beta_n^{\;(j)}||t| +  \frac{1}{\rho^2}2t^2\right) dt d\mathbb{P}_n(z) \\
&- \int_{-\infty}^{0} z \int_{z}^{0} \exp\left(\frac{1}{\rho^2}(c+1)\cdot|\sqrt{n}\beta_n^{\;(j)}||t| +  \frac{1}{\rho^2}2t^2\right) dt d\mathbb{P}_n(z)\\
&\leq  \int_{-\infty}^{\infty} z^2 \exp\left(\frac{1}{\rho^2(1+\rho^2)}z^2-\frac{1}{\rho^2(1+\rho^2)}\sqrt{n}\beta_n^{\;(j)}|z|\right) d\mathbb{P}_n(z)\\
&\leq \displaystyle\sup_{\mathbb{P}_n\in \mathcal{P}_{r,n}}\mathbb{E}_{\mathbb{P}_n}\Bigg[ (Z_{1,n}^{\;(j)})^2 \exp\left(\frac{2}{\rho^2}(Z_{1,n}^{\;(j)})^2+\frac{1}{\rho^2}(c+1)\cdot |\sqrt{n}\beta_n^{\;(j)}| |Z_{1,n}^{\;(j)}|\right)\Bigg]\leq \frac{D_8}{n}.
\end{aligned}
\end{equation*}
Proceeding in the same manner, we bound the second term $\textbf{(TB)}$ as
\begin{equation*}
\begin{aligned}
\sup_{\mathbb{P}_n\in \mathcal{P}_{r,n}}\mathbb{E}_{\mathbb{P}_n}\Bigg[ |Z_{1,n}^{\;(j)}|^3  \exp\left(\frac{2}{\rho^2}(Z_{1,n}^{\;(j)})^2+\frac{1}{\rho^2}(c+1)\cdot |\sqrt{n}\beta_n^{\;(j)}| |Z_{1,n}^{\;(j)}|\right)\Bigg]\leq \dfrac{D_{9}}{n^{3/2}}.
\end{aligned}
\end{equation*}
This completes Step 2 of our proof.
\end{proof}

\subsubsection{Proofs of main results}

We provide a proof for Theorem 1 and 2 below.
Supporting results for our proofs are collected in the preceding section.

\begin{proof} [Proof of Theorem 1]
Let us denote by $\mathcal{M}_b$ the collection of parameter sequences that satisfy Equation (4.1).
Note that 
\begin{equation}
\sup_{\mathcal{M}_b}\left(\mathbb{E}_{\mathcal{N}}\left[\bar{\Phi}\left(-\frac{1}{\rho}(\mathcal{Z}_n^{\;(j)} +\sqrt{n}\beta_n^{\;(j)})\right)\right]\right)^{-1} \leq E_1,
\label{bdd:den}
\end{equation}
where $E_1$ is a constant.
Now, observe that
\begin{equation*}
\begin{aligned}
\text{\normalfont SB}_{\mathbb{P}_n}(\rG_l)
%&\Big|\mathbb{E}_{\mathbb{P}_n}\left[ \rG_l(\mathcal{Z}_n^{\;(j)}; \sqrt{n}\beta^{\;(j)}_n) \right]- \mathbb{E}_{\mathcal{N}}\left[ \rG_l(\mathcal{Z}_n^{\;(j)}; \sqrt{n}\beta^{\;(j)}_n)\right]\Big|\\
%& \leq n\cdot \displaystyle\int_{-\infty}^{\infty}  \sup_{\alpha \in [0, 1]}  \mathbb{E}_{\mathbb{P}_n}\Bigg[\left(|t| +\frac{1}{\sqrt{n}} |\mathrm{e}_{1,n}^{\;(j)}|\right)\\
%& \;\;\;\;\;\;\;\;\;\;\;\;\;\;\;\;\;\;\;\;\;\;\;\;\;\;\;\;\;\;\;\;\;\;\;\;\;\;\;\;\;\;\;\;\;\;\;\;\;\;\;\;\;\;\times\Big|\cD^2\cS_{\rG_l}\left(\alpha t+(1-\alpha) \frac{1}{\sqrt{n}} \mathrm{e}_{1,n}^{\;(j)} + \mathcal{Z}_n^{\;(j)}[-1]\right)\Big|\Bigg] \mathrm{M}_1(t) dt\\
&\leq n\cdot \sup |\cD^2 \cS_{\rG_l}| \cdot \bigintsss_{-\infty}^{\infty} \mathbb{E}_{\mathbb{P}_n}\left[ |t| +\frac{1}{\sqrt{n}} |\mathrm{e}_{1,n}^{\;(j)}|\right]\mathrm{M}_1(t) dt.
\end{aligned}
\end{equation*}
Lemma \ref{fact:2} leads us to note that
$$\sup_{\mathcal{M}_b} \sup |\cD^2 \cS_{\rG_l}| <E_{l,1}.$$
%where $\mathcal{M}$ is the collection of all real-valued sequences $\{\beta_n^{\;(j)}: n\in \mathbb{N}\}$.
Thus, the Stein bound can further bounded as
\begin{equation*}
\begin{aligned}
%&\Big|\mathbb{E}_{\mathbb{P}_n}\left[ \rG_l(\mathcal{Z}_n^{\;(j)}; \sqrt{n}\beta^{\;(j)}_n) \right]- \mathbb{E}_{\mathcal{N}}\left[ \rG_l(\mathcal{Z}_n^{\;(j)}; \sqrt{n}\beta^{\;(j)}_n)\right]\Big|\\
\text{\normalfont SB}_{\mathbb{P}_n}(\rG_l)
&\leq \frac{E_{l,1}}{\sqrt{n}}\cdot \displaystyle\sup_n \displaystyle\sup_{\mathbb{P}_n\in \mathcal{P}_{b,n} } \left(\dfrac{1}{2}\mathbb{E}_{\mathbb{P}_n}\left[|\mathrm{e}_{1,n}|^3\right] + \mathbb{E}_{\mathbb{P}_n}\left[|\mathrm{e}_{1,n}|\right] \mathbb{E}_{\mathbb{P}_n}\left[\mathrm{e}_{1,n}^2\right]\right).
\end{aligned}
\end{equation*}
Then, the display in \eqref{bdd:den} together with Assumption 1 allows us to claim
$$\lim_n \sup_{\mathbb{P}_n\in \mathcal{P}_{b,n}} \rR^{(l)}_n =0 \  \text{ for } l\in [2].$$
As a result, we have that
$$\lim_n \sup_{\mathbb{P}_n\in \mathcal{P}_{b,n}}\Big|\widetilde{\mathbb{E}}_{\mathbb{P}_n}\left[\rH\circ\rP^{\;(j)}(\mathcal{Z}^{\;(j)}_n;  \sqrt{n}\beta^{\;(j)}_n)\right]- \widetilde{\mathbb{E}}_{\mathcal{N}}\left[\rH\circ\rP^{\;(j)}(\mathcal{Z}^{\;(j)}_n; \sqrt{n}\beta^{\;(j)}_n)\right]\Big| = 0.$$

\end{proof}

\begin{proof}[Proof of Theorem 2] 
Because the common denominator of our relative differences is bounded as 
\begin{equation*}
\begin{aligned}
\mathbb{E}_{\mathcal{N}}\left[\bar{\Phi}\left( -\frac{1}{\rho}{(\mathcal{Z}^{\;(j)}_n+\sqrt{n}\beta^{\;(j)}_n)}\right)\right] &= \bar{\Phi}\left(\frac{a_n\bar{\beta}}{\sqrt{1+\rho^2}}\right)\\
&\geq   \dfrac{E_2}{a_n  \bar{\beta}}\cdot\tE\left(a_n\bar{\beta},\frac{1}{1+\rho^2}\right),
\end{aligned}
\end{equation*} 
it suffices to show that
\begin{equation*}
\begin{aligned}
&\Scale[0.95]{\displaystyle\sup_{\mathbb{P}_n\in \mathcal{P}_{r,n}} a_n \bar{\beta}\cdot \left(\tE\left(a_n\bar{\beta},\frac{1}{1+\rho^2}\right)\right)^{-1}  \cdot \Bigg|\mathbb{E}_{\mathbb{P}_n}\Big[\rG_l(\mathcal{Z}_n^{\;(j)}; \sqrt{n}\beta_n^{\;(j)})\Big]-  \mathbb{E}_{\mathcal{N}}\Big[\rG_l(\mathcal{Z}_n^{\;(j)}; \sqrt{n}\beta_n^{\;(j)})\Big]\Bigg|}
\label{our goal}
\end{aligned}
\end{equation*}
converges to $0$ as $n$ tends to $\infty$. 
In the remaining proof, we provide a proof for the above statement.

Let $\bar{\rG}_l(z; \sqrt{n}\beta_n^{\;(j)})$, with the support set $\mathcal{S}_n$, be as defined in Lemma \ref{bar:G}.
Note that $\bar{\rG}_l(z; \sqrt{n}\beta_n^{\;(j)})$ agrees with $\rG_l(z; \sqrt{n}\beta_n^{\;(j)})$ whenever $z\in \mathcal{S}_n$. 
Thus, we observe that
$$
a_n \bar{\beta}\cdot \left(\tE\left(a_n\bar{\beta},\frac{1}{1+\rho^2}\right)\right)^{-1} \cdot \Bigg|\mathbb{E}_{\mathbb{P}_n}\Big[\rG_l(\mathcal{Z}_n^{\;(j)}; \sqrt{n}\beta_n^{\;(j)})\Big]-  \mathbb{E}_{\mathcal{N}}\Big[\rG_l(\mathcal{Z}_n^{\;(j)}; \sqrt{n}\beta_n^{\;(j)})\Big]\Bigg|
$$
is bounded by 
$$[\text{\textbf{T1}}_l] + [\text{\textbf{T2}}_l],$$
where
$$
[\text{\textbf{T1}}_l]:\ \  a_n \bar{\beta}\cdot \left(\tE\left(a_n\bar{\beta},\frac{1}{1+\rho^2}\right)\right)^{-1} \cdot  \Bigg|\mathbb{E}_{\mathbb{P}_n}\Big[\bar{\rG}_l(\mathcal{Z}_n^{\;(j)}; \sqrt{n}\beta_n^{\;(j)})\Big]-  \mathbb{E}_{\mathcal{N}}\Big[\bar{\rG}_l(\mathcal{Z}_n^{\;(j)}; \sqrt{n}\beta_n^{\;(j)})\Big]\Bigg|,
$$ 
 \begin{equation*}
\begin{aligned}
[\text{\textbf{T2}}_l]:\ \  &a_n \bar{\beta}\cdot \left(\tE\left(a_n\bar{\beta},\frac{1}{1+\rho^2}\right)\right)^{-1}\cdot \Bigg|\mathbb{E}_{\mathbb{P}_n}\Big[\rG_l(\mathcal{Z}_n^{\;(j)}; \sqrt{n}\beta_n^{\;(j)})\cdot \mathbf{1}_{\mathcal{S}^c_n}(\mathcal{Z}^{\;(j)}_n)\Big]\\
&\;\;\;\;\;\;\;\;\;\;\;\;\;\;\;\;\;\;\;\;\;\;\;\;\;\;\;\;\;\;\;\;\;\;\;\;\;\;\;\;\;\;\;\;\;\;\;\;\;\;\;\;\;\;\;\;\;\;\;\;\;\;\;\;\;\;\;\;\;\;\;\;\;\;\;\;\;\;\;\;\;\;\;\;\;\;-  \mathbb{E}_{\mathcal{N}}\Big[\rG_l(\mathcal{Z}_n^{\;(j)}; \sqrt{n}\beta_n^{\;(j)})\cdot \mathbf{1}_{\mathcal{S}^c_n}(\mathcal{Z}^{\;(j)}_n)\Big]\Bigg|.
\end{aligned}
\end{equation*}

Without loss of generality, we assume that $\bar{\rG}_l$ is non-negative valued and apply Theorem \ref{fact:4} to observe that
$$\sup_{\mathbb{P}_n\in \mathcal{P}_{r,n}} [\text{\textbf{T1}}_l] \leq \sup_{\mathbb{P}_n\in \mathcal{P}_{r,n}} a_n \bar{\beta}\cdot \left(\tE\left(a_n\bar{\beta},\frac{1}{1+\rho^2}\right)\right)^{-1} \text{\normalfont SB}_{\mathbb{P}_n}(\bar{\rG}_l) \leq \dfrac{a_n}{\sqrt{n}} \cdot E_{l,2}.$$
Otherwise, we can always write
 $$\bar{\rG}_l(z;\sqrt{n}\beta_n^{\;(j)})= \bar{\rG}_l^{+}(z;\sqrt{n}\beta_n^{\;(j)})- \bar{\rG}_l^{-}(z;\sqrt{n}\beta_n^{\;(j)}),$$
 where
$\Scale[0.95]{\bar{\rG}_l^{+}(z;\sqrt{n}\beta_n^{\;(j)}) = \max(\bar{\rG}_l(z;\sqrt{n}\beta_n^{\;(j)}), 0)},$ $\Scale[0.95]{\bar{\rG}_l^{-}(z;\sqrt{n}\beta_n^{\;(j)}) = -\min(\bar{\rG}_l(z;\sqrt{n}\beta_n^{\;(j)}), 0)}.$
We can then proceed similarly with the Stein bounds for the differences
\[\Big|\mathbb{E}_{\mathbb{P}_n}\left[\bar{\rG}_l^{+}(\mathcal{Z}_n^{\;(j)};\sqrt{n}\beta_n^{\;(j)})\right]- \mathbb{E}_{\mathcal{N}}\left[\bar{\rG}_l^{+}(\mathcal{Z}_n^{\;(j)};\sqrt{n}\beta_n^{\;(j)})\right]\Big|, \]
\[\Big|\mathbb{E}_{\mathbb{P}_n}\left[\bar{\rG}_l^{-}(\mathcal{Z}_n^{\;(j)};\sqrt{n}\beta_n^{\;(j)})\right]- \mathbb{E}_{\mathcal{N}}\left[\bar{\rG}_l^{-}(\mathcal{Z}_n^{\;(j)};\sqrt{n}\beta_n^{\;(j)})\right]\Big|.\]
As a result, we have
$$\lim_n\sup_{\mathbb{P}_n\in \mathcal{P}_{r,n}} [\text{\textbf{T1}}_l] =0.$$

To complete the proof, we bound the second term as follows:
\begin{equation*}
\begin{aligned}
\sup_{\mathbb{P}_n\in \mathcal{P}_{r,n}} [\text{\textbf{T2}}_l]  
&\leq  2E_{l,3}\cdot \displaystyle\sup_{\mathbb{P}_n\in \mathcal{P}_{r,n}} a_n \bar{\beta}\cdot \left(\tE\left(a_n\bar{\beta},\frac{1}{1+\rho^2}\right)\right)^{-1}\\
&\;\;\;\;\;\;\;\;\;\;\;\;\;\;\;\;\;\;\;\;\;\;\;\;\;\;\;\;\;\;\;\;\;\;\;\;\;\;\;\;\;\;\;\;\times \mathbb{E}_{\mathbb{P}_n}\left[\bar{\Phi}\left( -\frac{1}{\rho}{(\mathcal{Z}^{\;(j)}_n+\sqrt{n}\beta^{\;(j)}_n)}\right)\cdot \mathbf{1}_{\mathcal{S}^c_n}(\mathcal{Z}^{\;(j)}_n)\right]\\ 
&\leq 2E_{l,3}\cdot\displaystyle\sup_{\mathbb{P}_n\in \mathcal{P}_{r,n}} a_n \bar{\beta}\cdot \left(\tE\left(a_n\bar{\beta},\frac{1}{1+\rho^2}\right)\right)^{-1}\cdot \mathbb{P}_{\mathbb{P}_n}\left[ \mathcal{Z}^{\;(j)}_n\in \mathcal{S}^c_n \right].
%&\leq 2E_{l,4}\cdot\displaystyle\sup_{\mathbb{P}_n\in \mathcal{P}_{r,n}} a_n \bar{\beta}\cdot \left(\tE\left(a_n\bar{\beta},\frac{1}{1+\rho^2}\right)\right)^{-1}\cdot \tE(ca_n\bar{\beta},1).
\end{aligned}
\end{equation*}
Based on our assumptions, we have
$$
\displaystyle\sup_{n} \displaystyle\sup_{\mathbb{P}_n\in \mathcal{P}_{r,n}} \left(\tE(ca_n\bar{\beta},1)\right)^{-1}\cdot \mathbb{P}_{\mathbb{P}_n}\left[ \mathcal{Z}^{\;(j)}_n\in \mathcal{S}^c_n \right] 
< E_3.
$$
We fix $c>(1+\rho^2)^{-1}$.
Then, we can easily see that
$$\lim_n\sup_{\mathbb{P}_n\in \mathcal{P}_{r,n}} [\text{\textbf{T2}}_l] =0,$$
which proves the limit in our assertion.
As a consequence, we note that
$$\lim_n \sup_{\mathbb{P}_n\in \mathcal{P}_{r,n}}\Big|\widetilde{\mathbb{E}}_{\mathbb{P}_n}\left[\rH\circ\rP^{\;(j)}(\mathcal{Z}^{\;(j)}_n;  \sqrt{n}\beta_n^{\;(j)})\right]- \widetilde{\mathbb{E}}_{\mathcal{N}}\left[\rH\circ\rP^{\;(j)}(\mathcal{Z}_n^{\;(j)}; \sqrt{n}\beta_n^{\;(j)})\right]\Big| = 0.$$
\end{proof}

%%%%%%%%%%%%%%%%%%%%%%%%%%%%%%%%%%%%%%%%%%%%%%%%%%%%%%%%%%%%%%%%%%%%%%%
%%%%%%%%%%%%%%%%%%%%%%%%%%%%%%%%%%%%%%%%%%%%%%%%%%%%%%%%%%%%%%%%%%%%%%%

\subsection{Proofs for Section 5}

As before, uppercase letters are used to denote constants that do not depend on $n$. 
We begin by deriving the multivariate Stein bound in Lemma 2 and obtain smoothness properties for our pivot as stated in Proposition 7.
%Then, we detail out the theory of weak convergence in Section 5 of the paper.

\subsubsection{Proofs of Lemma 2 and Proposition 7}
\begin{proof}[Proof of Lemma 2]
Let $\mathrm{e}^{\star}_{i,n}$ be an independent copy of $\mathrm{e}_{i,n}$, and let $I$ be a uniform variable over $[n]$ that is independent of $\mathcal{Z}_n$.
Consider an exchangeable pair of variables $$\begin{pmatrix}{\mathcal{Z}^{\star}_n}' & {\mathcal{Z}_n}'\end{pmatrix}'$$ 
where we define:
$$\mathcal{Z}^{\star}_n= \mathcal{Z}_n - Z_{I,n} + Z^{\star}_{I,n}=\mathcal{Z}_n - \frac{1}{\sqrt{n} }\mathrm{e}_{I,n}+ \frac{1}{\sqrt{n}}\mathrm{e}^{\star}_{I,n}.$$
Let $\Theta$ be a random matrix whose $(j,k)^{\text{th}}$ entry assumes the value
$$\frac{1}{2}\mathbb{E}_{\mathbb{P}_n}\left[(\mathrm{e}_{I,n}^{\;(j)} - \mathrm{e}_{I,n}^{\star (j)})(\mathrm{e}_{I,n}^{(k)} - \mathrm{e}_{I,n}^{\star (k)}) \lvert  \mathcal{Z}_n\right] - I_{d,d}[j,k].$$
Based on Equation 11 in \cite{chatterjee2007multivariate}, we have
$$ \mathbb{E}_{\mathbb{P}_n}[\cS_\rg(\mathcal{Z}^{\star}_n)-\cS_\rg(\mathcal{Z}_n)]=0.$$
Using a second order Taylor series expansion, we now note that
\begin{equation*}
\begin{aligned}
n\cdot \mathbb{E}_{\mathbb{P}_n}[\cS_\rg(\mathcal{Z}^{\star}_n)-\cS_\rg(\mathcal{Z}_n)]
&= \mathbb{E}_{\mathbb{P}_n}[\rg(\mathcal{Z}_n)] - \mathbb{E}_{\mathcal{N}}[\rg(\mathcal{Z}_n)] \\
&+ \mathbb{E}_{\mathbb{P}_n}\left[\text{Tr}(\Theta  \cS''_{\rg}(\mathcal{Z}_n)) + n\cdot \mathcal{R}_n(\mathcal{Z}_n,  \mathcal{Z}_n^{\star})\right],
\end{aligned}
\end{equation*}
where $\mathcal{R}_n(\mathcal{Z}_n,  \mathcal{Z}_n^{\star})$ denotes the error term in the Taylor expansion.
Therefore, the difference in expectations
$$|\mathbb{E}_{\mathbb{P}_n}[\rg(\mathcal{Z}_n)] - \mathbb{E}_{\mathcal{N}}[\rg(\mathcal{Z}_n)]|$$ 
is bounded from above by
\begin{equation*}
\text{(UB)}: \ \ |\mathbb{E}_{\mathbb{P}_n}\left[\text{Tr}(\Theta  \cS''_{\rg}(\mathcal{Z}_n))\right]| + n\cdot | \mathbb{E}_{\mathbb{P}_n}\left[ \mathcal{R}_n(\mathcal{Z}_n,  \mathcal{Z}_n^{\star})\right]|.
\end{equation*}
Simplifying at first the error term from the Taylor expansion, observe that
$$n\cdot \mathbb{E}_{\mathbb{P}_n}[ \mathcal{R}_n(\mathcal{Z}_n,  \mathcal{Z}_n^{\star})]$$ is equal to
\begin{align*}
&\dfrac{D_0}{\sqrt{n}} \cdot \sum_{j, k, l} \mathbb{E}_{\mathbb{P}_n}\Bigg [ (\mathrm{e}_{1,n}^{\;(j)} - \mathrm{e}_{1,n}^{\star(j)})(\mathrm{e}_{1,n}^{(k)} - \mathrm{e}_{1,n}^{\star(k)}) (\mathrm{e}_{1,n}^{(l)} - \mathrm{e}_{1,n}^{\star(l)})\\ \nonumber
&\;\;\;\;\;\;\;\;\;\;\;\;\;\;\;\;\;\;\;\;\;\;\;\;\;\;\;\;\;\;\;\times\sup_{\alpha\in[0,1]} \cD^3\cS_{\rg}\left( \mathcal{Z}_n -(1-\alpha) \frac{1}{\sqrt{n}}\cdot (\mathrm{e}_{1,n} -  \mathrm{e}^{\star}_{1,n})\right)[j,k, l] \Bigg]. \nonumber
\end{align*}
By using the fact that: $|\mathrm{e}_{1,n}^{\;(j)} - \mathrm{e}_{1,n}^{\star(j)}| \leq (\|\mathrm{e}_{1,n}\| + \|\mathrm{e}^{\star}_{1,n}\|)$, we note that this term is further bounded by
\begin{equation}
\begin{aligned}
& \dfrac{D_1}{\sqrt{n}} \sum_{\lambda,\gamma\in \{0,1,2,3\}:\lambda+\gamma=3} \sum_{j, k, l} \mathbb{E}_{\mathbb{P}_n}\Bigg[\|\mathrm{e}_{1,n}\|^{\lambda}\|\mathrm{e}^{\star}_{1,n}\|^{\gamma} \sup_{\alpha\in [0,1]}\Big|\cD^3\cS_{\rg}\Big(\mathcal{Z}_n[-1] +\\
&\;\;\;\;\;\;\;\;\;\;\;\;\;\;\;\;\;\;\;\;\;\;\;\;\;\;\;\;\;\;\;\;\;\;\;\;\;\;\;\;\;\;\;\;\;\;\;\;\;\;\;\;\;\;\;\;\;\;\;\;\;\;\;\;\;\;\;\;\;\;\;\;\; \alpha\frac{1}{\sqrt{n}}\mathrm{e}_{1,n} + (1-\alpha)\frac{1}{\sqrt{n}}\mathrm{e}^{\star}_{1,n}\Big)[j, k, l]\Big|\Bigg].
\end{aligned}
\label{B1}
\end{equation}
Letting $\delta[ j,k]=0$ when $j\neq k$ and $1$ when $j=k$, we observe that
%$$\mathbb{E}_{\mathbb{P}_n}\left[\text{Tr}(\Theta \cS''_{\rg}(\mathcal{Z}_n))\right]$$
%is equal to
\begin{align*}
&\mathbb{E}_{\mathbb{P}_n}\left[\text{Tr}(\Theta \cS''_{\rg}(\mathcal{Z}_n))\right]\\
&=\dfrac{1}{2n}\cdot \sum_{j, k} \sum_{i=1}^n \mathbb{E}_{\mathbb{P}_n}\left[ \left(\mathbb{E}_{\mathbb{P}_n}\left[ (\mathrm{e}_{i,n}^{\;(j)} - \mathrm{e}_{i,n}^{\star(j)})(\mathrm{e}_{i,n}^{(k)} - \mathrm{e}_{i,n}^{\star(k)}) | \mathcal{Z}_n\right]  -2\cdot \delta[ j,k] \right) \cdot \cD^2\cS_{\rg}(\mathcal{Z}_n)[j,k]\right]\\
&= \dfrac{1}{2} \cdot \sum_{j,k} \mathbb{E}_{\mathbb{P}_n}\left[ (\mathrm{e}_{1,n}^{\;(j)}  \mathrm{e}_{1,n}^{(k)}- \delta[j,k] )\cdot \cD^2\cS_{\rg}(\mathcal{Z}_n)[j,k]\right]\\ \nonumber
&=  \dfrac{1}{2} \cdot \sum_{j,k} \mathbb{E}_{\mathbb{P}_n}\left[ (\mathrm{e}_{1,n}^{\;(j)}  \mathrm{e}_{1,n}^{(k)} - \delta[j,k] )\cdot \left(\cD^2\cS_{\rg}(\mathcal{Z}_n)[j,k] - \cD^2\cS_{\rg}(\mathcal{Z}_n[-1])[j,k]  \right)\right]\\
&\leq \dfrac{D_2}{2\sqrt{n}} \cdot \displaystyle\sum_{j,k,l} \mathbb{E}_{\mathbb{P}_n}\Bigg[  (\mathrm{e}_{1,n}^{\;(j)}  \mathrm{e}_{1,n}^{(k)} - \delta[j,k] ) \cdot \mathrm{e}_{1,n}^{(l)}\\
&\;\;\;\;\;\;\;\;\;\;\;\;\;\;\;\;\;\;\;\;\;\;\;\;\;\;\;\;\;\;\;\;\;\;\;\;\;\;\;\;\;\;\;\;\;\;\;\;\;\;\times \sup_{\alpha\in[0,1]}\Big |\cD^3\cS_{\rg}\left(\mathcal{Z}_n[-1]+\alpha \frac{1}{\sqrt{n}}\mathrm{e}_{1,n} \right)[j, k, l]\Big|\Bigg],
\end{align*}
and that the expression in the last display is further bounded from above by:
\begin{equation}
 \dfrac{D_2}{2\sqrt{n}} \cdot \displaystyle\sum_{j,k,l} \mathbb{E}_{\mathbb{P}_n}\Bigg[ (\|\mathrm{e}_{1,n}\|^3 + \|\mathrm{e}_{1,n}\|) \sup_{\alpha\in[0,1]}\Big |\cD^3\cS_{\rg}\left(\mathcal{Z}_n[-1]+\alpha \frac{1}{\sqrt{n}}\mathrm{e}_{1,n} \right)[j, k, l]\Big|\Bigg].
\label{B2}
\end{equation}
The Stein bound follows by combining the two bounds in \eqref{B1} and \eqref{B2}.
\end{proof}

Before stating a proof for Proposition 7, we introduce some notations.
For $j\in \mathbb{N}$, we will use the symbol $l_{j}(v)$ to denote a linear mapping in $v$. 
We use this symbol when we need to use the linearity of the mapping without specifying its actual form.

Our proof relies on an alternate represenation for our pivot.
For this purpose, let $T_n^{\;(j)}$ be a random variable such that 
$$\begin{pmatrix} (V_n^{\;(j)}) & (T_n^{\;(j)})' \end{pmatrix}'$$
given $U_n^{\;(j)}= \oU$ have the following density
\begin{align*}
p(v,t) &\propto \tE\left(Qt -R^{\;(j)}\begin{pmatrix} v & (\oU)' \end{pmatrix}' + r,  \frac{1}{\rho^2}\Sigma^{-1}\right)\cdot \tE\left(v-\sqrt{n}\beta_n^{\;(j)}, \frac{1}{\sigma_j^2}\right). \numberthis
\label{joint:lik}
\end{align*}
Also, define $\mathcal{T}_n^{\;(j)}$ to be a variable that is distributed as $T_n^{\;(j)} \lvert U_n^{\;(j)}=\oU$. 

\begin{proof}[Proof of Proposition 7] 
We start from re-writing our pivot as 
\begin{equation*}
\begin{aligned}
&\dfrac{\mathbb{E}\left[\mathbb{P}\left[V_n^{\;(j)}\geq \oV\; \lvert \; U_n^{\;(j)}= \oU, T_n^{\;(j)}= \oT \right]\cdot \mathbf{1}_{\oT>0_{p}} \; \lvert \; U_n^{\;(j)}=\oU\right]}{\mathbb{P}\left[T_n^{\;(j)} >0_p \; \lvert \; U_n^{\;(j)}=\oU\right]}.
\end{aligned}
\end{equation*}
Now, we use two facts.
First, 
$$V_n^{\;(j)} \lvert \; T_n^{\;(j)}= \oT, \; U_n^{\;(j)}= \oU \; \; \text{ and }\  \ T_n^{\;(j)} \lvert U_n^{\;(j)}=\oU$$
have Gaussian distributions with means 
$$l_{1}\left(\begin{pmatrix} {\oT}' & {\oU}' & {\sqrt{n}\beta_n^{\;(j)}}\end{pmatrix}'\right) \text{ and } l_{2}\left(\begin{pmatrix} {\oU}' & {\sqrt{n}\beta_n^{\;(j)}}\end{pmatrix}'\right),$$
respectively. 
Both distributions have a constant covariance matrix. 
Second, observe that
$$\begin{pmatrix} (\oV)' & (\oU)' \end{pmatrix}' = ({R}^{\;(j)})^{-1}\left(\Sigma^{1/2} \oZ+ \sqrt{n}\beta_n\right).$$
Using the above-stated facts, we note that
$\mathbb{P}\left[V_n^{\;(j)}\geq \oV\; \lvert \; U_n^{\;(j)}= \oU, T_n^{\;(j)}= \oT \right]$ takes the form
$$\bar{\Phi}\left(l_3\left(\begin{pmatrix} {\oT}' & \mathcal{Z}'_{n;\text{obs}}& \sqrt{n}\beta'_n \end{pmatrix}'\right)\right),$$
where $l_3$ is a linear mapping in its arguments.
Thus, our pivot simplifies as
$$
\dfrac{\mathbb{E}\left[\bar{\Phi}\left(l_3\left(\begin{pmatrix} {\oT}' & \mathcal{Z}'_{n;\text{obs}} & \sqrt{n}\beta'_n \end{pmatrix}'\right)\right)\cdot \mathbf{1}_{\oT>0_{p}} \; \Big\lvert \; U_n^{\;(j)}=\oU\right]}{\mathbb{P}\left[T_n^{\;(j)} >0_p \; \lvert \; U_n^{\;(j)}=\oU\right]}.
$$
Taking partial derivatives of our pivot, it follows that
 $$\|\cD^{p_0} \rP^{\;(j)} \left(\mathcal{Z}_n; \sqrt{n}\beta_n\right)\|\leq D_3\cdot \left( \mathbb{E}\left[\|\mathcal{T}_n^{\;(j)}- \mathbb{E}[\mathcal{T}_n^{\;(j)}]\|^{p_0}\; \Big\lvert  \; \mathcal{T}_n^{\;(j)}>0_{p} \right] + \|\mathbb{E}[\mathcal{T}_n^{\;(j)}]\|^{p_0}\right).$$
At last, we note that
$$\mathbb{E}[\mathcal{T}_n^{\;(j)}]= l_4\left(\begin{pmatrix}  \mathcal{Z}'_{n;\text{obs}} & \sqrt{n}\beta'_n \end{pmatrix}'\right)$$ 
for a linear mapping $l_4$ and that
\begin{equation*}
\begin{aligned}
\mathbb{E}\left[\|\mathcal{T}_n^{\;(j)}- \mathbb{E}[\mathcal{T}_n^{\;(j)}]\|^{p_0}\; \Big\lvert  \; \mathcal{T}_n^{\;(j)}>0_{p} \right] &\leq D_4\cdot \Big\| l_4\left(\begin{pmatrix}  \mathcal{Z}'_{n;\text{obs}} & \sqrt{n}\beta'_n \end{pmatrix}'\right)\Big\|^{p_0}.
%&\leq \sum_{\lambda, \gamma \in \{0\} \cup[p_0]: \lambda+\gamma\leq p_0} C_2^{\lambda, \gamma}\|\oZ\|^{\lambda} \|\sqrt{n}\beta_n\|^{\gamma}.
\end{aligned}
\end{equation*}
Our bound on the partial derivatives follows immediately.
\end{proof}

\subsubsection{Supporting results}

The supporting results in this section to prove Theorem 4.
\begin{lemma}
For $t, \alpha, \kappa \in [0,1]$, define
$$W_{t,\alpha, \kappa} = \sqrt{t}\left(\mathcal{Z}_n[-1] + \frac{\alpha}{\sqrt{n}}\mathrm{e}_{1,n}+ \frac{\kappa}{\sqrt{n}}\mathrm{e}^{\star}_{1,n}\right).$$
Let $\lambda, \gamma \in {0}\cup[3]$ such that $\lambda+\gamma\leq 3$.
Under the conditions stated in Assumptions 5 and 6, it holds that
\begin{equation*}
\begin{aligned}
& \left({\normalfont \tE}\left(a_n \bar{\beta}, \frac{1}{(1+\rho^2)}\cdot( \Lambda +\Delta)\right)\right)^{-1} \cdot \mathbb{E}_{\mathbb{P}_n}\Bigg[ \|\mathrm{e}_{1,n}\|^{\lambda}\|\mathrm{e}^{\star}_{1,n}\|^{\gamma} \cdot \sup_{\alpha, \kappa\in [0,1]}\displaystyle\int_{0}^{1}\sqrt{t} \\
&\;\;\;\;\;\;\;\;\;\;\;\;\;\;\;\;\;\;\;\;\;\;\;\;\;\;\;\;\;\;\;\;\;\;\;\;\;\;\;\;\;\;\;\;\;\;\;\;\;\;\times{\normalfont\tE }\left(\Sigma^{1/2}W_{t,\alpha, \kappa} + \sqrt{n}\beta_n-r, \frac{\Lambda + \Delta}{(1-t + \rho^2)} \right)  dt\Bigg]\leq D_5.
\end{aligned}
\end{equation*}
\label{bound:main:proof}
\end{lemma}

\begin{proof}
Fixing some notations for our proof, denote by $\text{Ev}_{\text{max}}$ the largest eigen value of $(1+\rho^2)^{-1}\cdot (\Lambda + \Delta)$.
Let
$$ 
\Theta(t) = \rho^{-2}\left(t(1-t+\rho^2)^{-1}\cdot \Sigma^{1/2}(\Lambda + \Delta) \Sigma^{1/2} + I\right)^{-1} \Sigma^{1/2} (\Lambda + \Delta).
$$
Let $\Theta_k(t)$ be the $k^{\text{th}}$ row of $\Theta_k$ and let $\|\Theta(t)\|_{\text{max}}= \displaystyle\max_{k\in [d]} \|\Theta_k(t)\|$.
Fix  
$$c_0 > \max\left(\text{Ev}_{\text{max}}^{1/2}(\|\bar{\beta}\|+1), \sup_{t\in [0,1]} \|\Theta(t)\|_{\text{max}} \cdot(\|\bar{\beta}\|+1)\right),$$
and let $\mathcal{R}= [-c_0 \cdot 1_d, c_0 \cdot 1_d] \subseteq \real^d$.

%Without any loss of generality, we focus on the case when $\bar{\lambda}=\bar{\kappa}=\breve{\kappa}=0$, $\breve{\lambda}=3$; the same line of reasoning will apply to other values of $\bar{\lambda}$,$\bar{\kappa}$,$\breve{\lambda}$, $\breve{\kappa}$.
%In this case, we 
%$$a_n^3 \mathbb{E}_{\mathbb{P}_n}\Big[ \|\mathrm{e}_{1,n}\|^{\lambda}\|\mathrm{e}^{\star}_{1,n}\|^{\gamma} \cdot \sup_{\alpha, \kappa\in [0,1]}\displaystyle\int_{0}^{1}\sqrt{t} {\normalfont\tE }\left(\Sigma^{1/2}W_{t,\alpha, \kappa} + \sqrt{n}\beta_n-r, \frac{1}{(1-t + \rho^2)} \cdot (\Lambda + \Delta)\right)  dt\Big].$$
Using the independence between $\mathcal{Z}_n[-1]$, $\mathrm{e}_{1,n}$, and $\mathrm{e}^{\star}_{1,n}$, we have
\begin{equation*}
\begin{aligned}
&\Scale[0.95]{\mathbb{E}_{\mathbb{P}_n}\Big[ \|\mathrm{e}_{1,n}\|^{\lambda}\|\mathrm{e}^{\star}_{1,n}\|^{\gamma} \cdot \displaystyle\sup_{\alpha, \kappa\in [0,1]}\displaystyle\int_{0}^{1}\sqrt{t} {\normalfont\tE }\left(\Sigma^{1/2}W_{t,\alpha, \kappa} + \sqrt{n}\beta_n-r, \frac{1}{(1-t + \rho^2)} \cdot (\Lambda + \Delta)\right)  dt\Big]}\\
&\Scale[0.95]{\leq\mathbb{E}_{\mathbb{P}_n}\Big[ \|\mathrm{e}_{1,n}\|^{\lambda}\|\mathrm{e}^{\star}_{1,n}\|^{\gamma} \cdot \displaystyle\sup_{\alpha, \kappa\in [0,1]}\displaystyle\int_{0}^{1}\sqrt{t} \tE\left(\Sigma^{1/2}W_{t,\alpha, \kappa} + \sqrt{n}\beta_n-r, \frac{1}{(1-t + \rho^2)} \cdot (\Lambda + \Delta)\right)}\\
&\Scale[0.95]{\;\;\;\;\;\;\;\;\;\;\;\;\;\;\;\;\;\;\times \mathbf{1}_{\mathcal{R}}((a_n)^{-1}\mathcal{Z}_n[-1]) dt\Big]+ \mathbb{E}_{\mathbb{P}_n}\Big[ \|\mathrm{e}_{1,n}\|^{\lambda}\|\mathrm{e}^{\star}_{1,n}\|^{\gamma} \Big] \mathbb{E}_{\mathbb{P}_n}\Big[\mathbf{1}_{\mathcal{R}^c}((a_n)^{-1}\mathcal{Z}_n[-1])\Big]}
\end{aligned}
\end{equation*}
\begin{equation*}
\begin{aligned}
&\Scale[0.95]{\leq \mathbb{E}_{\mathbb{P}_n}\left[\|\mathrm{e}_{1,n}\|^{\lambda}\|\mathrm{e}^{\star}_{1,n}\|^{\gamma} \exp(\chi\|\mathrm{e}_{1,n}\|)\right]} \\
&\Scale[0.95]{\;\;\;\;\;\times\displaystyle\bigintsss_{0}^{1}\sqrt{t}\cdot\mathbb{E}_{\mathbb{P}_n}\Bigg[  \tE\left(\sqrt{t}\Sigma^{1/2}\mathcal{Z}_n[-1] -a_n\bar{\beta} , \frac{1}{(1-t + \rho^2)}\cdot (\Lambda + \Delta)\right) \cdot\mathbf{1}_{\mathcal{R}}((a_n)^{-1}\mathcal{Z}_n[-1])  \Bigg]dt }\\
&\Scale[0.95]{+ \mathbb{E}_{\mathbb{P}_n}\Big[ \|\mathrm{e}_{1,n}\|^{\lambda}\|\mathrm{e}^{\star}_{1,n}\|^{\gamma} \Big] \mathbb{E}_{\mathbb{P}_n}\Big[\mathbf{1}_{\mathcal{R}^c}((a_n)^{-1}\mathcal{Z}_n[-1])\Big].}
\end{aligned}
\end{equation*} 

Let
%$$\Psi_t(z) = \tE\left( \sqrt{t}\Sigma^{1/2}z - \bar{\beta}, \frac{1}{(1-t + \rho^2)}\cdot (\Lambda + \Delta)\right).$$
$$\Psi_t(z)=  \frac{1}{1-t +\rho^2}(\sqrt{t}\Sigma^{1/2}z-\bar\beta)'(\Lambda +\Delta)(\sqrt{t}\Sigma^{1/2}z-\bar\beta),$$
and let $\mathcal{K}_t = \mathcal{R}$ for $t\in (0,1]$. 
An application of the Varadhan's principle under Assumptions 5 and 6 results in:
\begin{equation*}
\begin{aligned}
& \sup_n\sup_{\mathbb{P}_n\in \mathcal{P}_{r,n}}\sup_{t\in(0,1]}\left(\tE\left(a_n \bar{\beta}, \frac{1}{(1+\rho^2)}\cdot( \Lambda +\Delta)\right)\right)^{-1} \\
&\;\;\;\;\;\;\;\;\;\;\;\;\;\;\;\;\;\;\;\;\;\times \mathbb{E}_{\mathbb{P}_n}\left[\exp\Big(-a_n^2\Psi_t\Big(\frac{1}{a_n}\mathcal{Z}_n[-1]\Big)\Big)\cdot  \mathbf{1}_{\mathcal{K}_t}((a_n)^{-1}\mathcal{Z}_n[-1]) \right]<\infty,
\end{aligned}
\end{equation*}  
and leads us to observe that 
\begin{equation*}
\begin{aligned}
&\left(\tE\left(a_n \bar{\beta}, \frac{1}{(1+\rho^2)}\cdot( \Lambda +\Delta)\right)\right)^{-1}  \\
&\;\times \displaystyle\bigintsss_{0}^{1}\sqrt{t}\cdot\mathbb{E}_{\mathbb{P}_n}\Bigg[  \tE\left(\sqrt{t}\Sigma^{1/2}\mathcal{Z}_n[-1] -a_n\bar{\beta} , \frac{1}{(1-t + \rho^2)}\cdot (\Lambda + \Delta)\right) \cdot\mathbf{1}_{\mathcal{R}}((a_n)^{-1}\mathcal{Z}_n[-1])  \Bigg]dt\\
&\leq  \sup_n\sup_{\mathbb{P}_n\in \mathcal{P}_{r,n}}\sup_{t\in(0,1]}\left(\tE\left(a_n \bar{\beta}, \frac{1}{(1+\rho^2)}\cdot( \Lambda +\Delta)\right)\right)^{-1}\\
&\;\;\;\;\;\;\;\;\;\;\;\;\times \mathbb{E}_{\mathbb{P}_n}\left[\exp\Big(-a_n^2\Psi_t\Big(\frac{1}{a_n}\mathcal{Z}_n[-1]\Big)\Big)\cdot  \mathbf{1}_{\mathcal{K}_t}((a_n)^{-1}\mathcal{Z}_n[-1]) \right] \cdot  \displaystyle\bigintsss_{0}^{1} \sqrt{t} dt<\infty.
\end{aligned}
\end{equation*}  

Let $\mathcal{K}_0= \mathcal{R}^c$. 
Under the same assumptions, we also have:
\begin{equation*}
\begin{aligned}
&\left(\tE\left(a_n\bar{\beta}, \frac{1}{(1+\rho^2)}\cdot( \Lambda +\Delta)\right)\right)^{-1} \cdot  \mathbb{E}_{\mathbb{P}_n}\Big[\mathbf{1}_{\mathcal{R}^c}((a_n)^{-1}\mathcal{Z}_n[-1])\Big]\\
&\leq \sup_n\sup_{\mathbb{P}_n\in \mathcal{P}_{r,n}} \left(\tE\left(a_n\bar{\beta}, \frac{1}{(1+\rho^2)}\cdot( \Lambda +\Delta)\right)\right)^{-1} \cdot  \mathbb{E}_{\mathbb{P}_n}\Big[\mathbf{1}_{\mathcal{K}_0}((a_n)^{-1}\mathcal{Z}_n[-1])\Big] <\infty,
\end{aligned}
\end{equation*}  
which follows by noting that $c_0^2>\text{Ev}_{\text{max}}(\|\bar{\beta}\|+1)^2$. 
Our moment assumptions imply that
$$\sup_n\sup_{\mathbb{P}_n\in \mathcal{P}_{r,n}} \mathbb{E}_{\mathbb{P}_n}\left[\|\mathrm{e}_{1,n}\|^{\lambda}\exp(\chi\|\mathrm{e}_{1,n}\|)\right]<\infty, \; \sup_n\sup_{\mathbb{P}_n\in \mathcal{P}_{r,n}} \mathbb{E}_{\mathbb{P}_n}\left[\|\mathrm{e}^{\star}_{1,n}\|^{\gamma}\right]<\infty.$$
This completes our proof.
\end{proof}

\begin{proposition}
\label{rep}
\rm{
Define the functions
$\rg_1(z; \sqrt{n}\beta_n) = 1$ and $\rg_2(z; \sqrt{n}\beta_n) = \rH\circ \rP^{\;(j)}(z\; ;\sqrt{n}\beta_n)$.
Let $\mathcal{M}_r=\{\beta_n: \ n \in \mathbb{N}\}$ such that $\beta_n$ is parameterized as per Equation (5.2) for each $n$.
There exists a Lebesgue-almost everywhere differentiable function $\mathcal{A}_n(\cdot; \sqrt{n}\beta_n):\real^d \to \real$ and a subset $\mathcal{S}_n\subseteq \real^d$ such that
\begin{equation*}
\begin{aligned}
\bar{\rG}_l(z; \sqrt{n}\beta_n) &= \rg_l(z; \sqrt{n}\beta_n)\cdot  \tE\left(\Sigma^{1/2}z + \sqrt{n}\beta_n -r, \frac{1}{\rho^2}(\Lambda + \Delta)\right)\\
&\;\;\;\;\;\;\;\;\;\;\;\;\;\;\;\;\;\;\;\;\;\;\;\;\;\;\;\;\;\;\;\;\;\;\;\;\;\;\;\;\;\;\;\;\;\;\;\;\;\;\;\;\;\;\;\;\;\;\;\;\;\;\;\;\;\;\;\;\;\;\;\;\;\;\;\;\;\;\;\;\;\times\mathcal{A}_n(z; \sqrt{n}\beta_n^{\;(j)})\cdot \mathbf{1}_{z\in \mathcal{S}_n} 
\end{aligned}
\end{equation*}
agrees with $\rG_l(z; \sqrt{n}\beta_n)$ on its support set.
Additionally, we have
$$\lim_n  \; \dfrac{\mathbb{E}_{\mathcal{N}}\Big[\rF(\Sigma^{1/2}\mathcal{Z}_n + \sqrt{n}\beta_n) \cdot \mathbf{1}_{\mathcal{Z}_n \in \mathcal{S}_n^c}\Big]}{\mathbb{E}_{\mathcal{N}}\Big[\rF(\Sigma^{1/2}\mathcal{Z}_n + \sqrt{n}\beta_n)\Big]}, \text{ and } $$
%$$\lim_n\sup_{\beta_n \in \mathcal{M}_{r,n}}\;  \dfrac{\mathbb{E}_{\mathcal{N}}\Big[\rF(\Sigma^{1/2}\mathcal{Z}_n + \sqrt{n}\beta_n) \cdot \mathbf{1}_{\mathcal{Z}_n \in \mathcal{S}_n^c}\Big]}{\mathbb{E}_{\mathcal{N}}\Big[\rF(\Sigma^{1/2}\mathcal{Z}_n + \sqrt{n}\beta_n)\Big]}=0,$$
$$\displaystyle\sup_{ \mathcal{M}_r}\sup_{z\in \mathcal{S}_n} (a_n)^{|\mathcal{I}|} \|\cD^m\mathcal{A}_n(z; \sqrt{n}\beta_n)\| <\infty, \text{ for } \ m=0,1,2,3.$$
} 
\end{proposition}

To see a proof for Proposition \ref{rep}, we make two additional observations in Lemma \ref{rel:1} and Lemma \ref{rel:2}.

\begin{lemma}
\label{rel:1}
Suppose that $T_n^{\star} \sim \mathcal{N}(\vartheta_z, \rho^2 \bar{\Sigma})$, where $R = -\bar{\Sigma}Q'\Sigma^{-1/2}$, $\vartheta_z = -(R z + a_n \bar{\mu})$, and $\bar{\Sigma} = (Q' \Sigma^{-1} Q)^{-1}$.
It holds that
$$\rF(\Sigma^{1/2} z + \sqrt{n}\beta_n) \propto \text{\normalfont Exp}\left(\Sigma^{1/2}z + \sqrt{n}\beta_n -r, \frac{1}{\rho^2}\Lambda\right) \cdot \mathbb{P}\left[T^{\star}_n > 0\right].$$
\end{lemma}

The proof for Lemma \ref{rel:1} follows directly from the definition of the function $\rF$.

We need some more notations to state the next result.
Let $\mathcal{I}$ and its complement $\mathcal{J}= \mathcal{I}^c$ be as defined in Proposition 6.
Whenever $\mathcal{J}\neq \emptyset$,  we consider the partition  
$$\mathcal{J} =\mathcal{J}_1 \cup \mathcal{J}_2$$ such that
$$\bar{\Sigma}_{\mathcal{J}_1, \mathcal{I}}\bar{\Sigma}_{\mathcal{I}, \mathcal{I}}^{-1} \bar{\mu}^{(\mathcal{I})} - \bar{\mu}^{(\mathcal{J}_1)}>0, \text{ and } \bar{\Sigma}_{\mathcal{J}_2, \mathcal{I}}\bar{\Sigma}_{\mathcal{I}, \mathcal{I}}^{-1} \bar{\mu}^{(\mathcal{I})} - \bar{\mu}^{(\mathcal{J}_2)}=0.$$
For any $\mathcal{L} \subseteq \mathcal{J}_2$, define 
$$\bar{I}_{\mathcal{L}}= \mathcal{I} \cup \mathcal{L}$$ 
and correspondingly, let $\bar{J}_{\mathcal{L}}= \bar{I}_{\mathcal{L}}^c$.
Let
\begin{equation*}
\begin{aligned}
\bar{\Theta} &= \left(\bar{\Sigma}_{\bar{I}_{\mathcal{L}}, \bar{I}_{\mathcal{L}}}\right)^{-1}= \begin{bmatrix} \bar{\Theta}_{\mathcal{I} , \mathcal{I} } & \bar{\Theta}_{\mathcal{I} , \mathcal{L} } \\ \bar{\Theta}_{\mathcal{L} , \mathcal{I} }& \bar{\Theta}_{\mathcal{L} , \mathcal{L} }  \end{bmatrix},
\end{aligned}
\end{equation*}
where the submatrices are partitioned as per the indices in $\mathcal{I}$ and $\mathcal{L}$.
Also, define
\begin{equation*}
\begin{aligned}
\rl_1(z)&=  \bar{\Sigma}_{\mathcal{L}, \mathcal{I}} \bar{\Sigma}^{-1}_{\mathcal{I}, \mathcal{I}} (Rz)^{(\mathcal{I})} - (Rz)^{(\mathcal{L})},\\
\rl_2(z)&=  \bar{\Sigma}^{-1}_{\mathcal{I}, \mathcal{I}} (Rz)^{(\mathcal{I})} -  \bar{\Theta}_{\mathcal{I}, \mathcal{L}}\rl_1(z),\\
\rl_3(z) &=  \bar{\Sigma}_{\bar{J}_{\mathcal{L}}, \mathcal{I}} \bar{\Sigma}^{-1}_{\mathcal{I}, \mathcal{I}} (Rz)^{(\mathcal{I})} - (Rz)^{(\bar{J}_{\mathcal{L}})} - (\bar{\Sigma}_{\bar{J}_{\mathcal{L}}, \mathcal{I}} \bar{\Theta}_{\mathcal{I}, \mathcal{L}}+\bar{\Sigma}_{\bar{J}_{\mathcal{L}}, \mathcal{L}} \bar{\Theta}_{\mathcal{L}, \mathcal{L}})\rl_1(z).
\end{aligned}
\end{equation*}

\begin{lemma}
\label{rel:2}
Define
$$
\mathcal{S}^{\star}_{n, \mathcal{L}}= \left\{z: \;  \bar\Theta_{\mathcal{L}, \mathcal{L}}\rl_1(z)<0,\; \rl_2(z)> -ca_n \bar{\Sigma}^{-1}_{\mathcal{I}, \mathcal{I}}\bar{\mu}^{(\mathcal{I})},\rl_3(z) >-ca_n( \bar{\Sigma}_{\bar{J}_{\mathcal{L}}, \mathcal{I}} \bar{\Sigma}^{-1}_{\mathcal{I}, \mathcal{I}} \bar{\mu}^{(\mathcal{I})} - \bar{\mu}^{(\bar{J}_{\mathcal{L}})})  \right\},
$$
for $c\in (0,1)$.
For all $z\in \mathcal{S}^{\star}_{n, \mathcal{L}}$, we have that
$$\mathbb{P}\left[T^{\star}_n > 0\right]  = \text{\normalfont Exp}\left(\Sigma^{1/2}z+ \sqrt{n}\beta_n -r, \frac{1}{\rho^2}\Delta\right)\cdot \mathcal{A}^{\star}_{n,\mathcal{L}}(z; \sqrt{n}\beta_n)$$
for a Lebesgue-almost everywhere differentiable function $\mathcal{A}^{\star}_n(\cdot; \sqrt{n}\beta_n)$ which satisfies
 $$\displaystyle\sup_{ \mathcal{M}_r}\sup_{z\in \mathcal{S}^{\star}_n} (a_n)^{|\mathcal{I}|} \|\cD^m\mathcal{A}^{\star}_n(z; \sqrt{n}\beta_n)\| <\infty,\ m=0,1,2,3.$$
\end{lemma}

%Before outlining the proof for Proposition \ref{rep}, we discuss some notations and make two observations in Lemma \ref{rel:1} and \ref{rel:2}.
%For two vectors $u$ and $v$ of matching dimensions, $u\odot v$ represents their componentwise product.
%For a vector $u$ with all nonzero entries, $u^{-1}$ yields the vector $z$ such that $z^{\;(j)}= 1/u^{\;(j)}$.

\begin{proof}[Proof of Lemma \ref{rel:2}.]
Note that
\begin{equation*}
\begin{aligned}
\mathbb{P}\left[T^{\star}_n > 0\right]  &\propto \int_{t> -\vartheta_z} \tE\left(t, \frac{1}{\rho^2}\bar{\Sigma}^{-1}\right) dt\\
&\propto \tE\left(\vartheta_z^{(\bar{I}_{\mathcal{L}})}, \frac{1}{\rho^2}(\bar\Sigma_{\bar{I}_{\mathcal{L}}, \bar{I}_{\mathcal{L}}})^{-1}\right) \\
&\;\;\;\;\;\;\;\;\;\;\;\times \int_{\mathcal{H}_{\vartheta_z}} \tE\left(t, \frac{1}{\rho^2}\bar{\Sigma}^{-1}\right) \cdot  \exp\left(\frac{1}{\rho^2}(t^{(\bar{I}_{\mathcal{L}})})' \bar{\Sigma}^{-1}_{\bar{I}_{\mathcal{L}}, \bar{I}_{\mathcal{L}}} \vartheta_z^{(\bar{I}_{\mathcal{L}})} \right)dt,
%&= A_0\cdot \tE\left(\vartheta_z^{(\bar{I}_{\mathcal{L}})}, \frac{1}{\rho^2}(\Sigma_{\bar{I}_{\mathcal{L}}, \bar{I}_{\mathcal{L}}})^{-1}\right) \cdot \mathrm{I}_{\vartheta_z}.
\end{aligned}
\end{equation*}
where
$$\mathcal{H}_{\vartheta_z}= \left\{t: t^{(\bar{I}_{\mathcal{L}})}>0, t^{(\bar{J}_{\mathcal{L}})}> \bar{\Sigma}_{\bar{J}_{\mathcal{L}}, \bar{I}_{\mathcal{L}}} \bar{\Sigma}^{-1}_{\bar{I}_{\mathcal{L}}, \bar{I}_{\mathcal{L}}}\vartheta_z^{(\bar{I}_{\mathcal{L}})}- \vartheta_z^{(\bar{J}_{\mathcal{L}})} \right\}.$$

For two vectors $u$ and $v$ of matching dimensions, let $u\odot v$ represent their component-wise product, and let $u^{-1}$ denote the vector $z$ where $z^{\;(j)}= \frac{1}{u^{\;(j)}}$.
Define the following functions:
$$\rH_{0,n}(z) =\prod_j \frac{1}{\rho^2}e_j'(a_n\bar{\Sigma}^{-1}_{\mathcal{I}, \mathcal{I}}\bar{\mu}^{(\mathcal{I})}+ \rl_2(z)),$$
$$\rH_{1,n}(z)= \int_{\mathcal{H}_{\vartheta_z}}  \tE \left(\mathcal{Q}^{-1}(s) , \frac{1}{\rho^2} \bar{\Sigma}^{-1}\right)\exp\left( - (s^{(\mathcal{I})})' \mathbf{1}_{|\mathcal{I}|} +(s^{(\mathcal{L})})' \Theta_{\mathcal{L}, \mathcal{L}}\rl_1(z)  \right)ds^{(\bar{I}_{\mathcal{L}})} ds^{(\bar{J}_{\mathcal{L}})}.$$
By changing variables $$t\stackrel{\mathcal{Q}}{\longrightarrow} s,$$ where
$$\mathcal{Q}^{-1}(s) = \begin{pmatrix} \left(\frac{1}{\rho^2}(a_n\bar{\Sigma}^{-1}_{\mathcal{I}, \mathcal{I}}\bar{\mu}^{(\mathcal{I})}+ \rl_2(z))^{-1}\odot  s^{(\mathcal{I})}\right)' & (s^{(\mathcal{L})})' & (s^{(\bar{J}_{\mathcal{L}})})'\end{pmatrix}',$$
we observe that
$$
\int_{\mathcal{H}_{\vartheta_z}} \tE\left(t, \frac{1}{\rho^2}\bar{\Sigma}^{-1}\right) \cdot  \exp\left(\frac{1}{\rho^2}(t^{(\bar{I}_{\mathcal{L}})})' \bar{\Sigma}^{-1}_{\bar{I}_{\mathcal{L}}, \bar{I}_{\mathcal{L}}} \vartheta_z^{(\bar{I}_{\mathcal{L}})} \right)dt = (\rH_{0,n}(z))^{-1}\cdot \rH_{1,n}(z).
$$ 
Letting $\mathrm{H}_2(z) = \tE\left(\rl_1(z), \frac{1}{\rho^2}\Theta_{\mathcal{L}, \mathcal{L}}\right)$ and simplifying some more, we note that
\begin{equation*}
\begin{aligned}
\mathbb{P}\left[T^{\star}_n > 0\right]  &\propto \tE\left(\vartheta_z^{(\bar{I}_{\mathcal{L}})}, \frac{1}{\rho^2}(\bar\Sigma_{\bar{I}_{\mathcal{L}}, \bar{I}_{\mathcal{L}}})^{-1}\right)  \cdot (\rH_{0,n}(z))^{-1}\cdot \rH_{1,n}(z)\\
&\propto \tE\left(\Sigma^{1/2}z+ \sqrt{n}\beta_n -r, \frac{1}{\rho^2}\Delta\right)  \cdot (\rH_{0,n}(z))^{-1}\cdot \rH_{1,n}(z) \cdot \mathrm{H}_2(z).
\end{aligned}
\end{equation*}
The proof is now complete by letting
 $$\mathcal{A}^{\star}_{n, \mathcal{L}}(z; \sqrt{n}\beta_n) \propto (\rH_{0,n}(z))^{-1}\cdot \rH_{1,n}(z) \cdot \mathrm{H}_2(z).$$ 
\end{proof}

\begin{proof}[Proof of Proposition \ref{rep}.]
The desired representation follows by defining
$$\mathcal{A}_n(z; \sqrt{n}\beta_n^{\;(j)}) = \sum_{\mathcal{L}\subseteq \mathcal{J}_2}\mathcal{A}^{\star}_{n, \mathcal{L}}(z) \cdot \mathbf{1}_{z\in  \mathcal{S}^{\star}_{n, \mathcal{L}}}, \text{ and } \ \mathcal{S}_n = \underset{\mathcal{L}\subseteq \mathcal{J}_2}{\cup} \mathcal{S}^{\star}_{n, \mathcal{L}}.$$
\end{proof}

\subsubsection{Proofs of main results}
\begin{proof} [Proof of Theorem 3]
Denote by $\mathcal{M}_b$ the collection of parameters that satisfy Equation (5.1) and note that
$$
 \sup_{\mathcal{M}_b} \left(\mathbb{E}_{\mathcal{N}}\left[\rF \left( \Sigma^{1/2}\mathcal{Z}_n + \sqrt{n}\beta_n\right)\right]  \right)^{-1} \geq  E_1.
 $$
 Therefore, it is sufficient to show that
\begin{equation*}
\begin{aligned}
& \lim_n\sup_{\mathbb{P}_n \in \mathcal{P}_{b,n}} \text{\normalfont SB}_{\mathbb{P}_n}(\rG_l)=0
\end{aligned}
\end{equation*}
to establish convergence of our relative differences.

%Before applying the Stein bound in Lemma 2, we 
%$$\Big|\mathbb{E}_{\mathbb{P}_n}\left[\rG_l(\mathcal{Z}_n; \sqrt{n}\beta_n)\right] - \mathbb{E}_{\mathcal{N}}\left[\rG_l(\mathcal{Z}_n; \sqrt{n}\beta_n)\right]\Big|.$$
%let 
Based on Equation 5 in \cite{chatterjee2007multivariate}, we have
\begin{equation*}
\begin{aligned}
&\cD^3\cS_{\rG_l}\left(\mathcal{Z}_n[-1]+ \frac{\alpha}{\sqrt{n}}\mathrm{e}_{1,n}+ \frac{\kappa}{\sqrt{n}}\mathrm{e}^{\star}_{1,n}\right)[ j, k, l]\\
&\;\;\;\;\;\;\;\;\;\;\;\;\;\;\;\;\;\;\;\;\;\;\;\;\;\;\;\;\;\;\;\;=\dfrac{1}{2} \displaystyle\bigintsss_{0}^{1} \sqrt{t} \cdot  \mathbb{E}_{\mathcal{N}}\left[\cD^3\rG_l \left(W_{t,\alpha, \kappa}+\sqrt{1-t}Z;\sqrt{n}\beta_n\right)[j, k, l]\right]dt,
\end{aligned}
\label{simp}
\end{equation*}
where
$$W_{t,\alpha, \kappa} = \sqrt{t}\left(\mathcal{Z}_n[-1] + \frac{\alpha}{\sqrt{n}}\mathrm{e}_{1,n}+ \frac{\kappa}{\sqrt{n}}\mathrm{e}^{\star}_{1,n}\right).$$
%Simplifying the Stein bound, we note that
The smoothness properties of our pivot in Proposition 7 lead us to note that
$$\Big|\cD^{3}\rG_l \left(\mathcal{Z}_n;\sqrt{n}\beta_n\right)[j, k, l]\Big|\leq E_{l,2}\cdot\Bigg(\sum_{\lambda, \gamma \in \{0\} \cup[3]: \lambda+\gamma\leq 3}\|\mathcal{Z}_n\|^{\lambda} \|\sqrt{n}\beta_n\|^{\gamma}\Bigg).$$
Using both facts, it follows that the Stein bound for $\rG_l$ satisfies
\begin{equation*}
\begin{aligned}
\text{\normalfont SB}_{\mathbb{P}_n}(\rG_l) &\leq \frac{E_{l,3}}{\sqrt{n}}\cdot \sum_{\lambda, \gamma, \bar\lambda, \bar\gamma \in \{0\} \cup[3]: \lambda+\gamma\leq 3; \bar\lambda + \bar\gamma \leq 3 } \mathbb{E}_{\mathbb{P}_n}\Bigg[\|\mathrm{e}_{1,n}\|^{\lambda}\|\mathrm{e}^{\star}_{1,n}\|^{\gamma} \|\sqrt{n}\beta_n\|^{\bar\gamma}\\
&\;\;\;\;\;\;\;\;\;\;\;\;\;\;\;\;\;\;\;\;\;\;\;\;\;\;\;\;\;\;\;\;\;\;\;\;\;\;\;\;\;\;\;\;\;\;\;\;\;\;\;\;\;\;\;\;\;\;\;\;\;\;\;\;\times \sup_{\alpha, \kappa \in [0,1]}\Big\|\mathcal{Z}_n[-1] + \frac{\alpha}{\sqrt{n}}\mathrm{e}_{1,n}+ \frac{\kappa}{\sqrt{n}}\mathrm{e}^{\star}_{1,n}\Big\|^{\bar\lambda} \Bigg].
%&\leq \frac{B_{l,2}}{\sqrt{n}}\cdot \sup_n\sup_{\mathbb{P}_n \in \mathcal{P}_{b,n}} \mathbb{E}_{\mathbb{P}_n}\left[ \|\mathrm{e}_{1,n}\|^5 \right].
\end{aligned}
\end{equation*}

Due to the independence between $\mathcal{Z}_n[-1]$, $\mathrm{e}_{1,n}$, and $\mathrm{e}^{\star}_{1,n}$, we can further write that
$$\text{\normalfont SB}_{\mathbb{P}_n}(\rG_l)  \leq \frac{E_{l,4}}{\sqrt{n}}\cdot \sup_n\sup_{\mathbb{P}_n \in \mathcal{P}_{b,n}} \mathbb{E}_{\mathbb{P}_n}\left[ \|\mathrm{e}_{1,n}\|^6 \right].$$
Clearly, under Assumption 4,
$$\lim_n \sup_{\mathbb{P}_n\in \mathcal{P}_{b,n}} \rR^{(l)}_n =0.$$
As a result of Proposition 3, we conclude that
$$\lim_n \sup_{\mathbb{P}_n\in \mathcal{P}_{b,n}}\Big|\widetilde{\mathbb{E}}_{\mathbb{P}_n}\left[\rH\circ\rP^{\;(j)}(\mathcal{Z}_n;  \sqrt{n}\beta_n)\right]- \widetilde{\mathbb{E}}_{\mathcal{N}}\left[\rH\circ\rP^{\;(j)}(\mathcal{Z}_n; \sqrt{n}\beta_n)\right]\Big| = 0.$$
\end{proof}

%\begin{proof}[Proof of Proposition 7]
%Following the relation in Proposition 5, the probability on the left-hand side display is equal to:
%$$\mathbb{P}[\mathcal{T} > a_n \bar{\mu}], \text{ where } \mathcal{T}\sim \mathcal{N}(0, (1+\rho^2)\cdot \bar{\Sigma}).$$
%The proof of the two parts in the claim now closely follows the proofs for Proposition 2.1 and Corollary 4.1 in \citep{hashorva2003multivariate}; we thus omit the details here.
%\end{proof}

\begin{proof} [Proof of Theorem 4]
Let $\bar{\rG}_l(\cdot; \sqrt{n}\beta_n)$ be defined according to Lemma \ref{rep}.
Because $\bar{\rG}_l(z; \sqrt{n}\beta_n)$ agrees with $\rG_l(z; \sqrt{n}\beta_n)$ on its support $\mathcal{S}_n$, it is straightforward to see that 
$$\rR^{(l)}_n \leq  [\textbf{T1}_l] + [\textbf{T2}_l],$$
where
$$[\textbf{T1}_l]:\ \Big(\mathbb{E}_{\mathcal{N}}\Big[\rF(\Sigma^{1/2}\mathcal{Z}_n + \sqrt{n}\beta_n)\Big]\Big)^{-1}\cdot \Big|\mathbb{E}_{\mathbb{P}_n}\Big[\bar{\rG}_l(\mathcal{Z}_n; \sqrt{n}\beta_n)\Big]-  \mathbb{E}_{\mathcal{N}}\Big[\bar{\rG}_l(\mathcal{Z}_n; \sqrt{n}\beta_n)\Big]\Big|,$$
\begin{equation*}
\begin{aligned}
& [\textbf{T2}_l]:\ \Big(\mathbb{E}_{\mathcal{N}}\Big[\rF(\Sigma^{1/2}\mathcal{Z}_n + \sqrt{n}\beta_n)\Big]\Big)^{-1}\cdot \Big|\mathbb{E}_{\mathbb{P}_n}\Big[\rG_l(\mathcal{Z}_n; \sqrt{n}\beta_n)\cdot \mathbf{1}_{\mathcal{Z}_n\in \mathcal{S}^c_n}\Big]\\
&\;\;\;\;\;\;\;\;\;\;\;\;\;\;\;\;\;\;\;\;\;\;\;\;\;\;\;\;\;\;\;\;\;\;\;\;\;\;\;\;\;\;\;\;\;\;\;\;\;\;\;\;\;\;\;\;\;\;\;\;\;\;\;\;\;\;\;\;\;\;\;\;\;\;\;\;\;\;\;\;\;\;\;\;\;\;\;\;\;\;\;\;\;\;\;\;\;\;\;\;\;\;\;\;\;\;\;-  \mathbb{E}_{\mathcal{N}}\Big[\rG_l(\mathcal{Z}_n; \sqrt{n}\beta_n)\cdot \mathbf{1}_{\mathcal{Z}_n\in \mathcal{S}^c_n}\Big]\Big|.
\end{aligned}
\end{equation*}

The first term in the sum is bounded as
$$[\textbf{T1}_l] \leq  \Big(\mathbb{E}_{\mathcal{N}}\Big[\rF(\Sigma^{1/2}\mathcal{Z}_n + \sqrt{n}\beta_n)\Big]\Big)^{-1}\cdot \text{\normalfont SB}_{\mathbb{P}_n}(\bar\rG_l),$$
where $\text{\normalfont SB}_{\mathbb{P}_n}(\bar\rG_l)$ is the Stein bound in Lemma 2 with $\rg = \bar\rG_l$.
Let
$$W_{t,\alpha, \kappa} = \sqrt{t}\left(\mathcal{Z}_n[-1] + \frac{\alpha}{\sqrt{n}}\mathrm{e}_{1,n}+ \frac{\kappa}{\sqrt{n}}\mathrm{e}^{\star}_{1,n}\right).$$
Using the properties of our pivot in Proposition 7, the partial derivatives of the Stein function can be seen to satisfy
\begin{equation*}
\begin{aligned}
&\Big|\cD^3\cS_{\bar{\rG}_l}\left(\mathcal{Z}_n[-1]+ \frac{\alpha}{\sqrt{n}}\mathrm{e}_{1,n}+ \frac{\kappa}{\sqrt{n}}\mathrm{e}^{\star}_{1,n}\right)[j, k, l]\Big|\\
&= \dfrac{1}{2} \displaystyle\int_{0}^{1} \sqrt{t}\cdot  \mathbb{E}_{\mathcal{N}}\Bigg[\Big|\cD^3\bar{\rG}_l\Big(W_{t,\alpha, \kappa} +\sqrt{1-t} \mathcal{Z};\sqrt{n}\beta_n\Big)[j,k, l]\Big|\cdot \mathbf{1}_{W_{t,\alpha, \kappa}  +\sqrt{1-t} \mathcal{Z}\in \mathcal{S}_n}\Bigg]dt\\
&\leq \dfrac{E_{l,5}}{(a_n)^{|\mathcal{I}|}}\cdot \sum_{\bar{\lambda}, \bar{\kappa}, \breve{\lambda}, \breve{\kappa} \in {0}\cup[3]:\bar{\lambda}+\bar{\kappa}+\breve{\lambda}+\breve{\kappa}\leq 3} \displaystyle\int_{0}^{1}\sqrt{t} \cdot \frac{1}{n^{(\bar{\lambda}+\bar{\kappa})/2}}\|\mathrm{e}_{1,n}\|^{\bar{\lambda}} \|\mathrm{e}^{\star}_{1,n}\|^{\bar{\kappa}}  \|a_n\bar{\beta}\|^{\breve{\lambda}}\|\mathcal{Z}_n[-1]\|^{\breve{\kappa}}\ \\
&\;\;\;\;\;\;\;\;\;\;\;\;\;\;\;\;\;\;\;\;\;\;\;\;\;\;\;\;\;\;\;\;\;\;\;\;\;\;\;\;\;\;\;\;\;\;\;\;\;\;\;\;\;\;\;\;\;\;\;\;\times  \tE\left(\Sigma^{1/2}W_{t,\alpha, \kappa} + \sqrt{n}\beta_n-r, \frac{1}{(1-t + \rho^2)}\cdot (\Lambda + \Delta)\right) dt.
\end{aligned}
\end{equation*}
As a result, we have
\begin{equation*}
\begin{aligned}
\text{\normalfont SB}_{\mathbb{P}_n}(\bar\rG_l) &\leq\frac{E_{l,6}}{(a_n)^{|\mathcal{I}|}\sqrt{n}}\sum_{\substack{\lambda, \gamma,  \bar{\lambda}, \bar{\kappa}, \breve{\lambda}, \breve{\kappa} \in {0}\cup[3]: \\  \lambda+\gamma\leq 3,:\bar{\lambda}+\bar{\kappa}+\breve{\lambda}+\breve{\kappa}\leq 3}} \mathbb{E}_{\mathbb{P}_n}\Bigg[\|\mathrm{e}_{1,n}\|^{\lambda}\|\mathrm{e}^{\star}_{1,n}\|^{\gamma} \\
&\;\;\;\;\;\;\;\;\;\;\;\;\;\;\;\;\;\;\;\;\;\;\;\;\;\;\;\;\;\;\;\;\;\;\;\;\;\;\;\times  \sup_{\alpha, \kappa\in [0,1]}\displaystyle\int_{0}^{1} \frac{\sqrt{t}}{n^{(\bar{\lambda}+\bar{\kappa})/2}}\|\mathrm{e}_{1,n}\|^{\bar{\lambda}} \|\mathrm{e}^{\star}_{1,n}\|^{\bar{\kappa}}  \|a_n\bar{\beta}\|^{\breve{\lambda}}\|\mathcal{Z}_n[-1]\|^{\breve{\kappa}} \\
&\;\;\;\;\;\;\;\;\;\;\;\;\;\;\;\;\;\;\;\;\;\;\;\;\;\;\;\;\;\;\;\;\;\;\;\;\;\;\;\;\times \tE\left(\Sigma^{1/2}W_{t,\alpha, \kappa} + \sqrt{n}\beta_n-r, \frac{1}{(1-t + \rho^2)}\cdot (\Lambda + \Delta)\right) dt\Bigg].
\end{aligned}
\end{equation*}
Without losing generality, we focus on the case when $\bar{\lambda}=\bar{\kappa}=\breve{\kappa}=0$, $\breve{\lambda}=3$. 
The same line of reasoning applies to other values of $\bar{\lambda}$,$\bar{\kappa}$,$\breve{\lambda}$, $\breve{\kappa}$ and is therefore omitted from the proof.
%Hereafter, we simplify a prototype term in the sum (on the right-hand side display), noting that the same line of reasoning will apply to the remaining terms in the sum.
%Fix $\lambda, \gamma \in \{0\} \cup[3]$ such that $\lambda+\gamma\leq 3$, and say, $\bar{\lambda}=\bar{\kappa}=\breve{\kappa}=0$, $\breve{\lambda}=3$, which gives us the prototype term
%\begin{equation*}
%\begin{aligned}
%[\textbf{T}_l] &:\dfrac{E_{l,5}}{(a_n)^{|\mathcal{I}|}\sqrt{n}} \cdot \|a_n\bar{\beta}\|^3 \cdot  \mathbb{E}_{\mathbb{P}_n}\Bigg[ \|\mathrm{e}_{1,n}\|^{\lambda}\|\mathrm{e}^{\star}_{1,n}\|^{\gamma} \cdot \sup_{\alpha, \kappa\in [0,1]}\displaystyle\int_{0}^{1}\sqrt{t}\\
%&\;\;\;\;\;\;\;\;\;\;\;\;\;\;\;\;\;\;\;\;\;\;\;\;\;\;\;\;\;\;\;\;\;\;\;\;\;\;\;\;\times  \tE\left(\Sigma^{1/2}W_{t,\alpha, \kappa} + \sqrt{n}\beta_n-r, \frac{1}{(1-t + \rho^2)} \cdot (\Lambda + \Delta)\right)  dt\Bigg].
%\end{aligned}
%\end{equation*}
Using Lemma \ref{bound:main:proof}, we conclude that
$$\text{\normalfont SB}_{\mathbb{P}_n}(\bar\rG_l)  \leq E_{l,6} \cdot \frac{a_n^{3-|\mathcal{I}|}}{\sqrt{n}}\cdot \tE\left(a_n\bar{\beta}, \frac{1}{(1+\rho^2)} \cdot (\Lambda + \Delta)\right),$$ 
and as a consequence,
$$[\textbf{T1}_l] \leq  E_{l,6} \cdot \frac{a_n^{3-|\mathcal{I}|}}{\sqrt{n}} \Big(\mathbb{E}_{\mathcal{N}}\Big[\rF(\Sigma^{1/2}\mathcal{Z}_n + \sqrt{n}\beta_n)\Big]\Big)^{-1} \cdot \tE\left(a_n\bar{\beta}, \frac{1}{(1+\rho^2)} \cdot (\Lambda + \Delta)\right).$$
At last, note that
\begin{equation*}
\begin{aligned}
&\mathbb{E}_{\mathcal{N}}\left[\rF \left( \Sigma^{1/2}\mathcal{Z}_n + \sqrt{n}\beta_n\right)\right] \propto\frac{1}{(a_n)^{|\mathcal{I}|}}\cdot \text{\normalfont Exp}\left(a_n\bar{\beta}, \frac{1}{(1+\rho^2)}\cdot (\Lambda+ \Delta)\right),
\end{aligned}
\end{equation*}
which is stated in Corollary 2.
Thus, we have
$$\lim_n\displaystyle\sup_{\mathbb{P}_n\in \mathcal{P}_{r,n}} [\textbf{T1}_l]=0.$$

Next, we turn to the second term in our bound and observe that
\begin{equation*}
\begin{aligned}
\displaystyle\sup_{\mathbb{P}_n\in \mathcal{P}_{r,n}} [\textbf{T2}_l] &\leq  2\cdot \displaystyle\sup_{\mathbb{P}_n\in \mathcal{P}_{r,n}}\Big(\mathbb{E}_{\mathcal{N}}\Big[\rF(\Sigma^{1/2}\mathcal{Z}_n + \sqrt{n}\beta_n)\Big]\Big)^{-1} \cdot \mathbb{E}_{\mathbb{P}_n}\Big[\rG_l(\mathcal{Z}_n; \sqrt{n}\beta_n)\cdot \mathbf{1}_{\mathcal{Z}_n\in \mathcal{S}^c_n}\Big]\\
&\Scale[0.95]{\leq 2\rK\cdot \displaystyle\sup_{\mathbb{P}_n\in \mathcal{P}_{r,n}}\Big(\mathbb{E}_{\mathcal{N}}\Big[\rF(\Sigma^{1/2}\mathcal{Z}_n + \sqrt{n}\beta_n)\Big]\Big)^{-1}\cdot \mathbb{E}_{\mathbb{P}_n}\Big[\rF(\Sigma^{1/2}\mathcal{Z}_n + \sqrt{n}\beta_n) \cdot \mathbf{1}_{\mathcal{Z}_n\in \mathcal{S}_n^c}\Big]}\\
&\Scale[0.95]{\leq 2\rK \cdot \mathrm{L}_n \cdot \displaystyle\sup_n \displaystyle\sup_{\mathbb{P}_n\in \mathcal{P}_{r,n}}\dfrac{ \mathbb{E}_{\mathbb{P}_n}\Big[\rF(\Sigma^{1/2}\mathcal{Z}_n + \sqrt{n}\beta_n)\cdot \mathbf{1}_{\mathcal{Z}_n\in \mathcal{S}_n^c}\Big]}{\mathbb{E}_{\mathcal{N}}\Big[\rF(\Sigma^{1/2}\mathcal{Z}_n + \sqrt{n}\beta_n) \cdot \mathbf{1}_{\mathcal{Z}_n\in \mathcal{S}_n^c}\Big]} .}\\
\end{aligned}
\end{equation*}
where 
$$\mathrm{L}_n = \ \dfrac{\mathbb{E}_{\mathcal{N}}\Big[\rF(\Sigma^{1/2}\mathcal{Z}_n + \sqrt{n}\beta_n) \cdot \mathbf{1}_{\mathcal{Z}_n \in \mathcal{S}_n^c}\Big]}{\mathbb{E}_{\mathcal{N}}\Big[\rF(\Sigma^{1/2}\mathcal{Z}_n + \sqrt{n}\beta_n)\Big]}.
%=  \sup_{\beta_n \in \mathcal{M}_{r.n}}\;  \dfrac{\mathbb{P}_{\mathcal{N}}\left[ E_n= \oE, \; A_n= \oA,\; \mathcal{S}_n^c\right]}{\mathbb{P}_{\mathcal{N}}\left[ E_n= \oE, \; A_n= \oA\right]}
$$
We have $\displaystyle\lim_n \mathrm{L}_n =0$, based on Lemma \ref{rep}. 
Using Assumption 7, we conclude that
$$\lim_n\displaystyle\sup_{\mathbb{P}_n\in \mathcal{P}_{r,n}} [\textbf{T2}_l]=0.$$
This proves uniform convergence of our relative differences.
As a consequence, we claim that
$$\lim_n \sup_{\mathbb{P}_n\in \mathcal{P}_{r,n}}\Big|\widetilde{\mathbb{E}}_{\mathbb{P}_n}\left[\rH\circ\rP^{\;(j)}(\mathcal{Z}_n;  \sqrt{n}\beta_n)\right]- \widetilde{\mathbb{E}}_{\mathcal{N}}\left[\rH\circ\rP^{\;(j)}(\mathcal{Z}_n; \sqrt{n}\beta_n)\right]\Big| = 0.$$
\end{proof}

%%%%%%%%%%%%%%%%%%%%%%%%%%%%%%%%%%%%%%%%%

\begin{proof}[Proof of Proposition 8]
For a fixed value $V_n$, we let 
\begin{equation*}
\pi_{V_n}(T_n, A_n) = Q T_n + \begin{pmatrix} {\Lambda^{(\oE)}}' & {A_n}'\end{pmatrix}' - \begin{pmatrix}{V_n^{\;(\oE)}}' & {V_n^{\;(\oE^c)}}' \end{pmatrix}'.
\end{equation*}

We apply the change of variables defined in \eqref{linear:map:sel} and ignore the constant Jacobian to obtain the likelihood based on the distribution of $V_n$, $T_n$ and $A_n$.
We note that this likelihood is proportional to 
$$\rp_n(V_n)\cdot \bar\rp_n( \pi_{V_n}(T_n, A_n)| V_n).$$
Since the selection outcome is equivalent to $\{T_n\in \mathbb{R}^{p+}, A_n =\oA\}$,
we note that the conditional density for $V_n$ at $v$ is equal to
$$ 
\left(\mathbb{E}_{\mathbb{P}_n}\left[\bar\rF_n \left(V_n\right)\right]\right)^{-1} \cdot \rp_n(v)\cdot \bar{\rF}_n(v).
$$
Expressed in terms of the standardized variable $\mathcal{Z}_n$, the likelihood ratio in the claim is equal to
$$\left(\mathbb{E}_{\mathbb{P}_n}\left[\bar\rF_n \left( \Sigma^{1/2}\mathcal{Z}_n + \sqrt{n}\beta_n\right)\right]\right)^{-1}\bar\rF_n \left( \Sigma^{1/2}\oZ + \sqrt{n}\beta_n\right).$$
\end{proof}

\begin{proof} [Proof of Proposition 9] 
Observe that the expected difference
$$
\Big|\overline{\mathbb{E}}_{\mathbb{P}_n}\left[\rH\circ\rP^{\;(j)}(\mathcal{Z}_n;  \sqrt{n}\beta_n)\right]- \widetilde{\mathbb{E}}_{\mathcal{N}}\left[\rH\circ\rP^{\;(j)}(\mathcal{Z}_n; \sqrt{n}\beta_n)\right]\Big| 
$$
is bounded above by 
\begin{equation*}
\begin{aligned}
&\Big|\overline{\mathbb{E}}_{\mathbb{P}_n}\left[\rH\circ\rP^{\;(j)}(\mathcal{Z}_n;  \sqrt{n}\beta_n)\right] - \widetilde{\mathbb{E}}_{\mathbb{P}_n}\left[\rH\circ\rP^{\;(j)}(\mathcal{Z}_n;  \sqrt{n}\beta_n)\right]\Big| \\
&+ \Big|\widetilde{\mathbb{E}}_{\mathbb{P}_n}\left[\rH\circ\rP^{\;(j)}(\mathcal{Z}_n;  \sqrt{n}\beta_n)\right] - \widetilde{\mathbb{E}}_{\mathcal{N}}\left[\rH\circ\rP^{\;(j)}(\mathcal{Z}_n; \sqrt{n}\beta_n)\right]\Big|. 
\end{aligned}
\end{equation*}
Because
$$\lim_n \sup_{\mathbb{P}_n\in \mathcal{C}_n} \Big|\widetilde{\mathbb{E}}_{\mathbb{P}_n}\left[\rH\circ\rP^{\;(j)}(\mathcal{Z}_n;  \sqrt{n}\beta_n)\right] - \widetilde{\mathbb{E}}_{\mathcal{N}}\left[\rH\circ\rP^{\;(j)}(\mathcal{Z}_n; \sqrt{n}\beta_n)\right]\Big| =0,$$
we are left to prove that
\begin{equation}
\label{suff:ass}
\lim_n \sup_{\mathbb{P}_n\in \mathcal{C}_n}\Big|\overline{\mathbb{E}}_{\mathbb{P}_n}\left[\rH\circ\rP^{\;(j)}(\mathcal{Z}_n;  \sqrt{n}\beta_n)\right] - \widetilde{\mathbb{E}}_{\mathbb{P}_n}\left[\rH\circ\rP^{\;(j)}(\mathcal{Z}_n;  \sqrt{n}\beta_n)\right]\Big|=0.
\end{equation}
Simplifying the expected difference in \eqref{suff:ass}, we note that
\begin{equation*}
\begin{aligned}
& \lim_n \sup_{\mathbb{P}_n\in \mathcal{C}_n}\Big|\overline{\mathbb{E}}_{\mathbb{P}_n}\left[\rH\circ\rP^{\;(j)}(\mathcal{Z}_n;  \sqrt{n}\beta_n)\right] - \widetilde{\mathbb{E}}_{\mathbb{P}_n}\left[\rH\circ\rP^{\;(j)}(\mathcal{Z}_n;  \sqrt{n}\beta_n)\right]\Big|\\
&\leq 2\sup_z  |\rH(z)| \cdot  \lim_n \sup_{\mathbb{P}_n\in \mathcal{C}_n} \int \Big| \dfrac{\bar\rF_n(z;  \sqrt{n}\beta_n)}{\mathbb{E}_{\mathbb{P}_n}\left[\bar\rF_n(\mathcal{Z}_n;  \sqrt{n}\beta_n)\right]}  - \dfrac{\rF(z;  \sqrt{n}\beta_n)}{\mathbb{E}_{\mathbb{P}_n}\left[\rF(\mathcal{Z}_n;  \sqrt{n}\beta_n)\right]}\Big| d\mathbb{P}_n(z)\\
& \leq 2 \rK\cdot \lim_n \sup_{\mathbb{P}_n\in \mathcal{C}_n}\dfrac{\mathbb{E}_{\mathbb{P}_n}\left[|\bar\rF_n(\Sigma^{1/2}\mathcal{Z}_n +\sqrt{n}\beta_n) - \rF(\Sigma^{1/2}\mathcal{Z}_n +\sqrt{n}\beta_n)|\right]}{\mathbb{E}_{\mathbb{P}_n}\left[ \rF(\Sigma^{1/2}\mathcal{Z}_n +\sqrt{n}\beta_n)\right]}.
\end{aligned}
\end{equation*}
Our assertation holds whenever the limit in the final display is equal to $0$.
\end{proof}

\subsection{Proofs for Section 6}
\label{supp:theory}

\begin{proof}[Proof of Proposition 10]
The proof proceeds along the exact same lines as Proposition 1 by using the specific mapping:
\begin{equation*}
\begin{aligned}
\begin{pmatrix}{W_n^{\;(\oE)}}'  &  {W_n^{\;(\oE^c)}}' \end{pmatrix}'  &= Q T_n + \begin{pmatrix} A_{1,n}' & A_{2,n}'\end{pmatrix}' - \begin{pmatrix}{V_n^{\;(\oE)}}' & {V_n^{\;(\oE^c)}}' \end{pmatrix}'\\
&= \pi_{V_n}(T_n, A_{2,n})
\end{aligned}
\end{equation*}
in place of \eqref{linear:map:sel}.
\end{proof}

\begin{proof}[Proof of Proposition 11]
The proof is direct by substituting \eqref{linear:map:sel} in the proof of Proposition 1 with the specific mapping
\begin{equation}
\begin{aligned}
\begin{pmatrix}{W_n^{\;(\oE)}}'  &  {W_n^{\;(\oE^c)}}' \end{pmatrix}'  &= Q T_n + \begin{pmatrix} A_{1,n}' & A_{2,n}'\end{pmatrix}'  - P\begin{pmatrix}{V_n^{\;(\oE)}}' & {V_n^{\;(\oE^c)}}' \end{pmatrix}'\\
&= \pi_{PV_n}(T_n, A_{2,n}).
\end{aligned}
\label{cov:reg}
\end{equation}
\end{proof}

\end{document}